\documentclass [11pt] {article}
\usepackage{amssymb}
\usepackage{latexsym}
\usepackage{amsmath}
\usepackage{amscd}
\title{Phase Space Reduction of\\
Star Products on Cotangent Bundles}
\author{{\bf Niels Kowalzig
\!\!\thanks{nkowalzi@science.uva.nl}}\\[3mm]
Korteweg-de Vries Instituut voor Wiskunde\\Universiteit van
Amsterdam\\Plantage Muidergracht 24\\ NL-1018 TV Amsterdam\\
The Netherlands\\[3mm] {\bf Nikolai Neumaier
\!\!\thanks{neumaier@math.uni-frankfurt.de}\,,  Markus J.\ Pflaum
\!\!\thanks{pflaum@math.uni-frankfurt.de}}\\[3mm]
Fachbereich Mathematik\\ Universit\"at Frankfurt
\\Robert-Mayer-Stra\ss e 10\\D-60054 Frankfurt a.~M.\\
Germany\\[3mm]}
\date{March 2004}
\newcommand {\Mind} {{\mbox{\tiny $M$}}}
\newcommand {\Qind} {{\mbox{\tiny $Q$}}}
\newcommand {\QKotind} {{\mbox{\tiny $T^*Q$}}}
\newcommand {\tr} [1]{\mathsf{tr} \left({#1}\right)}
\newcommand {\cc} [1] {\overline {{#1}}}
\newcommand {\Cinf} [1] {\mathcal C^\infty ({#1})}
\newcommand {\Ginf} [1] {\Gamma^\infty ({#1})}
\newcommand {\pair} [2] {\langle {#1}, {#2}\rangle}
\newcommand {\JN} [1] {J_0({#1})}
\newcommand {\JP} [1] {J_+({#1})}
\newcommand {\Jbold} [1] {J({#1})}
\newcommand {\jN} [1] {j_0({#1})}

\newcommand {\jbold} [1] {j({#1})}
\newcommand {\im} {{\mathrm{i}}}
\newcommand {\id} {{\mathrm{id}}}
\newcommand {\Lie} {{\mathcal L}}
\newcommand {\ad} {{\mathrm{ad}}}
\newcommand {\Ad} {{\mathrm{Ad}}}
\newcommand {\red} {{\mbox{\tiny red}}}
\newcommand {\C} {{\mathrm C}}
\newcommand {\Pa} {{\mathrm P}}
\newcommand {\deltaH} {\delta_{\mbox{\rm \tiny H}}}
\newcommand {\Gind} {{\scriptscriptstyle G}}
\newcommand {\dR} {{\mbox{\rm\tiny dR}}}
\renewcommand {\d} {\mathrm{d}}

\newcommand{\PolFun} {\mathcal P (Q)}
\newcommand{\PolFunRed} {\mathcal P (\cc{Q})}
\newcommand{\repNu} [1] {\varrho_0\left({#1}\right)}

\newcommand{\repkap} [1] {\varrho_\kappa\left({#1}\right)}
\newcommand{\Fop} [1] {\mathsf F \left({#1}\right)}
\newcommand{\Pol} [1] {\mathsf P \left({#1}\right)}
\newcommand{\PolRed} [1] {\cc{\mathsf P} \left({#1}\right)}
\newcommand{\sKov} {\mathsf D}
\newcommand{\Hor} [1] {\mathsf{hor}_{#1}}
\newcommand{\Ver} {\mathsf{ver}}
\newcommand{\hor} {\mathrm{h}}
\newcommand{\starNu} {\mathbin{\star_0}}
\newcommand{\starWe}{\mathbin{\star_{1/2}}}
\newcommand{\starEi}{\mathbin{\star_1}}
\newcommand{\stark} {\mathbin{\star_\kappa}}
\newcommand{\starNuB} {\mathbin{\star^B_0}}

\newcommand{\starEiB}{\mathbin{\star^B_1}}
\newcommand{\starkB} {\mathbin{\star^B_\kappa}}
\newcommand{\starF} {\mathbin{\star_{{\mbox{\rm \tiny F}}}}}
\newcommand{\starkloc} [1] {\mathbin{\star^{#1}_\kappa}}
\newcommand{\starred} [1] {\mathbin{\star^{#1}}}
\newcommand{\bullred} [1] {\mathbin{\bullet^{#1}}}
\newcommand{\Nk} {{N_\kappa}}
\newenvironment {PROOF}{{\sc Proof:}}{{\hspace*{\fill}
                       $\square$}}
\newenvironment {INNERPROOF}{{\sc Proof:}}{{\hspace*{\fill}
                            $\bigtriangledown$}}
\renewcommand {\labelenumi} {{\it\roman{enumi}.)}}
\newtheorem {LEMMA} {Lemma} [section]
\newtheorem {SUBLEMMA} [LEMMA] {Sublemma}
\newtheorem {PROPOSITION} [LEMMA] {Proposition}
\newtheorem {THEOREM} [LEMMA] {Theorem}
\newtheorem {COROLLARY} [LEMMA] {Corollary}
\newtheorem {DEFINITION}[LEMMA] {Definition}
\newtheorem {REMARK}[LEMMA] {Remark}

\numberwithin{equation}{section}
\begin{document}
\maketitle
\begin{abstract}
\noindent
In this paper we construct star products on Marsden-Weinstein
reduced spaces in case both the original phase space and the
reduced phase space are (symplectomorphic to) cotangent bundles.
Under the assumption that the original cotangent bundle $T^*Q$
carries a symplectique structure of form $\omega_{B_0}=\omega_0 +
\pi^*B_0$ with $B_0$ a closed two-form on $Q$, is equipped by the
cotangent lift of a proper and free Lie group action on $Q$ and by
an invariant star product that admits a $G$-equivariant quantum
momentum map, we show that the reduced phase space inherits from
$T^*Q$ a star product. Moreover, we provide a concrete description
of the resulting star product in terms of the initial star product
on $T^*Q$ and prove that our reduction scheme is independent of the
characteristic class of the initial star product. Unlike other
existing reduction schemes we are thus able to reduce not only
strongly invariant star products. Furthermore in this article, we
establish a relation between the characteristic class of the
original star product and the characteristic class of the reduced
star product and provide a classification up to $G$-equivalence of
those star products on $(T^*Q,\omega_{B_0})$, which are invariant
with respect to a lifted Lie group action. Finally, we investigate
the question under which circumstances `quantization commutes with
reduction' and show that in our examples non-trivial restrictions
arise.
\end{abstract}
\clearpage
\tableofcontents
\section{Introduction}
\label{IntroSec}
Already in the first and fundamental article on deformation
quantization by Bayen et al.\ \cite{BayFla78},
the problem how to construct a star product on a reduced phase space
out of a known star product on the initial phase space has been
considered.
In particular, the example of the cotangent bundle $T^*S^{n-1}$ of the
$n-1$-sphere has been discussed, and for reduction
the Weyl-Moyal product on $T^* (\mathbb R^n\setminus \{0\})$ has been used.
Even after general existence proofs for deformation
quantizations on symplectic and Poisson manifolds have meanwhile appeared,
it remains an interesting question which relations one can establish
between star products on the original phase space and those on the
reduced phase space. In physics terms, this corresponds to the
question, whether extrinsic and intrinsic quantization are
equivalent, which in some sense would mean that `quantization
commutes with reduction' (cf.\ \cite{Got86} and \cite{Emm93a} for a
discussion of these topics in the framework of geometric
quantization and the conventional Hilbert space approach to quantum
mechanics). Moreover, the question of existence of star products on
symplectic stratified spaces is still unsolved, and it appears to be very
promising to attack this problem first for singular reduced phase spaces,
which have been studied in detail by Sjamaar and Lerman
\cite{SjaLer91}.

Considering particular examples, there are various explicit constructions
for phase spaces with additional structure, for instance
$\mathbb C P^n$ \cite{BorBriEmm96,Wal98} or more general complex
Grassmann manifolds \cite{Sch97} for which
a deformation quantization analogue of phase space reduction can be
constructed.
But these examples all seem to be well-tailored to special situations and
do not apply to arbitrary symplectic phase spaces.

In \cite{Fed98}, Fedosov introdcued a reduction scheme for
arbitrary symplectic manifolds $(M,\omega)$ with a compact, free
and symplectic Lie group action and a regular value of the momentum
map. Fedosov starts with a certain Fedosov star product obtained
from a certain $G$-invariant torsion free symplectic connection.
Adapting the original star product appropriately to the fibering
structure of the principal $G$-bundle $M_0 \to M_\red
= M_0 / G$, where $M_0$ denotes the inverse image of $0\in \mathfrak g^*$
with respect to the classical momentum map, he showed that one can
always achieve that the resulting star product on the reduced phase
space $(M_\red,\omega_\red)$ is equivalent to a canonical Fedosov
star product on it. Thus, Fedosov is able to prove a `reduction
commutes with quantization' theorem within his particular
situation. Another very general approach to reduction in the
framework of deformation quantization is the BRST-method as
presented by Bordemann et al.\ in \cite{BorHerWal00}. Using a
quantum BRST complex, the authors of this work are able to produce
quite an explicit formula for the reduced star product under the
following three assumptions:
\renewcommand {\labelenumi} {{\arabic{enumi}.)}}
\begin{enumerate}
\item the symmetry group acts properly and freely,
\item $0$ is a regular value of the momentum map,
\item the inital star product on the phase space with symmetry is
 strongly invariant
 (cf.\ Section \ref{InvDefQuaSubSec} for a definition).
\end{enumerate}
\renewcommand {\labelenumi} {{\it\roman{enumi}.)}}
By results obtained in \cite{GutRaw03} and \cite{MN03a}, which
reveal some obstructions in the characteristic class for a star product
to be strongly invariant, the last of these assumptions imposes a
restriction on the possible characteristic class of the original star product.
Note that implicitly, the same restriction appears in the Fedosov reduction
scheme.

The scope of the present paper is to develop a reduction scheme for
star products on cotangent bundles with respect to a symplectic
form which is the sum of the canonical symplectic form and the
pull-back of a closed two-form on the base manifold which can be
interpreted as a magnetic field. Additionally, we assume
that the reduced phase space is again a cotangent bundle or,
more precisely, symplectomorphic to a
cotangent bundle via a non-canonical diffeomorphism.
It is known that this latter assumption holds true (cf.\
\cite{Got86,Kum81,Mon83}) if the action is the
cotangent lift of a proper and free Lie group action on
the base manifold and if the momentum value for which the reduced
space is considered is an invariant element of the dual of the Lie algebra.
Our construction is adapted to the particular geometry of a cotangent bundle
and we have to restrict our reduction scheme to a certain class of star
products $\star$ namely those for which the space of formal functions which are
polynomials in the momenta form a $\star$-subalgebra.
Fortunately, this class of
star products is rich enough to obtain star products of arbitrary
characteristic class, hence we actually provide a reduction scheme which
does not depend on the characteristic class of the initial star product.

Our paper is organized as follows: In Section \ref{PrelimSec} we
recall from \cite{BorNeuWal98,BorNeuWal99,BorNeuPflWal03} the
construction of various star products on cotangent bundles which
have the common property that the formal functions polynomial in
the momenta form a subalgebra. We also collect some notions of
invariance with respect to Lie group actions in deformation
quantization and recall the definition of the deformation
quantization analogue of a $G$-equivariant classical momentum map.
The resulting quantum momentum maps  will play a fundamental role
in our framework of phase space reduction. In Section \ref{RedSec}
we then address the reduction of a certain class of star products
on $(T^*Q,\omega_{B_0})$ which includes the examples considered in
Section \ref{PrelimSec}. Under the assumptions imposed on the
inital data we establish a relation between the Poisson algebra of
functions on $T^*(Q/G)$ polynomial in the momenta and the Poisson
algebra of horizontal invariant functions on $T^*Q$ which are
polynomial in the momenta as well. Moreover, we succeed to
construct a deformation of this classical correspondence which
enables us to define an associative product on the polynomial
functions on $T^*(Q/G)$ which is induced by a star product on
$(T^*Q,\omega_{B_0})$. It turns out that this product can be
uniquely extended to the whole space $\Cinf{T^*(Q/G)}[[\nu]]$. We
thus obtain a star product with respect to some symplectic form on
the reduced phase space which differs from the canonical symplectic
structure  by an additional magnetic field term and which depends
on $B_0$, the chosen classical momentum map, the curvature of the
chosen connection, and the momentum value used for the phase space
reduction. Also in this section we investigate the behaviour of our
reduction scheme with respect to natural operations on star product
algebras like isomorphisms, automorphisms, and derivations and we
give conditions on which these transfer to the reduced star
products. In Section \ref{InvQMMSec} we relax our assumptions
somewhat to arbitrary Lie group actions on $Q$ and return to
consider the examples of Section \ref{PrelimSec}. We derive first
necessary and sufficient conditions on the geometric data which
guarantee the considered star products to be invariant with respect
to the lifted Lie group action on the base manifold. Furthermore,
we find additional conditions which guarantee that these products
admit a $G$-equivariant quantum Hamiltonian and thus a
$G$-equivariant quantum momentum map in the sense of Xu
\cite{Xu98}. In these particular cases it turns out that the
quantum momentum maps are polynomials in the momenta, a result
which is important for our reduction scheme to work. In addition,
if there is a $G$-invariant torsion free connection on $Q$, we can
give a classification of star products on $(T^*Q, \omega_{B_0})$
invariant with respect to a lifted group action up to
$G$-equivalence. In Section \ref{CharClassSec} we succeed in
computing the characteristic class of the reduced star products.
This clarifies how the choices made in the course of our reduction
scheme affect the equivalence class of the resulting reduced star
product. Moreover, we compare the resulting star products to
naturally given star products on the reduced cotangent space in
Section \ref{AppExaSec} and derive conditions for `reduction
commutes with quantization' results to hold. We conclude Section
\ref{AppExaSec} by comparing our construction to known results and
examples considered in the literature
\cite{BayFla78,BorHerWal00,Fed98}.
\paragraph{Acknowledgements}
N.\ N.\ is most indebted to M.\ Bordemann for pointing out
reference \cite{Bor04} and for many clarifying discussions. N.\ N.\
and M.\ P.\ gratefully acknowledge financial support by the
Deutsche Forschungsgemeinschaft.
\section{Preliminaries and Notation}
\label{PrelimSec}
In this section we first recall some notation and several canonical
constructions of star products  on cotangent bundles (see
\cite{BorNeuWal98,BorNeuWal99,BorNeuPflWal03} for further details).
Then we collect various notions of invariance in deformation
quantization with respect to Lie group actions and recall the
definition of the quantum analogue of a $G$-equivariant classical
momentum map (cf.\ \cite{ArnCor83,BerBieGut98,Xu98}). In our
framework, this notion will turn out to be fundamental for the
formulation of phase space reduction.
\subsection{Constructions of Star Products on Cotangent Bundles}
\label{StarProdConsSubSec}
Throughout this article, $Q$ will always denote a smooth
$n$-dimensional manifold. Recall that the cotangent bundle
$\pi : T^*Q \to Q$ is equipped with the canonical symplectic
form $\omega_0 = - \d \theta_0$, where $\theta_0$ denotes the canonical
one-form. The zero section of $T^*Q$ is denoted by $i: Q \to T^*Q$ by means of
which we consider $Q$ as embedded into $T^*Q$.
Local coordinates on $Q$ will be denoted by
$x^1, \ldots, x^n$, the induced coordinates on $T^*Q$ by
$q^1, \ldots, q^n, p_1, \ldots, p_n$.

Given $k$ one-forms $\beta_1,\ldots,\beta_k\in \Ginf{T^*Q}$ we define
a fiberwise acting differential operator $\Fop{\beta_1\vee
\ldots\vee\beta_k} :\Cinf{T^*Q} \rightarrow \Cinf{T^*Q}$ of order
$k$ by
\begin{equation}\label{FOpDefEq}
\left(\Fop{\beta_1\vee \ldots\vee\beta_k} f\right)(\zeta_x)
:= \left.\frac{\partial^k}{\partial t_1 \cdots \partial t_k}
\right|_{t_1=\ldots = t_k =0} f (\zeta_x + t_1 \beta_1 (x)+
\ldots + t_k \beta_k (x)), \quad \zeta_x \in T_x^*Q.
\end{equation}
Clearly,
$\mathsf F$ extends to an injective algebra morphism from
$\Ginf{\bigvee T^*Q}$ into the algebra of differential operators on
$\Cinf{T^*Q}$.

To every contravariant symmetric tensor field $T \in
\Ginf{\bigvee^kT^*Q}$ one can assign a smooth function $\Pol{T}\in
\Cinf{T^*Q}$ by
$(\Pol{T})(\zeta_x):= \frac{1}{k!} i_s(\zeta_x)\ldots
i_s(\zeta_x)T(x)$.
Note that $\Pol{T}$ is a homogeneous polynomial of degree $k$ in the
momenta. Let us denote the space of these fiberwise polynomial functions
by $\mathcal P^k(Q)$. Obviously, $\mathsf P$ then
extends to an isomorphism between the $\mathbb Z$-graded commutative
algebras $\Ginf{\bigvee TQ}$ and
$\PolFun = \bigoplus_{k=0}^\infty \mathcal P^k(Q)$.

Fixing a torsion free connection $\nabla$ on the base manifold $Q$
one can assign to every formal function $f\in \Cinf{T^*Q}[[\nu]]$ a
formal series with values in the differential operators on
$\Cinf{Q}$ by
\begin{equation}\label{RhoNullDefEq}
\repNu{f}\chi : = i^* \Fop{\exp (- \nu\sKov) \chi} f =
\sum_{l=0}^\infty \frac{(-\nu)^l}{l!} i^*\left(
\frac{\partial^l f}{\partial p_{j_1}\cdots \partial p_{j_l}}
\right) i_s(\partial_{x^{j_1}})\ldots i_s(\partial_{x^{j_l}})
\frac{1}{l!} \sKov^l \chi.
\end{equation}
Here, $\sKov$ denotes the operator of symmetric covariant derivation
which in local coordinates is given by $\sKov = \d x^i \vee
\nabla_{\partial_{x^i}}$.

It has been shown in \cite{BorNeuWal98} that the restriction of
$\varrho_0$ to $\PolFun[[\nu]]$ is injective and that the image of
$\PolFun[[\nu]]$ under $\varrho_0$ is closed with respect to
composition of differential operators. Hence, one can define an
associative product on $\PolFun[[\nu]]$ by
\begin{equation}\label{StarNullPolDefEq}
F \starNu F' := {\varrho_0}^{-1}\left( \repNu{F}\repNu{F'}
\right), \quad \text{ $F,F' \in\PolFun[[\nu]]$}.
\end{equation}
Moreover, it has been shown that $\starNu$ can be described by
bidifferential operators. Hence, one can uniquely extend $\starNu$
to a product on $\Cinf{T^*Q}[[\nu]]$ yielding a star product on
$(T^*Q,\omega_0)$ which also will be denoted by $\starNu$ and which
will be called the standard ordered star product (corresponding to
$\nabla$). By definition of $\starNu$ it is obvious that
$\varrho_0$ defines a representation of
$(\Cinf{T^*Q}[[\nu]],\starNu)$ on $\Cinf{Q}[[\nu]]$.

Now we want to define further star products depending on a
so-called order parameter $\kappa \in[0,1]$.
To this end consider a smooth positive density
$\upsilon$ on $Q$. Then $\upsilon$ and
the connection $\nabla$ define a one-form $\alpha \in \Ginf{T^*Q}$ by
\begin{equation}\label{alphaDefEq}
\nabla_X \upsilon = \alpha (X) \upsilon, \quad X \in \Ginf{TQ}.
\end{equation}
It is immediate to check that this one-form satisfies
\begin{equation}\label{dalphaEq}
\d \alpha = - \tr{R},
\end{equation}
where $R$ denotes the curvature tensor of $\nabla$ and $\tr{R}$ the
trace of the curvature endomorphism. Let us denote by
$\Hor{\nabla}$ the horizontal lift (with respect to $\nabla$) of
vector fields on $Q$ to vector fields on $T^*Q$ (cf.\ also
Definition \ref{HorVerDef}). Using local coordinates we then define
a differential operator $\Delta$ on $\Cinf{T^*Q}$ by
\begin{equation}\label{DeltaDefEq}
\Delta := \Delta_0 + \Fop{\alpha} := \Fop{\d x^i}
\Lie_{\Hor{\nabla}(\partial_{x^i})} + \Fop{\alpha} =
\frac{\partial^2}{\partial q^k\partial p_k} + p_l
\pi^*\Gamma^l_{jk}\frac{\partial^2}{\partial p_j\partial p_k}+
\pi^*(\Gamma^l_{lk}+ \alpha_k)\frac{\partial}{\partial p_k}.
\end{equation}
Here, $\Gamma^l_{jk}$ denote the Christoffel symbols of $\nabla$
and $\alpha_k$ the components of $\alpha$ in the chosen local chart.
After some immediate computation it turns out that $\Delta$ is independent
of the choice of local coordinates.
In \cite{BorNeuPflWal03} we also considered the formal series $\Nk$
of differential operators given by
\begin{equation}\label{NkappaDefEq}
\Nk := \exp (-\kappa\nu \Delta).
\end{equation}
Then this operator induces the $\kappa$-ordered star product
\begin{equation}\label{starkappaDefEq}
f \stark f' := \Nk^{\!\!\!-1}((\Nk f)\starNu (\Nk f'))
\end{equation}
which is obtained from $\starNu$ by the equivalence transformation
$\Nk^{\!\!\!-1}$. In the cases $\kappa =1$ and $\kappa
= 1/2$ the corresponding star products are called star product of
anti-standard ordered type and Weyl ordered star product. For a
further discussion of these star products we again refer the reader
to \cite{BorNeuWal98,BorNeuWal99,BorNeuPflWal03}, where one can
particularly find a motivation for the above definitions by means
of the applicability of the GNS construction to $\starWe$. For
later use let us recall the following factorization property of
$\Nk$ from \cite[Lemma 3.6]{BorNeuPflWal03}:
\begin{equation}\label{NkFactEq}
\Nk = \exp (- \kappa\nu \Delta_0) \exp\left(
-\Fop{\frac{\exp(\kappa\nu\sKov)-\id}{\sKov}\alpha}
\right).
\end{equation}
At this point let us mention two further important formulas which
explicitly determine the $\stark$-left- and the
$\stark$-right-multiplication with a formal function pulled-back
from $Q$ (cf.\ \cite[Prop.\ 3.2]{BorNeuPflWal03}):
\begin{eqnarray}\label{starklinksmultEq}
\pi^* \chi \stark f   &=& \Fop{\exp(\kappa\nu\sKov)\chi}f,\\
\label{starkrechtsmultEq} f \stark \pi^*\chi &=&
\Fop{\exp((\kappa-1)\nu\sKov)\chi}f,
\end{eqnarray}
where $\chi \in \Cinf{Q}[[\nu]]$ and $f\in \Cinf{T^*Q}[[\nu]]$.
Let us also note that by
$\repkap{f}:=
\repNu{\Nk f}$ one obtains a representation of
$(\Cinf{T^*Q}[[\nu]],\stark)$ on $\Cinf{Q}[[\nu]]$. Due to the
properties of $\varrho_0$ and $\Nk$ the restriction of
$\varrho_\kappa$ to $\PolFun[[\nu]]$ is injective as well.

Up to now we have only considered so-called homogeneous star
products on $(T^*Q,\omega_0)$, i.e.\ star products for which
$\mathcal H = \Lie_{\xi_0} + \nu \partial_\nu$ is a derivation.
Here ${\xi_0} \in \Ginf{T (T^*Q)}$ denotes the canonical Liouville
vector field defined by $i_{\xi_0} \omega_0 = - \theta_0$. From
this property and the fact that the bidifferential operator
describing $\starWe$ at order $2$ in the formal parameter is
symmetric it is easy to deduce (cf.\ \cite[Thm.\
4.6]{BorNeuPflWal03}) that the characteristic class of the star
products $\stark$ is given by $c(\stark) =c(\starWe)=[0]$.

Now we want to recall a construction which yields star products with
characteristic class $\frac{1}{\nu}
[\pi^*B]$, where $B\in Z^2_\dR(Q)[[\nu]]$ is an arbitrary formal
series of closed two-forms on $Q$. The thus obtained star products
comprise deformation quantizations of the symplectic manifold
$(T^*Q,\omega_{B_0}=\omega_0+\pi^*B_0)$, where $B_0$ -- which is
assumed to be real -- denotes the term of zeroth order in the
formal parameter of $B$.

For $A \in \Ginf{T^*Q}[[\nu]]$ with real $A_0$ and all
$\kappa\in[0,1]$ consider the operator $\mathcal A_k :
\Cinf{T^*Q}[[\nu]]\to \Cinf{T^*Q}[[\nu]]$ defined by
\begin{equation}\label{AkappaDefEq}
\mathcal A_\kappa := t^*_{-A_0} \exp\left(-
\Fop{\frac{\exp(\kappa\nu\sKov)-
\exp((\kappa-1)\nu\sKov)}{\nu\sKov}A - A_0}
\right),
\end{equation}
where $t_{-A_0}$ denotes the fiberwise translation by $-A_0$ that is
$t_{-A_0}(\zeta_x) = \zeta_x - A_0(x)$.
The important property of $\mathcal A_k$ is
that it defines an automorphism of $\stark$ if and only if $\d A
=0$ (cf.\ \cite[Thm.\ 3.4]{BorNeuPflWal03}). To define star products
reflecting the presence of a magnetic
field $B$ on $Q$ consider a good open cover $\{O_j\}_{j\in I}$
of $Q$ together with local formal potentials $A^j \in
\Ginf{T^*O_j}[[\nu]]$ of $B$. This means $B|_{O_j} = \d A^j$, where in
addition $A^j_0$ is chosen to be real. Let us denote by $\mathcal A^j_\kappa$
the operator determined by Eq.\ (\ref{AkappaDefEq}) with $A^j$ used
instead of $A$. Then one defines an
associative product $\starkloc{j}$ on $\Cinf{T^*O_j}[[\nu]]$ by
\begin{equation}\label{starBlocDefEq}
f \starkloc{j} f' := \mathcal A^j_\kappa \left( ((\mathcal
A^j_\kappa)^{-1} f) \stark ((\mathcal A^j_\kappa)^{-1}
f')\right),\quad\textrm{ }f,f' \in \Cinf{T^*O_j}[[\nu]].
\end{equation}
Now one makes the crucial observation that the operator $(\mathcal
A^k_\kappa)^{-1}\mathcal A^j_\kappa$ corresponds via
(\ref{AkappaDefEq}) to the closed formal one-form $A^j|_{O_j\cap
O_k} - A^k|_{O_j\cap O_k}$. Hence, this operator is an automorphism
of $(\Cinf{T^*(O_j\cap O_k)}[[\nu]],\stark)$ and one can define a
star product $\starkB$ on $\Cinf{T^*Q}[[\nu]]$ by setting
\begin{equation}\label{starkappaBDefEq}
\left.f \starkB f'\right|_{T^*O_j} = \left.f\right|_{T^*O_j}
\starkloc{j} \left.f' \right|_{T^*O_j}, \quad f,f'\in\Cinf{T^*Q}[[\nu]].
\end{equation}
These star products
do not depend on the particular choice of the covering nor on the
choice of the potentials but only on $B$ (and of course on
$\stark$). By  a straightforward computation one checks that
$\starkB$ is a star product on $(T^*Q, \omega_{B_0})$ for all
$\kappa\in [0,1]$. Moreover, in \cite[Thm.\ 4.6]{BorNeuPflWal03} it
has been shown that the characteristic class of $\starkB$ is given by
$c(\starkB) = \frac{1}{\nu} [\pi^*B]$. Thus there exists for every equivalence
class of star products on $(T^*Q, \omega_{B_0})$ a
representative $\starkB$.
\subsection{Notions of Invariance in Deformation Quantization}
\label{InvDefQuaSubSec}
Let $G$ be a Lie group and denote by $\mathfrak g = \mathrm{Lie}(G)$
its Lie algebra. In addition, assume  $\varphi : G
\times M \to M$ to be a (left) action of $G$ on a
manifold $M$ equipped with a symplectic form $\omega$, and denote for every
$g\in G$  by $\varphi_g : M \to M$ the map $m \mapsto \varphi(g,m)$.
Then a star product $\star$ on $(M,\omega)$ is called invariant with respect to
$\varphi : G \times M \to M$ or, for short, $G$-invariant
in case no confusion can arise, if every
$\varphi_g^*$ is an automorphism of the star product, i.e.\ if
\begin{equation}\label{GenGInvDefEq}
\varphi_g^*(f \star f') =
\varphi_g^*f \star \varphi_g^*f' \quad
\text{ for all $f,f' \in \Cinf{M}[[\nu]]$, $g\in G$.}
\end{equation}
In other words this means that $r$ defined by $r(g)f :=
\varphi_{g^{-1}}^*f$ defines a left action on
$(\Cinf{M}[[\nu]],\star)$ by automorphisms. By anti-symmetrization
of Eq.\ (\ref{GenGInvDefEq}) with respect to $f$ and $f'$ one
checks that the action $\varphi$ necessarily has to be symplectic.

Given a $G$-invariant star product $\star$ on a symplectic manifold
$(M,\omega)$, differentiation of $r$ yields an action $\rho$ of the Lie
algebra $\mathfrak g=\mathrm{Lie}(G)$ on $\Cinf{M}[[\nu]]$ by
derivations of $\star$. Explicitly,
\begin{equation}
\rho(\xi) f = - \Lie_{\xi_\Mind}f \quad \text{ for all $f\in \Cinf{M}[[\nu]]$,
  $\xi\in \mathfrak g$},
\end{equation}
where $\xi_\Mind$ denotes the fundamental vector field associated to $\xi$.
Clearly, in the $G$-invariant case, every fundamental vector field
$\xi_\Mind$ is symplectic.

Now let $\star$ be a $G$-invariant star product and denote by
$C^1(\mathfrak g,\Cinf{M})$ the space of linear forms on $\mathfrak g$
with values in $\Cinf{M}$. Then an element
$J = J_0 + J_+\in C^1(\mathfrak g,\Cinf{M})[[\nu]]$ with real-valued
$J_0\in C^1(\mathfrak g,\Cinf{M})$ and
$J_+\in \nu C^1(\mathfrak g,\Cinf{M})[[\nu]]$ is
called a quantum Hamiltonian for $r$, if
\begin{equation}\label{QHamDefEq}
\rho(\xi) = -\Lie_{\xi_\Mind} =
\frac{1}{\nu}\ad_\star(\Jbold{\xi})\quad
\text{ for all $\xi \in \mathfrak g$}.
\end{equation}
In other words this means that the Lie derivative with respect to the
generating vector fields is a quasi-inner (or essentially inner) derivation of
$\star$.
$J$ is called a $G$-equivariant quantum Hamiltonian, if it is a
quantum Hamiltonian and additionally satisfies
\begin{equation}\label{GEquQHamDefEq}
\varphi_g^* \Jbold{\xi} = \Jbold{\Ad(g^{-1})\xi}\quad
\text{for all $\xi \in \mathfrak g$, $g \in G$}.
\end{equation}
A quantum Hamiltonian $J$ is called a quantum momentum map if in
addition
\begin{equation}\label{QMMapDefEq}
\frac{1}{\nu}(\Jbold{\xi}\star\Jbold{\eta}- \Jbold{\eta}\star
\Jbold{\xi}) = \Jbold{[\xi,\eta]}
\quad \text{ for all $\xi,\eta\in \mathfrak g$}.
\end{equation}

Clearly, differentiation of (\ref{GEquQHamDefEq}) with respect to
$g$ at $e$ shows that a $G$-equivariant quantum Hamiltonian always
defines a quantum momentum map. Note that the converse generally
does not hold true. In the sequel we will refer to a
$G$-equivariant quantum Hamiltonian as a $G$-equivariant quantum
momentum map.

The zeroth order parts of (\ref{QHamDefEq}) and
(\ref{GEquQHamDefEq}) mean that $J_0$ is a $G$-equivariant
classical momentum map for $r$, i.e.\ $\xi_\Mind$ is a Hamiltonian
vector field with Hamiltonian function $\JN{\xi}$ and the smooth
mapping $\Check J_0 : M \to \mathfrak g^*$ defined by $\pair{\Check
J_0 (m)}{\xi}= \JN{\xi}(m)$
is $\Ad^*$-equivariant.
Note that we follow the convention to denote by
$\Ad^*(g)=(\Ad(g^{-1}))^*$ the coadjoint action of $g$.

Now recall that a $G$-invariant star product is called strongly
$G$-invariant, if $J = J_0$ defines a quantum Hamiltonian, where
$J_0$ is a $G$-equivariant classical momentum map.

Finally, let us briefly introduce the notions of isomorphisms and
equivalence transformations in the $G$-invariant framework. A star
product $\star$ on $(M,\omega)$ is called $G$-isomorphic to the
star product $\star'$ on $(M,\omega')$, if one can find an
isomorphism from $(\Cinf{M}[[\nu]],\star)$ to
$(\Cinf{M}[[\nu]],\star')$ which commutes with $r(g)$ for every
$g\in G$. If in addition $\omega = \omega'$ and there exists an
equivalence transformation from $(\Cinf{M}[[\nu]],\star)$ to
$(\Cinf{M}[[\nu]],\star')$ which commutes with $r(g)$ for every
$g\in G$, one calls $\star$ a star product $G$-equivalent to
$\star'$. From these definitions it is obvious how to define the
notions of $G$-automorphisms and $G$-self equivalences. As an
immediate consequence note that for a $G$-isomorphism $\mathcal T$
from $(\Cinf{M}[[\nu]],\star)$ to $(\Cinf{M}[[\nu]],\star')$ and
$J$  a $G$-equivariant quantum momentum map for $\star$ the
transformed map $J'(\xi):= \mathcal T
\Jbold{\xi}$ is a $G$-equivariant quantum momentum map for $\star'$.
But one has to observe that the notion of strong $G$-invariance is
not preserved under $G$-isomorphisms, in general, that means that
for a strongly $G$-invariant $\star$, a $G$-isomorphic $\star'$
need not be strongly $G$-invariant.
\section{Reduction of Star Products on Cotangent Bundles}
\label{RedSec}
In this section, we present a general procedure for the phase
space reduction of a certain class of star products on cotangent
bundles. As we show in Section \ref{InvQMMSubSec},
this reduction method applies in particular to the star products
$\stark$ and $\starkB$ constructed in the previous section.
Since we will be concerned with different
cotangent bundles $T^*Q$, $T^*\cc{Q}$ at the same instance, we
adopt the convention to use $\cc{\stackrel{}{{}^{}}\quad}$ to
indicate objects related to the cotangent bundle
$\cc{\pi}:T^*\cc{Q}\to \cc{Q}$.
\subsection{Classical Phase Space Reduction of
$(T^*Q,\omega_0+ \pi^*B_0)$: Geometric and Algebraic Properties}
\label{classRedSubSec}
Given a Lie group $G$ we denote by $\{e_i\}_{1\leq i
\leq \dim{(G)}}$ a basis of its Lie algebra $\mathfrak g$ and
by $\{e^i\}_{1\leq i \leq \dim{(G)}}$ the corresponding dual basis
of $\mathfrak g^*$. Assume that $\phi : G \times Q \to Q$ is a left
action on the base manifold $Q$. This action gives rise to a
$G$-action $\Phi : G \times T^*Q \to T^*Q$ on $T^*Q$, the cotangent
lift of $\phi$. Explicitly, it is defined by $\Phi(g,\zeta_x) :=
\Phi_g(\zeta_x) := (T^*\phi_{g^{-1}})(\zeta_x)$ and satisfies
$\pi \circ \Phi_g = \phi_g \circ \pi$. We then have
$\Phi_g^*\theta_0 = \theta_0$ and consequently $\Phi_g^*\omega_0 =
\omega_0$. Let us denote the fundamental vector fields of $\phi$ by
$\xi_\Qind\in\Ginf{TQ}$ and those of $\Phi$ by $\xi_\QKotind \in
\Ginf{T(T^*Q)}$. Clearly, these vector fields are $\pi$-related,
i.e.\ $T \pi\circ \xi_\QKotind = \xi_\Qind \circ \pi$ for all $\xi
\in \mathfrak g$. Recall that a $G$-equivariant classical momentum
map for the canonical symplectic form $\omega_0$ is given by
$\JN{\xi} =\theta_0(\xi_\QKotind)=
\Pol{\xi_\Qind}$, where $\Pol{\xi_\Qind}\in \mathcal P^1(Q)$
denotes the function linear in the momenta corresponding to the
vector field $\xi_\Qind$. In canonical coordinates this reads
$\JN{\xi} = p_i \pi^*\xi_\Qind^i$.

For the symplectic form $\omega_{B_0} = \omega_0 + \pi^*B_0$ with
$B_0\in Z^2_\dR(Q)$ the following well-known result gives first
conditions which have to be satisfied in order to be able to
construct the reduced phase space.
\begin{LEMMA}\label{omegaBInvCMMLem}
Let $\Phi$ act on $(T^*Q,\omega_{B_0})$ as above. Then the following holds
true:
\begin{enumerate}
\item
$\Phi$ is a symplectic action with respect to $\omega_{B_0}$ if and
only if $B_0$ is $G$-invariant.
\item
If $B_0$ is $G$-invariant, then there is a $G$-equivariant
classical momentum map for $\Phi$ if and only if there is a
real-valued element $j_0\in C^1(\mathfrak g,\Cinf{Q})$ such that
\begin{equation}\label{CMMExiCondEq}
\d \jN{\xi} = i_{\xi_\Qind} B_0 \quad\text{ and }\quad
\phi_g^* \jN{\xi} = \jN{\Ad(g^{-1})\xi}
\quad\text{ for all $g\in G$, $\xi \in \mathfrak g$}.
\end{equation}
In this case $\JN{\xi} = \Pol{\xi_\Qind} + \pi^*\jN{\xi}$ defines a
$G$-equivariant classical momentum map which is unique up to
elements of the space ${\mathfrak g^*}^G$ of invariants.
\item
If the relations (\ref{CMMExiCondEq}) are satisfied, one has in
particular
\begin{equation}\label{jNullLieKlamEq}
\jN{[\xi,\eta]} = B_0(\xi_\Qind,\eta_\Qind)\quad
\text{ for all $\xi,\eta\in \mathfrak g$}.
\end{equation}
\end{enumerate}
\end{LEMMA}
\begin{PROOF}
The proof of i.) is obvious, so let us show ii.).
To this end check first that $J_0$ defines a classical
Hamiltonian for $\Phi$ if and only if $\d (\JN{\xi} -
\Pol{\xi_\Qind}) = \pi^*i_{\xi_\Qind}B_0$ for all $\xi \in
\mathfrak g$. After application of $i^*$ one notes that this is equivalent to
the existence of $j_0$ such that the first condition in Eq.\
(\ref{CMMExiCondEq}) is fulfilled. Now observe that the canonical
momentum map $J^0$ for the case $B_0=0$, which is given by
$J^0(\xi)=\Pol{\xi_\Qind}$, is $G$-equivariant. Therefore,
$G$-equivariance of $J_0$ is equivalent to the $G$-equivariance of
$j_0$. The statement about the ambiguity of $J_0$ is a general fact
which holds true for arbitrary Hamiltonian $G$-spaces. Assertion
iii.) is obtained by differentiating $\phi_g^* \jN{\xi} =
\jN{\Ad(g^{-1})\xi}$ with respect to $g$ at $e$ and using the first condition
in Eq.\ (\ref{CMMExiCondEq}).
\end{PROOF}

From now on we will assume that the above conditions for the
existence of a $G$-equivariant classical momentum map are satisfied
and that the action of $G$ on $Q$ is proper and free. This implies
in particular that the orbit space of the $G$-action is smooth and
even that $p : Q \to \cc{Q}= Q/G$ is a left principal $G$-bundle.
Now fix an element $\mu_0\in{\mathfrak g^*}^G$ and recall that it
gives rise to the reduced phase space $\Check J_0^{-1} (\mu_0)
\big/ G$, where $\Check J_0 :T^*Q \to \mathfrak g^*$ is defined by
$\pair{\Check J_0 (\zeta_x)}{\xi}=\JN{\xi}(\zeta_x)$. Since $\mu_0$
is a  regular value of $\Check J_0$, classical Marsden--Weinstein
reduction applies and  the reduced phase space naturally carries
the structure of a symplectic manifold. The induced symplectic form
on the reduced space will be denoted by $\omega_{\mu_0}$. It is
uniquely characterized by the relation
\begin{equation}
  \label{MarWeiIndSymp}
     \pi_{\mu_0}^* \omega_{\mu_0} = i_{\mu_0}^* \omega_{B_0},
\end{equation}
where $i_{\mu_0}: \Check J_0^{-1}(\mu_0) \rightarrow T^*Q$
denotes the inclusion
and $\pi_{\mu_0}$ the projection of $\Check J_0^{-1}(\mu_0)$
onto the orbit space.

In case $B_0=0$, $\JN{\xi}= \Pol{\xi_\Qind}$ it is well known that the
reduced phase space is symplectomorphic to $T^*\cc{Q}$ equipped with a
symplectic structure of the form $\cc{\omega}_0 +
\cc{\pi}^*b_0$, where $b_0$ is a closed two-form on
$\cc{Q}$ which vanishes for $\mu_0=0$.
Note that the construction of an appropriate symplectomorphism
is not canonical unless $\mu_0 =0$ (cf.\ \cite{Got86,Kum81,Mon83}).

An analogous result holds in the case of non-vanishing $B_0$. Since
neither a precise statement of this nor a proof seems to have been
appeared in the literature we briefly present it here. To this end
we need some tools from the theory of (left) principal $G$-bundles.
Let $\gamma$ be a connection one-form on the principal bundle $p :
Q \to \cc{Q}$.  Note that $\gamma$ transforms according to the rule
$\phi_g^*\gamma = \Ad(g) \gamma$ since $G$ acts from the left on
$Q$. Moreover, let $\lambda:= \d\gamma -
\frac{1}{2}[\gamma,\gamma]_\wedge$ denote the corresponding
curvature form; observe the minus sign in front of the bracket
which is due to the fact that we work with a left principal
$G$-bundle. Recall that $\lambda$ is a $\mathfrak g$-valued horizontal
two-form on $Q$.  Finally, we associate to every smooth map
$\Check\jmath : Q \to \mathfrak g^*$ a
one-form $\Gamma_{\Check\jmath} \in \Ginf{T^*Q}$ by
\begin{equation}\label{GamDefEq}
\Gamma_{\Check\jmath} := \pair{\Check\jmath}{\gamma} .
\end{equation}
Then
$\phi_g^* \Gamma_{\Check\jmath}= \Gamma_{\Ad^*(g^{-1})\phi_g^*\Check\jmath}$
for all $g\in G$. Thus, if $\Check\jmath$ is $G$-equivariant, then
$\Gamma_{\Check\jmath}$ is $G$-invariant as well.
Now we are prepared to formulate our result.
\begin{THEOREM}{\bf (cf.\ \cite{Got86,Kum81,Mon83})}
\label{BNullRedThm}
With notations and assumptions from above, the reduced phase
space $\Check J_0^{-1}(\mu_0) / G$ with symplectic
form $\omega_{\mu_0}$ induced by $\omega_{B_0}$
is symplectomorphic to $(T^*\cc{Q},\omega_{b_0} ) =
(T^*(Q/G), \cc{\omega}_0+
\cc{\pi}^*b_0)$, where $b_0$ is the uniquely determined
closed two-form on $\cc{Q}$ which satisfies $p^*b_0 = B_0 + \d
\Gamma_{\Check j_0 - \mu_0}$:
\begin{equation}
((T^*Q)_{\mu_0}, \omega_{\mu_0}) := (\Check J_0^{-1} (\mu_0) /
G,\omega_{\mu_0}) \cong (T^*\cc{Q}, \omega_{b_0}).
\end{equation}
\end{THEOREM}
\begin{PROOF}
Let us consider the one-form $\Gamma_{\Check j_0-\mu_0}$, where the
smooth mapping $\Check j_0 :Q\to \mathfrak g^*$ is defined by
$\pair{\Check j_0(x)}{\xi}=\jN{\xi}(x)$. This one-form induces
a fiber translation $t_{\Gamma_{\Check j_0-\mu_0}}$ given by
$t_{\Gamma_{\Check j_0-\mu_0}}(\zeta_x) = \zeta_x +
\Gamma_{\Check j_0-\mu_0}(x)$. Due to
the properties of the connection one-form $t_{\Gamma_{\Check j_0-\mu_0}}$
maps $\Check J_0^{-1}(\mu_0)$ to
$(T^*Q)^0:=\{\zeta_x\in T^*Q\,|\,\zeta_x(\xi_\Qind(x))=0 \text{ for all }
\xi\in \mathfrak g\}$.
Next observe that $t_{\Gamma_{\Check j_0-\mu_0}}$ commutes
with every $\Phi_g$ by the
$G$-invariance of $\mu_0$ and the $G$-equivariance of $j_0$.
Therefore, $t_{\Gamma_{\Check j_0-\mu_0}}$ passes to the quotient
and defines a diffeomorphism $\Psi_{\mu_0}$ from $\Check
J_0^{-1}(\mu_0) / G$ to $(T^*Q)^0/G \cong T^*(Q/G)$. In order to
determine the symplectic form that is carried over to
$T^*\cc{Q}=T^*(Q/G)$ via $(\Psi_{\mu_0}^{-1})^*$, one first has to
compute $t^*_{-\Gamma_{\Check j_0-\mu_0}}\omega_{B_0}= \omega_0 +
\pi^*(B_0 + \d \Gamma_{\Check j_0-\mu_0})$. Now an easy computation
using the relations between $B_0$ and $j_0$ yields that
\begin{equation}\label{MagHorRel}
  B_0 + \d \Gamma_{\Check j_0-\mu_0} = \mathrm{H} (B_0) + \pair{\Check j_0 -
 \mu_0}{\lambda},
\end{equation}
where $\mathrm{H}(B_0) = B_0 -\sum_{i=1}^{\dim{(G)}}
\Gamma_{e^i}\wedge i_{{e_i}_\Qind}B_0 - \frac{1}{2}
\sum_{i,k=1}^{\dim{(G)}} \Gamma_{e^i}\wedge \Gamma_{e^k}
i_{{e_i}_\Qind} i_{{e_k}_\Qind}B_0$ denotes the totally horizontal
part of $B_0$. But from Eq.\ (\ref{MagHorRel}) it follows
that $B_0 + \d \Gamma_{\Check j_0-\mu_0}$ is horizontal. Due to
the fact that it is $G$-invariant and that $p$ is a surjective
submersion, this implies the existence of a unique closed two-form
$b_0$ as stated in the theorem. But then $\Psi_{\mu_0}$ induces
a symplectic form on $T^*\cc{Q}$. Explicitly, this form is given by
$\omega_{b_0}=\cc{\omega}_0 + \cc{\pi}^*b_0$ which can be derived
from the defining equation
$i_{\mu_0}^*\omega_{B_0} = \pi_{\mu_0}^*\omega_{\mu_0}$
and the commutative diagram
\begin{equation}\label{ComClassRedDiag}
\begin{CD}
\Check
J_0^{-1}(\mu_0) @>t_{\Gamma_{\Check j_0-\mu_0}}\circ i_{\mu_0}>>
(T^*Q)^0 @> \pi \circ i^0>> Q \\ @VV\pi_{\mu_0}V @VV \pi^0 V @VV p
V\\ \Check J_0^{-1}(\mu_0)/G @>\Psi_{\mu_0}>> (T^*Q)^0/G \cong
T^*\cc{Q}@>
\cc{\pi} >> \cc{Q},
\end{CD}
\end{equation}
where $i_{\mu_0},i^0$ are the inclusions into $T^*Q$ and
$\pi_{\mu_0},\pi^0$ the projections onto the respective orbit spaces.
\end{PROOF}

Now consider the space $\Cinf{T^*Q}^G$ of
$G$-invariant smooth functions on $T^*Q$ and the space
\[
 I_{\mu_0,\cdot}:=\Big\{f \in \Cinf{T^*Q}^G\, \Big|\;
 f = \sum_{i = 1}^{\dim{(G)}}h^i (\JN{e_i} - \pair{\mu_0}{e_i})
 \text{ with } h^i\in\Cinf{T^*Q} \Big\}.
\]
Since $I_{\mu_0,\cdot}$ is a Poisson
ideal in $\Cinf{T^*Q}^G$, the pointwise product and the
Poisson bracket $\{\,\,,\,\,\}_{B_0}$ corresponding to $\omega_{B_0}$
induce the structure of a Poisson algebra on the quotient $\Cinf{T^*Q}^G
/I_{\mu_0,\cdot}$ by
\begin{equation}\label{QuoPoiStrEq}
[f]_{\mu_0,\cdot}
\cdot_\red [f']_{\mu_0,\cdot} := [f f']_{\mu_0,\cdot},\quad
\{[f]_{\mu_0,\cdot},[f']_{\mu_0,\cdot}\}_\red:=
[\{f,f'\}_{B_0}]_{\mu_0,\cdot}.
\end{equation}
The thus obtained Poisson algebra is known to be isomorphic to
$(\Cinf{T^*\cc{Q}}, \{ \,\,,\,\,\}_{b_0})$, where $\{
\,\,,\,\,\}_{b_0}$ denotes the Poisson bracket corresponding to the
symplectic form $\omega_{b_0}$. Unfortunately, there is no canonical
construction of such an isomorphism, but
restricting to functions polynomial in the
momenta, we are able to find a natural isomorphism between the
Poisson subalgebras $\PolFun^G /I_{\mu_0,\cdot}^{\mathcal
P}\subseteq \Cinf{T^*Q}^G /I_{\mu_0,\cdot}$ and $\PolFunRed\subseteq
\Cinf{T^*\cc{Q}}$ which depends only on the choice of the connection
$\gamma$. Hereby we have used the abbreviations $\PolFun^G:=
\PolFun\cap \Cinf{T^*Q}^G$ and $I_{\mu_0,\cdot}^{\mathcal P}:=
\PolFun\cap I_{\mu_0,\cdot}$.

Let us now provide the details. Due to the choice of the connection $\gamma$
the tangent bundle of $Q$ can be written as the direct sum of the horizontal
bundle $HQ$ and the (trivial) vertical bundle $VQ$.
Clearly, this decomposition induces a decomposition of the symmetric powers of
$TQ$. Hence the space of sections $\Ginf{\bigvee TQ}$ can be
written as
\[
\Ginf{\mbox{$\bigvee$} TQ} = \bigoplus_{k=0}^\infty
\Ginf{\mbox{$\bigvee^k$} HQ}
\oplus \bigoplus_{k=1}^\infty \bigoplus_{r=1}^k \Ginf{\mbox{
$\bigvee^{k-r}$} HQ \vee \mbox{$\bigvee^r$} VQ}.
\]
Obviously, $\Ginf{\bigvee^l HQ \vee \bigvee^r VQ}$ is bigraded by the
horizontal degree $l$ and the vertical degree $r$. We will refer to the spaces
$\bigoplus_{k=0}^\infty \Ginf{\bigvee^k HQ}$ and
$\bigoplus_{k=1}^\infty \bigoplus_{r=1}^k \Ginf{\bigvee^{k-r} HQ
\vee \bigvee^r VQ}$ as the space of totally horizontal sections and
partially vertical sections in $\bigvee TQ$, respectively. Moreover,
we denote by $\mathrm{H}$ the projection onto the totally horizontal sections
and by $\mathrm{PV}$ the projection onto the partially vertical
sections. Since $\{{e_i}_\Qind\}_{1\leq i \leq
\dim{(G)}}$ is a set of basis sections of the vertical bundle, there
exist for every $T \in \mathrm{PV}(\Ginf{\bigvee TQ})$
uniquely determined tensor fields $\mathrm R^i(T) \in
\Ginf{\bigvee TQ}$ such that $T = \sum_{i=1}^{\dim{(G)}}\mathrm R^i
(T) \vee {e_i}_\Qind$. Hence $\mathrm R^i:\mathrm{PV}
(\Ginf{\bigvee TQ}) \to \Ginf{\bigvee TQ}$ is a well-defined
mapping that extends to all of $\Ginf{\bigvee TQ}$ by setting
$\mathrm R^i (T):=0$ for $T \in \mathrm{H}(\Ginf{\bigvee TQ})$.
Using the isomorphism $\mathsf P : \Ginf{\bigvee TQ} \to \PolFun$ we then get
the following decomposition of $\PolFun$ into the spaces of so-called
totally horizontal and partially vertical polynomial functions:
\[
\PolFun = \hor(\PolFun) \oplus \mathrm{pv}(\PolFun).
\]
Hereby, we have used
$\hor = \mathsf P \circ \mathrm{H} \circ \mathsf P^{-1}$ and
$\mathrm{pv} = \mathsf P \circ \mathrm{PV} \circ
\mathsf P^{-1}$. Under the
isomorphism $\mathsf P$ the mapping $\mathrm R^i$ transforms
to $\mathrm r^i: \PolFun \to \PolFun$ which explicitly is given
by
\begin{equation}\label{RestAbbDefEq}
\mathrm r^i (F) = \left\{
    \begin {array} {cl}
    \frac{1}{r} \Fop{\Gamma_{e^i}} F & \mbox { if }
    F \mbox{ is vertical of degree } r \geq 1,\\
    0 & \mbox { if } F \mbox{ is vertical of degree } 0.
    \end {array}
    \right.
\end{equation}
Thus, every $F\in \PolFun$ can be written as
\[
F = \hor(F) + \sum_{i=1}^{\dim{(G)}}\mathrm r^i(F)
\Pol{{e_i}_\Qind}.
\]
Now consider again the $G$-invariant one-form
$\Gamma_{\Check j_0 - \mu_0}$ on $T^*Q$ defined in the proof of
Theorem \ref{BNullRedThm}.
Observe that by $\gamma(Y)=0$ for every $Y\in \Ginf{HQ}$ the equality
$t^*_{\Gamma_{\Check j_0-\mu_0}} F = F$ is satisfied for every
totally horizontal polynomial function $F$.
Hence every $F\in \PolFun$ can be written in the form
\begin{equation}\label{ClassicalDecompEq}
F = \hor(t^*_{- \Gamma_{\Check j_0 -\mu_0}}F) +
\sum_{i=1}^{\dim{(G)}} t^*_{\Gamma_{\Check j_0- \mu_0}}
\mathrm r^i(t^*_{-\Gamma_{\Check j_0 - \mu_0}}F)
(\JN{e_i} - \pair{\mu_0}{e_i}),
\end{equation}
where we have used that $t^*_{\Gamma_{\Check j_0 - \mu_0}}
\Pol{\xi_\Qind}= \Pol{\xi_\Qind} + \pi^*\jN{\xi} - \pair{\mu_0}{\xi}$
for $\xi \in \mathfrak g$. After these rather technical
preparations we can now prove the following result.
\begin{PROPOSITION}\label{ClassRedProp}
\begin{enumerate}
\item
The space $\PolFun$ decomposes into the direct sum
\begin{equation}\label{PolDecGenEq}
\PolFun = \hor(\PolFun) \oplus \Big\{ F \in \PolFun\,\Big|\,
 F = \sum_{i=1}^{\dim{(G)}} H^i(\JN{e_i} - \pair{\mu_0}{e_i})
 \text{ with } H^i \in \PolFun \Big\}
\end{equation}
and this decomposition is $G$-invariant. Moreover, the projections
onto the respective subspaces are given by $\hor_{\Check j_0
- \mu_0}:=\hor \circ t^*_{-\Gamma_{\Check j_0
-\mu_0}}$ and $\mathrm{pv}_{\Check j_0 - \mu_0}:=
t^*_{\Gamma_{\Check j_0-\mu_0}}\circ \mathrm{pv}\circ
t^*_{-\Gamma_{\Check j_0 -\mu_0}}$.
\item
According to i.) the space of $G$-invariant polynomial functions
$\PolFun^G$ decomposes into the direct sum
\begin{equation}\label{InvPolDecEq}
\PolFun^G = \hor(\PolFun^G) \oplus I^{\mathcal P}_{\mu_0,\cdot}.
\end{equation}
\item
The space $\hor(\PolFun^G)$ of totally horizontal $G$-invariant
polynomial functions becomes a Poisson algebra with the usual
pointwise product of functions and the Poisson bracket
$\{\,\,,\,\,\}_{\Check j_0 - \mu_0}$ defined by
\begin{equation}\label{PoiBraMuDefEq}
\{F, F'\}_{\Check j_0 - \mu_0}:= \hor_{\Check j_0 - \mu_0}
(\{F,F'\}_{B_0})= \hor (t^*_{-\Gamma_{\Check j_0 -\mu_0}}
\{F,F'\}_{B_0}),\quad F,F' \in \hor(\PolFun^G).
\end{equation}
\item
As a Poisson algebra, $(\hor(\PolFun^G), \{\,\,,\,\,\}_{\Check
j_0-\mu_0})$ is isomorphic to $\PolFunRed$ with the
Poisson bracket induced by $\omega_{b_0}
= \cc{\omega}_0 + \cc{\pi}^* b_0$, where $b_0$
denotes the uniquely determined closed two-form on $\cc{Q}$ such
that $p^*b_0 = B_0 + \d \Gamma_{\Check j_0 - \mu_0}$. If
${}^\hor:\Ginf{\bigvee T\cc{Q}} \to \Ginf{\bigvee TQ}$ denotes the
horizontal lift (which is obtained by extension from
$\Ginf{T\cc{Q}}$ to $\Ginf{\bigvee T\cc{Q}}$ as homomorphism with
respect to $\vee$, particularly $\chi^\hor = p^*\chi$ for $\chi\in
\Cinf{\cc{Q}}$),
an explicit Poisson algebra isomorphism is given by
\begin{equation}\label{LiftIsoEq}
l : \PolFunRed \ni \PolRed{t} \mapsto \Pol{t^\hor} \in
\hor(\PolFun^G), \quad t \in \Ginf{\mbox{$\bigvee$} T\cc{Q}}.
\end{equation}
\item
Finally, $(\hor(\PolFun^G), \{\,\,,\,\,\}_{\Check j_0-\mu_0})$ is
isomorphic to $\PolFun^G / I^{\mathcal P}_{\mu_0,\cdot}$ with
Poisson algebra structure defined by Eq.\ (\ref{QuoPoiStrEq}).
An isomorphism is given by
\begin{equation}
\hor(\PolFun^G) \ni F \mapsto [F]_{\mu_0,\cdot} \in
\PolFun^G / I^{\mathcal P}_{\mu_0,\cdot}.
\end{equation}
\end{enumerate}
\end{PROPOSITION}
\begin{PROOF}
Statements i.) and ii.) are obvious from the above
considerations. Since the kernel of $\hor_{\Check j_0
- \mu_0}|_{\PolFun^G}$ is a Poisson ideal in $\PolFun^G$,
it is straightforward to verify that
$\{\,\,,\,\,\}_{\Check j_0 - \mu_0}$  defines a Poisson
bracket on $\hor(\PolFun^G)$. For the proof of iv.), first note
that $l:\PolFunRed \to\hor(\PolFun^G)$ is an isomorphism of
commutative algebras. Thus it remains to show that the
Poisson bracket obtained on $\PolFunRed$ by pull-back by $l$ of the
Poisson bracket $\{\,\,,\,\,\}_{\Check j_0 - \mu_0}$ on
$\hor(\PolFun^G)$ coincides with the Poisson bracket
induced by $\omega_{b_0}$. To check this, it is enough to compute
the induced bracket of
$\PolRed{t}$ and $\PolRed{s}$ for $t,s \in \Ginf{T\cc{Q}}$.
By a straightforward computation one obtains
\[
\{l(\PolRed{t}), l(\PolRed{s})\}_{\Check j_0-\mu_0} =
- \Pol{[t,s]^\hor} - \pi^*((B_0 +
\d \Gamma_{\Check j_0 - \mu_0})(t^\hor, s^\hor)).
\]
Since $Tp\,t^\hor = t \circ p$ and analogously for $s$, we get $(B_0 +
\d \Gamma_{\Check j_0 - \mu_0})(t^\hor, s^\hor) = p^*(b_0(t,s))$
by definition of $b_0$.
Then, using $l(\cc{\pi}^*\chi)=\pi^* p^*\chi$ for $\chi\in
\Cinf{\cc{Q}}$, we find
\[
l^{-1}( \{l(\PolRed{t}), l(\PolRed{s})\}_{\Check j_0-\mu_0} ) = -
\PolRed{[t,s]}- \cc{\pi}^*(b_0(t,s)).
\]
But this proves iv.) since this last expression coincides with the
Poisson bracket of $\PolRed{t}$ and $\PolRed{s}$ with respect
to the symplectic form $\omega_{b_0}
= \cc{\omega}_0 + \cc{\pi}^* b_0$. For the proof of v.), one again
has to check the compatibility of the Poisson brackets, but
this is straightforward observing that the kernel of $\hor_{\Check
j_0 - \mu_0}$ restricted to $\PolFun^G$ is the ideal
$I_{\mu_0,\cdot}^{\mathcal P}$.
\end{PROOF}
\subsection{Reduction of a Certain Class of Star Products on
Cotangent Bundles}
\label{starRedSubSec}
In view of the reduction of the classical structures considered in
the preceding section, it should  be possible to analogously
determine the reduction of certain star products on
$(T^*Q,\omega_{B_0})$ in order to obtain star products on
$(T^*\cc{Q},\omega_{B_0})$. To this end recall first our general
assumptions that the $G$-action on $Q$ is proper and free and that
$\mu_0 \in \mathfrak g^*$, which is a regular value of $\Check
J_0$, is chosen to be invariant with respect to the coadjoint
action of $G$. For the quantum reduction we have to assume in
addition that we are given a star product $\star$ on
$(T^*Q,\omega_{B_0})$, a $G$-equivariant quantum momentum map $J$
and a deformation $\mu$ of the classical momentum value $\mu_0$
such that the following properties hold true:
\begin{itemize}
\item
$\star$ is $G$-invariant, i.e.\ invariant with respect to the
lifted action $\Phi$ of $G$ on $T^*Q$.
\item
$J = J_0 + J_+\in C^1(\mathfrak g,\Cinf{T^*Q})[[\nu]]$ is a
$G$-equivariant quantum momentum map for $\star$, where $J_0$
denotes a $G$-equivariant classical momentum map of the form
$\JN{\xi} = \Pol{\xi_\Qind}+ \pi^*\jN{\xi}$ as in Lemma
\ref{omegaBInvCMMLem}.
\item
$\PolFun[[\nu]]$ is a $\star$-subalgebra and $\Jbold{\xi}\in
\PolFun[[\nu]]$ for all $\xi\in \mathfrak g$.
\item
The quantum momentum value has the form
$\mu = \mu_0 + \mu_+$ with $\mu_+\in\nu\mathfrak g_c^*[[\nu]]$ and
$\mathfrak g^*_c$ the complexification of $\mathfrak g^*$. Moreover,
$\mu$ is invariant with respect to the coadjoint action of $G$.
\end{itemize}

The third of the above assumptions will enable us to compute the
reduced star product by means of polynomial functions, only.
Note at this point that later in Corollary \ref{GenQMMExiCor} ii.)
we will prove that the assumption $\Jbold{\xi}\in \PolFun[[\nu]]$
for all $\xi\in \mathfrak g$ is actually no
additional assumption but a consequence of the fact that
$\PolFun[[\nu]]$ is a $\star$-subalgebra.

As for the algebraic part of the classical reduction
consider $\Cinf{T^*Q}^G[[\nu]]$ and the subspace
\[
 I_{\mu,\star}:= \Big\{f \in\Cinf{T^*Q}^G[[\nu]]\, \Big| \;
 f = \sum_{i=1}^{\dim{(G)}} h^i \star
 (\Jbold{e_i} - \pair{\mu}{e_i}) \text{ with }
 h^i \in \Cinf{T^*Q}[[\nu]]\Big\} .
\]
Then $I_{\mu,\star}$ is a two-sided ideal of $\Cinf{T^*Q}^G[[\nu]]$,
since $J$ is a quantum Hamiltonian, and the quotient space
$\Cinf{T^*Q}^G[[\nu]]/I_{\mu,\star}$
becomes an associative algebra by
\begin{equation}\label{QuotDefProdDefEq}
[f]_{\mu,\star} \star_\red [f']_{\mu,\star} := [f \star
f']_{\mu,\star}.
\end{equation}
In order to interpret the associative product
$\star_\red$ as a star product on the reduced phase space, one
has to find a $\mathbb C[[\nu]]$-module isomorphism
between $\Cinf{T^*Q}^G[[\nu]]/I_{\mu,\star}$ and $(\Cinf{T^*Q}^G /
I_{\mu_0,\cdot})[[\nu]]$. Since the latter space is isomorphic to
$\Cinf{T^*\cc{Q}}[[\nu]]$ one then defines an associative product
on $\Cinf{T^*\cc{Q}}[[\nu]]$ by declaring the isomorphism in question
to be a homomorphism of associative algebras.
If the thus constructed product is a star product indeed, then
we have obtained the reduced star product we are looking for.
At this point one has to mention that there is no canonical construction of
such an isomorphism for the space of all smooth functions, but like
in the classical case one can find a natural isomorphism
between the subspaces obtained by restriction to polynomial functions.

Let us now provide the details for the construction of this
latter isomorphism. To this end observe first that by the third
assumption above the product $\star_\red$ can be restricted
to $\PolFun^G[[\nu]]/ I^{\mathcal
P}_{\mu,\star} \subseteq \Cinf{T^*Q}^G[[\nu]]/I_{\mu,\star}$, where
$I^{\mathcal P}_{\mu,\star} := \PolFun [[\nu]] \cap I_{\mu,\star}$.
We now claim that there is a naturally constructed $\mathbb
C[[\nu]]$-module isomorphism between $\PolFun^G[[\nu]]/ I^{\mathcal
P}_{\mu,\star}$ and $(\PolFun^G/ I^{\mathcal
P}_{\mu_0,\cdot})[[\nu]]$ where the latter space is isomorphic to
$\PolFunRed [[\nu]]$ by Proposition \ref{ClassRedProp} iv.) and v.).

Clearly, the projections $\hor_{\Check j_0 - \mu_0}$,
$\mathrm{pv}_{\Check j_0 - \mu_0}$ and the maps $\mathrm
r^i_{\Check j_0- \mu_0}:= t^*_{\Gamma_{\Check j_0-\mu_0}}
\circ \mathrm r^i \circ t^*_{-\Gamma_{\Check j_0-\mu_0}}$ extend by
$\mathbb C[[\nu]]$-linearity from $\PolFun$ to $\PolFun[[\nu]]$.
Hence, by Proposition \ref{ClassRedProp} i.), every $F\in
\PolFun[[\nu]]$ decomposes uniquely into a sum of the form
\[
F = \hor_{\Check j_0 - \mu_0}(F) + \mathrm{pv}_{\Check j_0
-\mu_0}(F) = \hor_{\Check j_0 - \mu_0}(F) +
\sum_{i=1}^{\dim{(G)}} \mathrm r^i_{\Check j_0 -\mu_0}(F)
(\JN{e_i}- \pair{\mu_0}{e_i}).
\]
After defining $\triangle_{\mu,\star} : \PolFun[[\nu]] \to
\PolFun[[\nu]]$ by the equation
\begin{equation}\label{triangleDefEq}
\triangle_{\mu,\star} F := \frac{1}{\nu} \sum_{i=1}^{\dim{(G)}}
\left(\mathrm r^i_{\Check j_0 - \mu_0}(F)
(\JN{e_i}- \pair{\mu_0}{e_i}) - \mathrm r^i_{\Check j_0 -
\mu_0}(F)\star (\Jbold{e_i}- \pair{\mu}{e_i})
\right),
\end{equation}
the above decomposition can be rewritten as
\[
F = \hor_{\Check j_0 - \mu_0}(F) + \sum_{i=1}^{\dim{(G)}}
\mathrm r^i_{\Check j_0 -\mu_0}(F)\star (\Jbold{e_i}-
\pair{\mu}{e_i}) + \nu
\triangle_{\mu,\star} F.
\]
Now repeat this ad infinitum and decompose at every step the
remaining term not of the form $\hor_{\Check j_0 - \mu_0}(F') +
\sum_{i=1}^{\dim{(G)}} \mathrm r^i_{\Check j_0-
\mu_0}(F')\star (\Jbold{e_i}- \pair{\mu}{e_i})$. This procedure finally
yields
\begin{equation}\label{DeformedDecompEq}
F = \hor_{\Check j_0 - \mu_0} \left(\frac{\id}{\id - \nu
\triangle_{\mu,\star}}F\right) +
\sum_{i=1}^{\dim{(G)}} \mathrm r^i_{\Check j_0 - \mu_0}\left(
\frac{\id}{\id - \nu \triangle_{\mu,\star}}F\right)
\star (\Jbold{e_i}- \pair{\mu}{e_i}).
\end{equation}
Like in the classical case we obtain:
\begin{LEMMA}\label{DefDecLem}
\begin{enumerate}
\item
The space $\PolFun[[\nu]]$ decomposes into the direct sum
\begin{equation}\label{DefPolDecGenEq}
  \begin{split}
   &\PolFun[[\nu]] = \\
   & \hspace{0.5em} \hor(\PolFun)[[\nu]] \oplus\!
   \Big\{ \!F \in \PolFun[[\nu]]\,\Big|\;
   F =\!\!\! \sum_{i=1}^{\dim{(G)}} \!\!\! H^i\star(\Jbold{e_i} -\pair{\mu}{e_i})
   \text{ with } H^i \in \PolFun[[\nu]]
   \Big\}
  \end{split}
\end{equation}
and this decomposition is $G$-invariant. Moreover, the projections
onto the respective subspaces are given by
\[
 \hor_{\Check j_0
 - \mu_0}\circ \frac{\id}{\id - \nu \triangle_{\mu,\star}}
 \quad \text{and}
 \quad \sum_{i=1}^{\dim{(G)}} \!\mathrm r^i_{\Check j_0 - \mu_0}\left(
 \frac{\id}{\id - \nu \triangle_{\mu,\star}}F\right)
 \star (\Jbold{e_i}- \pair{\mu}{e_i}), \quad F\in \PolFun[[\nu]].
\]
\item
According to i.) the space of formal series of $G$-invariant
polynomial functions $\PolFun^G[[\nu]]$ decomposes into the direct
sum
\begin{equation}\label{DefInvPolDecEq}
\PolFun^G[[\nu]] = \hor(\PolFun^G) [[\nu]]\oplus
I^{\mathcal P}_{\mu,\star}.
\end{equation}
\end{enumerate}
\end{LEMMA}
\begin{PROOF}
The fact that every element of $\PolFun[[\nu]]$ can be decomposed
as stated in i.) is obvious from the above considerations. To see
that the sum in (\ref{DefPolDecGenEq}) is direct, one just has to
observe that the lowest order in $\nu$ of $\sum_{i=1}^{\dim{(G)}}
H^i\star(\Jbold{e_i} -\pair{\mu}{e_i})$ is given by
$\sum_{i=1}^{\dim{(G)}} H^i_0(\JN{e_i} -\pair{\mu_0}{e_i})$ (where
we have written $H^i =
\sum_{k=0}^\infty \nu^k H^i_k$). But then assuming that
$\sum_{i=1}^{\dim{(G)}} H^i\star(\Jbold{e_i}
-\pair{\mu}{e_i})$ is horizontal, the fact that the sum in
Eq.\ (\ref{PolDecGenEq}) is direct implies that $H^i_0=0$ for all
$1\leq i \leq \dim{(G)}$. Repeating this argument order by order,
we finally get $H^i=0$ for all $1\leq i \leq \dim{(G)}$ proving
that the decomposition is direct. Using the invariance properties
of $\mu$, $J$ and $\star$, it is easy to check that $\Phi_g^*$ maps
elements of the form $\sum_{i=1}^{\dim{(G)}} H^i\star(\Jbold{e_i}
-\pair{\mu}{e_i})$ to elements of the same form. Thus,
the above decomposition turns out to be $G$-invariant, since
obviously $\Phi_g^*$ maps totally horizontal polynomials to
totally horizontal polynomials. ii.) is a direct consequence of i.)
and the definition of $I^{\mathcal P}_{\mu,\star}$.
\end{PROOF}

After these preparations we obtain one of the main results of
this section.
\begin{THEOREM}\label{StarProdRedThm}
\begin{enumerate}
\item
The space $\hor(\PolFun^G)[[\nu]]$ of formal series of totally
horizontal $G$-invariant polynomial functions becomes an
associative algebra by means of the product $\bullred{J,\mu}$
defined by
\begin{equation}\label{AstMuDefEq}
F\bullred{J,\mu} F':= \hor_{\Check j_0 - \mu_0}\left(
\frac{\id}{\id - \nu \triangle_{\mu,\star}}(F \star F')\right),
\qquad F,F' \in \hor(\PolFun^G)[[\nu]].
\end{equation}
\item
The pull back of $\bullred{J,\mu}$ to $\PolFunRed[[\nu]]$ via the
isomorphism $l : \PolFunRed[[\nu]] \to \hor(\PolFun^G)[[\nu]]$
defined in Eq.\ (\ref{LiftIsoEq}) gives rise to a star product
$\starred{J,\mu}$ on $\PolFunRed[[\nu]]$, where the underlying
Poisson bracket is induced by the symplectic form
$\omega_{b_0} = \cc{\omega}_0 + \cc{\pi}^* b_0$ given in
Proposition \ref{ClassRedProp} iv.).
\item
$(\hor(\PolFun^G)[[\nu]],\bullred{J,\mu})$ is isomorphic to
$\PolFun^G[[\nu]] / I^{\mathcal P}_{\mu,\star}$ with the
associative algebra structure defined in Eq.\
(\ref{QuotDefProdDefEq}). An isomorphism is given by
\begin{equation}
\hor(\PolFun^G)[[\nu]] \ni F \mapsto [F]_{\mu,\star} \in
\PolFun^G[[\nu]] / I^{\mathcal P}_{\mu,\star}.
\end{equation}
\item
The star product $\starred{J,\mu}$ on $\PolFunRed[[\nu]]$ can be
described by bidifferential operators, hence can be uniquely
extended to a star product on $\Cinf{T^*\cc{Q}}[[\nu]]$, which
also will be denoted by $\starred{J,\mu}$.
\end{enumerate}
\end{THEOREM}
\begin{PROOF}
Using the fact that the kernel of the projection $\hor_{\Check j_0-
\mu_0}\circ \frac{\id}{\id - \nu \triangle_{\mu,\star}}$
restricted to $\PolFun^G[[\nu]]$ is a two-sided ideal in
$\PolFun^G[[\nu]]$, it is straightforward to see that the
composition $\bullred{J,\mu}$ defined in (\ref{AstMuDefEq}) on
$\hor(\PolFun^G)[[\nu]]$ is associative. For the proof of ii.),
first note that $F \bullred{J,\mu} 1 =
\hor_{\Check j_0 - \mu_0}( \frac{\id}{\id - \nu
\triangle_{\mu,\star}} F) = F = 1\bullred{J,\mu} F$ for all
$F \in \hor(\PolFun^G)[[\nu]]$ since $\triangle_{\mu,\star} F=0$
and $t^*_{-\Gamma_{\Check j_0 - \mu_0}}F=F$ for totally horizontal
polynomial functions, implying that $f \starred{J,\mu} 1 = f = 1
\starred{J, \mu} f$ for $f\in \PolFunRed[[\nu]]$. Now, an immediate
computation yields
\[
F\bullred{J,\mu} F' - F' \bullred{J,\mu} F = \nu \hor_{\Check
j_0-\mu_0} (\{F,F'\}_{B_0})+ O(\nu^2), \quad \text{ for all } F,F' \in
\hor(\PolFun^G),
\]
but according to statement iv.) of Proposition \ref{ClassRedProp}
this implies that for $f,f' \in \PolFunRed$ the lowest order in $f
\starred{J,\mu} f' - f' \starred{J,\mu} f$ is given by the Poisson
bracket corresponding to $\omega_{b_0}$. Assertion iii.) is obvious
from Lemma \ref{DefDecLem} ii.) and the definition of the
associative product on $\hor(\PolFun^G)[[\nu]]$ resp.\
$\PolFun^G[[\nu]] / I^{\mathcal P}_{\mu,\star}$. Since $\star$ is a
differential star product and $\triangle_{\mu,\star}$ a
differential operator, it is obvious that the product
$\bullred{J,\mu}$ on $\hor(\PolFun^G)[[\nu]]$ can be described by
bidifferential operators. But this also implies that the
corresponding star product on $\PolFunRed[[\nu]]$ is given by
bidifferential operators. Since bidifferential operators are
completely determined by their values on polynomial functions,
$\starred{J,\mu}$ extends to a star product on
$\Cinf{T^*\cc{Q}}[[\nu]]$ in a unique way.
\end{PROOF}

\begin{REMARK}\label{PullBackRem}
Clearly, pulling back the star product $\starred{J,\mu}$ on
$\Cinf{T^*\cc{Q}}[[\nu]]$ to $\Cinf{(T^*Q)_{\mu_0}}[[\nu]]$ via the
symplectomorphism $\Psi_{\mu_0}$ constructed in the proof of
Theorem \ref{BNullRedThm} one obtains a star product
$\star_{\Psi_{\mu_0}}^{J,\mu}:= \Psi_{\mu_0}^*\starred{J,\mu}$ on
the symplectic manifold $(\Check J_0^{-1}(\mu_0)/G,\omega_{\mu_0})$.
All results we derive in the sequel about the
star products $\starred{J,\mu}$ will then
transfer almost litteraly to the star products
$\star_{\Psi_{\mu_0}}^{J,\mu}$.
\end{REMARK}

In the remaining section we prove several properties of the
reduction scheme introduced above. This will partially clarify the
dependence of $\starred{J,\mu}$ on the chosen $\mu \in {\mathfrak
g^*}^G + \nu{\mathfrak g^*_c}^G[[\nu]]$. Moreover, we thus explain
the relation between automorphisms resp.\ derivations of $\star$
and those of $\starred{J,\mu}$. In particular, we show that our
reduction is natural with respect to $G$-isomorphisms of star
products which preserve $\PolFun[[\nu]]$ and satisfy an appropriate
compatibility relation for the momentum values.
\begin{PROPOSITION}\label{IsoAutoDerRedProp}
\begin{enumerate}
\item
Let $\bullred{J,\mu}$ and ${\bullet'}^{J',\mu'}$ denote the
products on $\hor(\PolFun^G)[[\nu]]$ obtained by reduction of
$(\star,J,\mu)$ and $(\star',J',\mu')$. Moreover, let $\mathcal T$ be
an isomorphism (equivalence transformation) from $\star$ to
$\star'$ which preserves $\PolFun[[\nu]]$ and which commutes with every
$\Phi_g^*$. Define
$\mathsf T:\hor(\PolFun^G)[[\nu]]\to \hor(\PolFun^G)[[\nu]]$ by
\begin{equation}\label{mathsfTDefEq}
\mathsf T F := \hor_{\Check j'_0 - \mu'_0} \left(
\frac{\id}{\id- \nu \triangle'_{\mu',\star'}}\mathcal T F\right),
\quad F\in \hor(\PolFun^G)[[\nu]].
\end{equation}
Then $\mathsf T$
is an  isomorphism (equivalence transformation) from
$\bullred{J,\mu}$ to ${\bullet'}^{J',\mu'}$, if $\mu' = \mu +
\delta \mu$, where $\delta\mu\in {\mathfrak g^*}^G + \nu {\mathfrak
g^*_c}^G[[\nu]]$ is given by $\pair{\delta\mu}{\xi} = J'(\xi) -
\mathcal T\Jbold{\xi}$. The inverse is given by
\begin{equation}
\label{mathsfTInvEq}
 \mathsf T^{-1}  F = \hor_{\Check j_0 - \mu_0} \left( \frac{\id}{\id- \nu
 \triangle_{\mu,\star}} \mathcal T^{-1} F\right), \quad
 F\in \hor(\PolFun^G)[[\nu]].
\end{equation}
If $J'(\xi) = \mathcal T\Jbold{\xi}$, then $\mathsf T$ is an
isomorphism (equivalence transformation) from $\bullred{J,\mu}$ to
${\bullet'}^{J',\mu}$.
\item
Let $\bullred{J,\mu}$ and $\bullred{J,\mu'}$ denote the products on
$\hor(\PolFun^G)[[\nu]]$ obtained by reduction of $(\star,J,\mu)$
and $(\star,J,\mu')$. Assume $\mathcal A$ to be an automorphism
of $\star$ (starting with $\id$) which preserves $\PolFun[[\nu]]$
and commutes with every $\Phi_g^*$. Define $\mathsf
A:\hor(\PolFun^G)[[\nu]]\to
\hor(\PolFun^G)[[\nu]]$ by
\begin{equation}
\mathsf A F := \hor_{\Check j_0 - \mu'_0} \left(
\frac{\id}{\id- \nu \triangle_{\mu',\star}}\mathcal A F\right),
\quad F\in \hor(\PolFun^G)[[\nu]].
\end{equation}
Then $\mathsf A$ is an isomorphism (equivalence transformation) from
$\bullred{J,\mu}$ to $\bullred{J,\mu'}$, if $\mu' = \mu +
\delta \mu$, where $\delta\mu\in {\mathfrak g^*}^G + \nu
{\mathfrak g^*_c}^G[[\nu]]$ is given by $\pair{\delta\mu}{\xi} =
J(\xi) - \mathcal A\Jbold{\xi}$.  The inverse is given by
\begin{equation}
  \mathsf A^{-1} F =
  \hor_{\Check j_0 - \mu_0}
  \left( \frac{\id}{\id- \nu \triangle_{\mu,\star}}
  \mathcal A^{-1} F\right), \quad F\in \hor(\PolFun^G)[[\nu]].
\end{equation}
If $\Jbold{\xi} = \mathcal  A\Jbold {\xi}$,
then $\mathsf A$ is an automorphism of $\bullred{J,\mu}$ (starting with $\id$).
\item
Let $\mathcal D$ denote a $\mathbb C[[\nu]]$-linear derivation of
$\star$ which preserves $\PolFun[[\nu]]$ and commutes with every
$\Phi_g^*$. Define $\mathsf D:\hor(\PolFun^G)[[\nu]]
\to \hor(\PolFun^G)[[\nu]]$ by
\begin{equation}
\mathsf D F := \hor_{\Check j_0 -\mu_0} \left(
\frac{\id}{\id- \nu \triangle_{\mu,\star}}\mathcal D F\right), \quad
F\in \hor(\PolFun^G)[[\nu]].
\end{equation}
Then $\mathsf D$
is a $\mathbb C[[\nu]]$-linear derivation of $\bullred{J,\mu}$, if
$\mathcal D \Jbold{\xi}=0$ for all $\xi \in \mathfrak g$.
\item
Every isomorphism resp.\ equivalence transformation, automorphism
(starting with $\id$) or derivation constructed
according to i.), ii.) or iii.) transfers to an isomorphism resp.\
equivalence transformation, automorphism (starting with $\id$), or
derivation on $\PolFunRed[[\nu]]$ with the reduced star product as
algebra structure. Moreover, the induced map extends in a unique
way to an isomorphism resp.\ equivalence transformation,
automorphism (starting with $\id$), or derivation on
$\Cinf{T^*\cc{Q}}[[\nu]]$ with the reduced star product as algebra
structure.
\item
All the above constructions can be localized in the sense that
starting from a local isomorphism resp.\ equivalence transformation,
automorphism (starting with $\id$) or derivation on
$\Cinf{T^*U}[[\nu]]$ resp.\ $\mathcal P(U)[[\nu]]$, where
$U\subseteq Q$ is an open subset of $Q$, one obtains a local
isomorphism resp.\ equivalence transformation, automorphisms (starting
with $\id$) or derivation on $\mathcal P(\cc{U})[[\nu]]$ resp.\
$\Cinf{T^*\cc{U}}[[\nu]]$, where $\cc{U}
= p (U)\subseteq \cc{Q}$.
\end{enumerate}
\end{PROPOSITION}
\begin{PROOF}
For the proof of the first three statements it suffices to show
i.) since ii.) is just a special case of i.) and the proof
of iii.) is only a slight adaption of the argument which proves i.).
Since $\mathcal T$ commutes with $\Phi_g^*$ it also commutes
with $\Lie_{\xi_\QKotind}$ for every $\xi\in \mathfrak g$. By definition of
a quantum Hamiltonian this implies
that $\frac{1}{\nu}\ad_{\star'}(J'(\xi) - \mathcal T
\Jbold{\xi})=0$.
But this entails that for every such isomorphism $\mathcal T$ one
has $J'(\xi) - \mathcal T \Jbold{\xi} \in \mathbb C[[\nu]]$, hence
$\pair{\delta\mu}{\xi} = J'(\xi) - \mathcal T\Jbold{\xi}$ defines
an element $\delta\mu \in \mathfrak g^* + \nu \mathfrak
g^*_c[[\nu]]$. Using the $G$-equivariance of the quantum momentum
maps and the $G$-invariance of $\mathcal T$ it is evident that
$\delta\mu$ is $G$-invariant. Moreover, it is straightforward to
see that the condition $\mu' =
\mu + \delta\mu$ implies that $\mathcal T$ maps elements of
$I^{\mathcal P}_{\mu,\star}$ to elements of ${I'}^{\mathcal
P}_{\!\!\mu',\star'}$. Using the definition of the products
$\bullred{J,\mu}$ and ${\bullet'}^{J',\mu'}$ on $\hor (\PolFun^G)[[\nu]]$
an easy computation then shows that
$\mathsf T$ as defined in Eq.\ (\ref{mathsfTDefEq}) satisfies the
stated properties and that its inverse is given by
Eq.\ (\ref{mathsfTInvEq}). Statement iv.) is obvious by observing that the
induced mappings on $\PolFunRed[[\nu]]$ can be described as formal
series of differential operators possibly composed with the
pull-back of a diffeomorphism of $T^*\cc{Q}$ preserving
$\PolFunRed$. This latter situation occurs in case $\mathcal T$ and
$\mathcal A$ do not start with $\id$.
Statement v.) is evident from the fact that
we are only concerned with local operators.
\end{PROOF}

Note that, evidently, the above construction of mappings on the
reduced star product algebras is unfortunately not rich enough to
yield all these mappings but only describes under which
circumstances such mappings on the original star product on
$\Cinf{T^*Q}[[\nu]]$ descend to such mappings on them.

Now we want to consider star products which
possess certain additional properties and will show under which
preconditions these properties transfer to the reduced star
products. Let us first recall some notions of special star products
(cf.\ e.g.\ \cite{BorNeuPflWal03,GutRaw03}):
\begin{enumerate}
\item
A star product $\star_{\mbox{\rm\tiny s}}$ resp.\
$\star_{\mbox{\rm\tiny as}}$ on $T^*Q$ is said to be of standard
ordered resp.\ anti-standard ordered type, if for all
$ f \in \Cinf{T^*Q}[[\nu]]$ and $\chi \in \Cinf{Q}[[\nu]]$
\begin{equation}
\pi^* \chi \star_{\mbox{\rm\tiny s}} f = \pi^*\chi \,f \quad
\textrm{resp.}\quad f \star_{\mbox{\rm\tiny as}} \pi^*\chi =
f \,\pi^*\chi.
\end{equation}
\item
A star product is called of Vey type or a natural star product, if
the bidifferential operator describing the star product at order
$r$ in the formal parameter is of order $r$ in each argument.
\end{enumerate}

\begin{LEMMA}
\begin{enumerate}
\item
If $\star_{\mbox{\rm\tiny s}}$ resp.\ $\star_{\mbox{\rm\tiny as}}$
is a star product of standard ordered resp.\ anti-standard ordered
type on $(T^*Q,\omega_{B_0})$, then for every quantum momentum map
$J$ and every possible choice of a momentum value $\mu$ the reduced
star product ${\star_{\mbox{\rm\tiny s}}}^{\!\!J,\mu}$ resp.\
${\star_{\mbox{\rm\tiny as}}}^{\!\!\!\!J,\mu}$ on
$(T^*\cc{Q},\omega_{b_0})$ is also of standard ordered resp.\
anti-standard ordered type.
\item
If $\star$ is a star product of Vey type on $(T^*Q,\omega_{B_0})$,
then for every quantum momentum map $J$ and every possible choice
of a momentum value $\mu$ the reduced star product
$\starred{J,\mu}$ on $(T^*\cc{Q},\omega_{b_0})$ is of Vey type as
well.
\end{enumerate}
\end{LEMMA}
\begin{PROOF}
For the proof of i.) we only consider the case of a standard ordered
star product since the proof for anti-standard ordered star
products is completely analogous. First, observe that for $\chi
\in \Cinf{\cc{Q}}$ we have $l(\PolRed{\chi}) =
l(\cc{\pi}^*\chi) = \Pol{p^*\chi}= \pi^*p^*\chi$. Thus, an easy
computation yields
\begin{eqnarray*}
l(\PolRed{\chi}) {\bullet_{\mbox{\rm\tiny s}}}^{\!\!J,\mu}
l(\PolRed{t}) &=& \hor_{\Check j_0 - \mu_0}\left( \frac{\id}{\id -
\nu
\triangle_{\mu,\star_{\mbox{\rm\tiny s}}}}
\left(\pi^*p^*\chi \star_{\mbox{\rm\tiny s}} \Pol{t^\hor}
\right)\right)\\
&=&\hor_{\Check j_0 - \mu_0}\Bigg( \frac{\id}{\id - \nu
\triangle_{\mu,\star_{\mbox{\rm\tiny s}}}}
\Big(\underbrace{\pi^*p^*\chi \,\Pol{t^\hor}}_{=
\Pol{(\chi t)^\hor}= l(\PolRed{\chi t})}\Big)\Bigg) =
l(\PolRed{\chi} \PolRed{t}).
\end{eqnarray*}
This implies that $\cc{\pi}^*\chi {\star_{\mbox{\rm\tiny
s}}}^{\!\!J,\mu} \PolRed{t}= \cc{\pi}^*\chi \,\PolRed{t}$ for all
$t\in \Ginf{\bigvee T\cc{Q}}$, hence
${\star_{\mbox{\rm\tiny s}}}^{\!\!J,\mu}$ is of standard ordered
type. The proof of ii.) consists of a rather lengthy but straightforward
argument counting the order of differentiation in every
order of the formal parameter within the projection $\hor_{\Check j_0
- \mu_0}\circ\frac{\id}{\id
- \nu \triangle_{\mu,\star}}$.
\end{PROOF}

Now recall (cf.\ e.g.\ \cite[Def.\ 3]{Neu02}) that a star product
$\star$ is called Hermitian, if the operation of complex
conjugation $\C$, where we set $\C \nu := -\nu$, is an
anti-automorphism of $\star$. Analogously, $\star$ is said to have
the $\nu$-parity property, if $\Pa := (-1)^{\deg_\nu}$ is an
anti-automorphism of $\star$, where $\deg_\nu:= \nu \partial_\nu$.
Finally, a star product which is Hermitian and has the $\nu$-parity
property is called star product of Weyl type. In contrast to the
above lemma, where the properties of the original star product
transfer to the reduced star product without any further
conditions, the property of being Hermitian and the $\nu$-parity
are not stable with respect to reduction, in general, unless
certain additional conditions on the quantum momentum map $J$ and
the momentum value $\mu$ are satisfied.
\begin{LEMMA}
\begin{enumerate}
\item
Let $\star$ be a Hermitian star product on $(T^*Q,\omega_{B_0})$.
If the relation
\begin{equation}\label{RedHermCondEq}
\C \left(\JP{\xi} - \pair{\mu_+}{\xi}- \frac{\nu}{2}
\tr {\ad(\xi)}\right) = \JP{\xi} - \pair{\mu_+}{\xi} - \frac{\nu}{2}
\tr {\ad(\xi)}
\end{equation}
is satisfied for all $\xi \in \mathfrak g$, then the reduced star product
$\starred{J,\mu}$ is Hermitian as well.
\item
Let $\star$ be a star product on $(T^*Q,\omega_{B_0})$ which has the
$\nu$-parity property. If the equation
\begin{equation}\label{RedPariCondEq}
\Pa \left(\JP{\xi} - \pair{\mu_+}{\xi} - \frac{\nu}{2}
\tr {\ad(\xi)}\right) = \JP{\xi} - \pair{\mu_+}{\xi}-
\frac{\nu}{2} \tr {\ad(\xi)}
\end{equation}
holds true for all $\xi \in \mathfrak g$, then the reduced star product
$\starred{J,\mu}$ also has the $\nu$-parity property.
\item
Let $\star$ be a star product of Weyl type on
$(T^*Q,\omega_{B_0})$. If Eqs.\ (\ref{RedHermCondEq}) and
(\ref{RedPariCondEq}) are satisfied, then the reduced star product
$\starred{J,\mu}$ is of Weyl type, too.
\end{enumerate}
\end{LEMMA}
\begin{PROOF}
To prove i.) we observe first that $\C J$ defines a $G$-equivariant
quantum Hamiltonian since $\star$ is assumed to be Hermitian.
Suppose that Eq.\ (\ref{RedHermCondEq}) holds true. Then it is easy
to verify that $\C (\nu\triangle_{\mu,\star}F) =
\nu\triangle_{\mu,\star}\C F$ for all $F\in \PolFun^G[[\nu]]$.
Since $\hor_{\Check j_0-\mu_0}$ and $l$ obviously commute with $\C$,
this implies that $\starred{J,\mu}$ is Hermitian.
The proof of ii.) is completely analogous to the one of i.)
replacing $\C$ by $\Pa$ and iii.) follows by combination of i.) and ii.).
\end{PROOF}

Now we consider homogeneous star products on
$(T^*Q,\omega_0)$, i.e.\ star products for which $\mathcal H =
\Lie_{\xi_0} + \nu\partial_\nu$ is a derivation. Observe
that for every $\kappa \in [0,1]$ both the star product $\stark$ and
the star product $\starkB$, $B = \nu B_1$ constructed in
Section \ref{StarProdConsSubSec} are of this kind. Therefore,
the following results directly apply to some of the examples we
will discuss in more detail in Sections \ref{InvQMMSubSec}
and \ref{ExaRedSubSec}.
\begin{LEMMA}\label{RedHomoCondLem}
\begin{enumerate}
\item
Let $\star$ denote a homogeneous star product on $(T^*Q,\omega_0)$
and recall that the $G$-equivariant quantum momentum map $J$ has
the form
 $\Jbold{\xi}= \Pol{\xi_\Qind} + \pair{\Tilde \mu_0}{\xi}+ \JP{\xi}$ with
$\Tilde \mu_0\in {\mathfrak g^*}^G$.
Then the reduced star product $\starred{J,\mu}$ on
$(T^*\cc{Q},\cc{\omega}_0)$ turns out to be homogeneous as well, if $J$ and
the momentum value $\mu$ satisfy
\begin{equation}\label{RedHomoCondEq}
\mathcal H \Jbold{\xi} - \Jbold{\xi} =
\nu\partial_\nu\pair{\mu}{\xi}-\pair{\mu}{\xi}
\quad\text{ for all } \xi \in \mathfrak g.
\end{equation}
\item
Under the preconditions of i.), there exists another
$G$-equivariant quantum momentum map $J'$ of particular form
$J'(\xi)= \Pol{\xi_\Qind} + \nu J_1(\xi)
=\Pol{\xi_\Qind} + \nu \pi^* j_1(\xi)$ with $j_1 \in
C^1(\mathfrak g,\Cinf{Q})$ and another momentum value $\mu' = \nu
\mu_1$ such that $\starred{J',\mu'}$ coincides with
$\starred{J,\mu}$.
\end{enumerate}
\end{LEMMA}
\begin{PROOF}
First observe that Eq.\ (\ref{RedHomoCondEq}) implies
$\Tilde\mu_0=\mu_0$, since $\Lie_{\xi_0} \Pol{\eta_\Qind} =
\Pol{\eta_\Qind}$. Therefore $\starred{J,\mu}$ is a star product
with respect to the canonical symplectic form on $T^*\cc{Q}$. Using
(\ref{RedHomoCondEq}) it is easy to show that the mapping $\mathcal
H$, which evidently commutes with every $\Phi_g^*$,
preserves $I_{\mu,\star}^{\mathcal P}$. This implies that the mapping
$\mathsf H: \hor(\PolFun^G)[[\nu]]\to
\hor(\PolFun^G)[[\nu]]$, which is defined by
\[
\mathsf H F := \hor\left( \frac{\id}{\id - \nu
\triangle_{\mu,\star}}\mathcal H F
\right),\quad  F \in \hor(\PolFun^G)[[\nu]],
\]
is a derivation with respect to the product $\bullred{J,\mu}$.
Since according to Eq.\ (\ref{RestAbbDefEq}) the map $\mathrm r^i$,
$1\leq i\leq \dim{(G)}$ is homogeneous of degree $-1$ with respect
to $\xi_0$, an easy computation shows that $\mathcal H$ commutes
with $\nu \triangle_{\mu,\star}$. Thus  $\mathsf H F =
\nu\partial_\nu F + \hor (\Lie_{\xi_0}F)$ holds true. From this
observation it is evident that the composition $l^{-1} \circ
\mathsf H \circ l$, which is a derivation of the star product
$\starred{J,\mu}$, equals $\cc{\mathcal H}= \nu
\partial_\nu + \Lie_{\cc{\xi}_0}$, where $\cc{\xi}_0$ denotes the
canonical Liouville vector field on $(T^*\cc{Q},\cc{\omega}_0)$.
This proves i.). For the proof of ii.) one has to analyse Eq.\
(\ref{RedHomoCondEq}). By $\Jbold{\xi}
\in \PolFun[[\nu]]$ and since every  eigenvalue of
$\Lie_{\xi_0}$ has to be a non-negative integer, one observes that $J$
has the form
$\Jbold{\xi} = \Pol{\xi_\Qind} + \pair{\Tilde
\mu_0}{\xi}+ \nu \pi^*j_1(\xi) + \sum_{r=2}^\infty \nu^r
\pair{\mu_r}{\xi}$. Using $\Tilde \mu_0=\mu_0$ we then obtain
the relation $\Jbold{\xi} - \pair{\mu}{\xi} =
\Pol{\xi_\Qind} + \nu\pi^*j_1(\xi) -
\nu \pair{\mu_1}{\xi}$.
This implies that $\starred{J,\mu}$
equals $\starred{J',\mu'}$ which proves ii.).
\end{PROOF}

\begin{REMARK}\label{MomentValueRem}
Proposition ii.) of the preceding lemma shows in particular
that two reduced star products $\starred{J,\mu}$
and $\starred{J',\mu'}$ coincide, if $\Jbold{\xi} -
\pair{\mu}{\xi} = J'(\xi) -\pair{\mu'}{\xi}$. This follows from
the fact that only the difference of the quantum momentum map and
the momentum value enter the construction of the reduced star
products. In view of this observation it is equivalent to either
fix a quantum momentum map and vary the momentum values or vary the
quantum momentum map and choose the momentum value to be zero. In
the general considerations of the following sections we will
nevertheless treat $J$ and $\mu$ as independent parameters of the
construction, but in concrete examples it will turn out to be
convenient to restrict the consideration to one fixed quantum
momentum map and to vary the momentum values arbitrarily.
\end{REMARK}
\section{Invariant Star Products on $T^*Q$ and
Quantum Moment Maps}
\label{InvQMMSec}
In this section we prepare the grounds for the phase space
reduction of the star products $\stark$ and $\starkB$ defined in
Section \ref{StarProdConsSubSec}. To this end we will provide
conditions on the data entering the construction of $\stark$ and
$\starkB$ which guarantee that the obtained star products are
$G$-invariant. Furthermore, we provide conditions which are
necessary and sufficient for the existence of a $G$-equivariant
quantum momentum map. Note that all results derived in this section
hold for the lifted action of an arbitrary Lie group action on the
base manifold $Q$ and that the assumption made earlier for phase
space reduction, namely that the action is proper and free, is not
needed here. Only assuming that there exists a $G$-invariant
torsion free connection on $Q$ we will show in particular that
every $G$-invariant star product on $(T^*Q,\omega_{B_0})$ is
$G$-equivalent to some star product $\starNuB$. This will actually
turn out to be one of the key results for the computation of the
characteristic class of the reduced star products in Section
\ref{CharClassSec}. In the course of these investigations, we also
obtain a complete classification up to $G$-equivalence of star
products on $(T^*Q,\omega_{B_0})$ which are invariant with respect
to lifted group actions. This turns out to be a slight refinement
of the classification results in \cite{BerBieGut98} for our special
geometric situation.
\subsection{Invariance of $\stark$ and $\starkB$ and their Quantum
Moment Maps}\label{InvQMMSubSec} In order to derive necessary
and sufficient conditions for the $G$-invariance of the star
products $\stark$, we prove the following:
\begin{PROPOSITION}
\label{RepEquiProp}
Let $\phi$ denote a diffeomorphism of $Q$, let
$\Phi=T^*({\phi^{-1}})$ be the lift to the cotangent bundle
and assume that the connection $\nabla$ is
invariant with respect to $\phi$, i.e.\ that
$\phi^*\nabla_X Y =\nabla_{\phi^*X} \phi^*Y$ holds true
for all $X,Y \in \Ginf{TQ}$. Then there exists for every
$\kappa \in [0,1]$ a uniquely determined formal series of
differential operators $\mathcal S_{\kappa,\phi}$ on $\Cinf{T^*Q}$
which starts with $\id$ and commutes with $\mathcal H=
\Lie_{\xi_0}+ \nu \partial_\nu$ such that
\begin{equation}\label{repkapalmostequiEq}
\phi^* \repkap{f} (\phi^{-1})^* = \repkap{\mathcal S_{\kappa,\phi}
\Phi^* f} \quad \text{ for all $f \in \Cinf{T^*Q}[[\nu]]$}.
\end{equation}
In addition $\mathcal S_{\kappa,\phi}$, is explicitly given by
\begin{equation}\label{SkapphiEq}
\mathcal S_{\kappa,\phi} =\exp\left(-
\Fop{\frac{\exp (\kappa\nu \sKov)-\id}{\sKov}\left(\phi^*
\alpha-\alpha\right)}\right).
\end{equation}
Moreover, $\mathcal S_{0,\phi}$ and $\mathcal S_{1,\phi}$ are
automorphisms of $\starNu$ and $\starEi$, respectively.
Furthermore, $\mathcal S_{\kappa,\phi}$ is an automorphism of
$\stark$ for $\kappa\neq 0,1$ if and only if $\sKov(\phi^*\alpha-
\alpha)=0$.
\end{PROPOSITION}
\begin{PROOF}
A straightforward computation using the $\phi$-invariance of $\nabla$
and that $\mathsf F$ satisfies Eq.\ (\ref{FopEquiEq})
yields $\phi^* \repkap{f}(\phi^{-1})^*
= \repkap{\Nk^{\!\!\!-1} \Phi^*\Nk f}$. Next observe that $\Phi^*\Delta_0 f=
\Delta_0 \Phi^* f$ by Eqs.\ (\ref{HorVerEquiEq}).
Using the factorization property of $\Nk$ given in Eq.\ (\ref{NkFactEq}),
this proves that the formal series of differential operators
$\mathcal S_{\kappa,\phi}$ given in (\ref{SkapphiEq}) satisfies
(\ref{repkapalmostequiEq}). The facts
that $\mathcal S_{\kappa,\phi}$ starts with $\id$ and that it
commutes with $\mathcal H$ are obvious from its explicit form. For
the proof of the uniqueness assume that ${\mathcal
S'}_{\!\!\kappa,\phi}$ is a second formal series of differential
operators having the same properties as $\mathcal S_{\kappa,\phi}$.
Restricting to the space of functions polynomial in the
momenta we obtain $\repkap{{\mathcal
S'}_{\!\!\kappa,\phi}F}=\repkap{\mathcal S_{\kappa,\phi}F}$ for all
$F\in\PolFun[[\nu]]$. Now observe that the homogeneity of
${\mathcal S'}_{\!\!\kappa,\phi}$ and $\mathcal S_{\kappa,\phi}$
implies that they map elements $F\in\PolFun[[\nu]]$ to elements of
$\PolFun[[\nu]]$. Since the restriction of $\varrho_\kappa$
to $\PolFun[[\nu]]$ is injective one thus concludes ${\mathcal
S'}_{\!\!\kappa,\phi}F=\mathcal S_{\kappa,\phi}F$. But this implies
${\mathcal S'}_{\!\!\kappa,\phi}=\mathcal S_{\kappa,\phi}$ since
differential operators on $\Cinf{T^*Q}$ are completely determined
by their values on $\PolFun$, proving the uniqueness of $\mathcal
S_{\kappa,\phi}$. In case $\kappa = 0$,
$\mathcal S_{0,\phi}$ is an automorphism of $\starNu$ since it
coincides with $\id$. For $\kappa=1$ and $A = \nu (\phi^*\alpha - \alpha)$
consider the operator $\mathcal A_1$ defined by Eq.\ (\ref{AkappaDefEq}).
Then $\mathcal A_1$ coincides with $\mathcal S_{1,\phi}$
and is an automorphism of $\starEi$ if and only if
$\d (\phi^*\alpha - \alpha)=0$. Now recall that $\d \alpha = - \tr{R}$ and
that $\phi^*R=R$ due to the $\phi$-invariance of $\nabla$. Then
$\phi^*\alpha - \alpha$ is obviously closed, hence $\mathcal S_{1,\phi}$
is an automorphism of $\starEi$. Finally, let us consider the case
$\kappa\neq 0,1$. Assuming that $\sKov(\phi^*\alpha -\alpha)=0$,
the map $\mathcal S_{\kappa,\phi}= \exp\left(-\kappa\nu
\Fop{\phi^*\alpha -\alpha} \right)$ coincides with
$\mathcal A_\kappa$ for $A= \kappa\nu (\phi^*\alpha - \alpha)$,
hence it is an automorphism of $\stark$ due to the closedness of
$\phi^*\alpha - \alpha$. Conversely, let us assume that $\mathcal
S_{\kappa,\phi}$ is an automorphism of $\stark$. Then $\mathcal
S_{\kappa,\phi}(\mathcal A_\kappa)^{-1}$ is again an automorphism
of $\stark$, where the latter is taken for $A= \kappa\nu
(\phi^*\alpha - \alpha)$. Now a straightforward expansion of this
automorphism yields $\mathcal S_{\kappa,\phi}(\mathcal
A_\kappa)^{-1}=\exp\left(\frac{\kappa(\kappa-1)}{2}
\nu^2\Fop{\sKov(\phi^*\alpha-\alpha)}+ O(\nu^3)\right)$.
But since every star product automorphism starting with $\id$ is of
the form $\exp(\nu\mathrm D)$, where $\mathrm D$ is a derivation of
the star product, and since the lowest order of a derivation is
given by a symplectic vector field, one concludes that
$\frac{\kappa(\kappa-1)}{2} \nu^2\Fop{\sKov(\phi^*\alpha-\alpha)}$
must be zero since otherwise it can never define a symplectic
vector field. Together with the injectivity of $\mathsf F$ and the
precondition $\kappa\neq 0,1$ we thus obtain
$\sKov(\phi^*\alpha-\alpha)=0$ and the proposition is proved.
\end{PROOF}

Using the explicit formulas for the $\stark$-left- and
$\stark$-right-multiplication with functions $\pi^*\chi$,
$\chi \in \Cinf{Q}$, we also get the following.
\begin{LEMMA}\label{InvConnInvLem}
Let $\phi$ denote a diffeomorphism of $Q$ and let $\stark$ be
invariant with respect to $\Phi= T^*(\phi^{-1})$, i.e.\
$\Phi^*(f\stark f' ) = \Phi^*f \stark \Phi^* f'$ for all $f,f'\in
\Cinf{T^*Q}[[\nu]]$. Then the connection $\nabla$ is invariant with
respect to $\phi$.
\end{LEMMA}
\begin{PROOF}
Let us consider the case $\kappa=0$ first.
Using Eq.\ (\ref{starkrechtsmultEq}) and property (\ref{FopEquiEq})
for $\mathsf F$, a comparison of the terms of second order in the formal
parameter within the equation $\Phi^*(f \starNu \pi^*\chi) =
\Phi^*f \starNu \pi^*\phi^*\chi$ yields
$\Fop{(\phi^*\sKov) \d \chi'} =
\Fop{\sKov \d \chi'}$ for all $\chi'\in \Cinf{Q}$. Since
$\mathsf F$ is injective, we thus conclude that $(\phi^*\sKov)
\d \chi' = \sKov \d \chi'$. Evaluating
this for local coordinate functions $x^l$ and using $\phi^*\sKov =
\sKov - \d x^i\vee \d x^j \vee i_s(S_\phi(\partial_{x^i},
\partial_{x^j}))$, we thus obtain that the tensor field $S_\phi$ from
Lemma \ref{pullbackConnLem} has to vanish. Hence $\nabla$ has
to be invariant with respect to $\phi$. For $\kappa \neq 0$ we
consider the second order terms of $\Phi^*(\pi^*\chi\stark f ) =
\pi^*\phi^*\chi\stark \Phi^*f$. Then Eq.\ (\ref{starklinksmultEq}) yields
$\Fop{\kappa^2(\phi^*\sKov) \d \chi'} =\Fop{\kappa^2\sKov \d
\chi'}$ for all $\chi'\in \Cinf{Q}$. Since $\kappa\neq 0$, we may
conclude as above that this implies $\phi$-invariance of $\nabla$.
\end{PROOF}

Combining the results of Proposition \ref{RepEquiProp} and Lemma
\ref{InvConnInvLem} we get:
\begin{THEOREM}\label{starkInvTheo}
\begin{enumerate}
\item
The star products $\starNu$ and $\starEi$ are $G$-invariant if and
only if the connection $\nabla$ is $G$-invariant.
\item
For $\kappa\neq 0,1$ the star product $\stark$ is $G$-invariant if
and only if the connection $\nabla$ is $G$-invariant and
\begin{equation}
\sKov (\phi_g^*\alpha -\alpha ) =0\quad \textrm{for all }
g \in G.
\end{equation}
\item
In either case we have
\begin{equation}
\label{RepkappaGeq}
\phi_g^* \repkap{f} \phi_{g^{-1}}^* = \repkap{\mathcal
S_{\kappa,\phi_g}\Phi_g^*f},
\end{equation}
where according to Eq.\ (\ref{SkapphiEq}) the automorphism
$\mathcal S_{\kappa,\phi_g}$ of $\stark$ is given by
\begin{equation}\label{SkapPhigEq}
\mathcal S_{\kappa,\phi_g} = \exp\left(-
\Fop{\frac{\exp (\kappa\nu \sKov)-\id}{\sKov}\left(\phi_g^*
\alpha-\alpha\right)}\right).
\end{equation}
More explicitly, this means $\mathcal S_{0,\phi_g}= \id$, $\mathcal
S_{\kappa,\phi_g}=
\exp\left(-\Fop{\kappa\nu \left(\phi_g^*\alpha -\alpha\right)}
\right)$ for $\kappa\neq 0,1$, and $\mathcal S_{1,\phi_g}
= \exp\left(- \Fop{\frac{\exp (\nu \sKov)-\id}{\sKov}\left(\phi_g^*
\alpha-\alpha\right)}\right)$.
\end{enumerate}
\end{THEOREM}
\begin{PROOF}
Let us assume that $\nabla$ is $G$-invariant. Then apply
(\ref{repkapalmostequiEq}) with $\phi = \phi_g$ and use the
representation property of $\varrho_\kappa$ to check
\[
\repkap{\mathcal S_{\kappa,\phi_g}\Phi_g^* (f \stark f')} = \repkap{
\mathcal S_{\kappa,\phi_g}\Phi_g^* f \stark \mathcal S_{\kappa,
\phi_g} \Phi_g^* f'}, \quad f,f'\in \Cinf{T^*Q}[[\nu]].
\]
Restricting to $F,F' \in \PolFun[[\nu]]$, the injectivity of
$\varrho_\kappa$ and the fact that $\mathcal
S_{\kappa,\phi_g}\Phi_g^*$ preserves $\PolFun[[\nu]]$ imply that
\[
\mathcal S_{\kappa,\phi_g}\Phi_g^*(F \stark F') =
\mathcal S_{\kappa,\phi_g}\Phi_g^*F \stark
\mathcal S_{\kappa,\phi_g}\Phi_g^* F'\quad\textrm{for all }F,F'
\in\PolFun[[\nu]].
\]
Since $\stark$ is described by bidifferential operators and since
these are completely determined by their values on
$\PolFun[[\nu]]$, this yields that $\mathcal
S_{\kappa,\phi_g}\Phi_g^*$ is an automorphism of $\stark$. For
$\kappa =0$ resp.\ $\kappa=1$ we already know by Proposition
\ref{RepEquiProp} and the $G$-invariance of $\nabla$ that $\mathcal
S_{0,\phi_g}$ resp.\ $\mathcal S_{1,\phi_g}$ is an automorphism of
$\starNu$ resp.\ $\starEi$, and so is $\Phi_g^*$. In case
$\kappa\neq 0,1$ the additional condition $\sKov(\phi_g^*\alpha
-\alpha)=0$ is equivalent to $\mathcal S_{\kappa,\phi_g}$ being an
automorphism of $\stark$, but this is equivalent to $\Phi_g^*$
being an automorphism of $\stark$. Conversely, let us assume that
$\stark$ is $G$-invariant. Then Lemma \ref{InvConnInvLem} implies
that $\nabla$ is also $G$-invariant. But now we can use the above
consideration to conclude that for $\kappa\neq 0,1$ the additional
equation $\sKov(\phi_g^*\alpha -\alpha)=0$ must hold. Together this
proves claims i.) and ii.) of the theorem. Assertion iii.) is
obvious by the proof of i.), ii.), and Eq.\ (\ref{SkapphiEq}).
\end{PROOF}

Finally, we consider the star products
$\starkB$. Before we can state the generalization of Theorem
\ref{starkInvTheo} we provide some rather trivial but nevertheless
crucial results:
\begin{LEMMA}\label{starkBlinkrechtsMultLem}
\begin{enumerate}
\item
For every $\kappa\in [0,1]$ and formal series $B\in
Z^2_\dR(Q)[[\nu]]$ of closed two-forms on $Q$ with real $B_0$,
$\starkB$-left-multiplication by $\pi^*\chi$, $\chi\in
\Cinf{Q}[[\nu]]$ coincides with  $\stark$-left-multiplication
by $\pi^*\chi$. Analogously,  $\starkB$-right-multiplication
by $\pi^*\chi$ coincides with $\stark$-right-multiplication
by $\pi^*\chi$.
\item
For every $\kappa \in [0,1]$ a mapping $\mathcal A_k$ of the form given in
Eq.\ (\ref{AkappaDefEq}) is an automorphism of $\stark$, if and
only if it is an automorphism of $\starkB$.
\end{enumerate}
\end{LEMMA}
\begin{PROOF}
By definition of $\starkB$ we have $\pi^*\chi\starkB f|_{T^*O_j}=
\mathcal A^j_\kappa \left( ((\mathcal
A^j_\kappa)^{-1} \pi^*\chi|_{T^*O_j}) \stark ((\mathcal
A^j_\kappa)^{-1} f|_{T^*O_j})\right)$. By the explicit form of
$\mathcal A^j_\kappa$ it is obvious that $(\mathcal
A^j_\kappa)^{-1} \pi^*\chi= \pi^*\chi$ holds true and that $(\mathcal
A^j_\kappa)^{-1}$ commutes with $\Fop{\beta}$ for every $\beta
\in\Ginf{\bigvee T^*Q}[[\nu]]$. Using these observations together with the
expression for $\pi^*\chi\stark f$ given in Eq.\
(\ref{starklinksmultEq}) one obtains $\pi^*\chi\starkB f = \pi^*\chi
\stark f$. The proof for $\starkB$-right-multiplication is
completely analogous. Assertion ii.) is obvious from the fact
that mappings as given in Eq.\ (\ref{AkappaDefEq}) commute
with the local isomorphisms $\mathcal A^j_\kappa$ from $\stark$ to
$\starkB$.
\end{PROOF}

After these preparations we can state one of the main results of
this section:
\begin{THEOREM}\label{starkBinvThm}
\begin{enumerate}
\item
The star products $\starNuB$ and $\starEiB$ are $G$-invariant, if
and only if $\nabla$ is $G$-invariant and $\phi_g^*B =B$ for all $g
\in G$.
\item
For $\kappa\neq 0,1$ the star product $\starkB$ is $G$-invariant, if
and only if $\nabla$ and $B$ are $G$-invariant and
$\sKov(\phi_g^*\alpha-\alpha)=0$ for all $g\in G$.
\end{enumerate}
\end{THEOREM}
\begin{PROOF}
First assume that $\nabla$ and $B$ are $G$-invariant. Additionally,
assume for $\kappa \neq 0,1$ that
$\sKov(\phi_g^* \alpha-\alpha)=0$ for all $g\in G$. Then the star products
$\stark$ are all $G$-invariant by Theorem \ref{starkInvTheo}.
But this implies that
\begin{equation}\label{PhigmitAjkappaconjEq}
\begin{split}
&\mathcal A_\kappa^j \Phi_g^* (\mathcal A_\kappa^j)^{-1}\\ &=
t^*_{-(A^j_0 - \phi_g^*A^j_0)} \exp\left(
-\Fop{\frac{\exp(\kappa\nu\sKov)-\exp((\kappa-1)\nu\sKov)}
{\nu\sKov} (A^j -\phi_g^*A^j) -(A^j_0 - \phi_g^*A^j_0)}
\right)\Phi_g^*
\end{split}
\end{equation}
is a local automorphism of $(\Cinf{T^*O_j}[[\nu]],\starkB)$. Obviously,
the equality $\phi_g^*B=B$ entails the relation $\d (A^j
-\phi_g^*A^j)=0$. By Lemma \ref{starkBlinkrechtsMultLem} ii.) this
implies that
\[
  t^*_{-(A^j_0 - \phi_g^*A^j_0)}
  \exp\left( -\Fop{\frac{\exp(\kappa\nu\sKov)-\exp((\kappa-1)\nu\sKov)}
  {\nu\sKov} (A^j -\phi_g^*A^j) -(A^j_0 - \phi_g^*A^j_0)} \right)
\]
is a local automorphism of
$(\Cinf{T^*O_j}[[\nu]],\starkB)$. But then $\Phi_g^*$ is also an
automorphism of $\starkB$, proving one direction of i.) and ii.).
For the converse statement assume that $\starkB$ is $G$-invariant.
Then Lemma \ref{starkBlinkrechtsMultLem} i.) and Lemma
\ref{InvConnInvLem} imply that $\nabla$ is
$G$-invariant. For the cases $\kappa=0$ and $\kappa =1$ this
implies that $\starNu$ and $\starEi$ are invariant. Together with the
above considerations this entails by Lemma
\ref{starkBlinkrechtsMultLem} ii.) that
\[
t^*_{-(A^j_0 - \phi_g^*A^j_0)} \exp\left( -\Fop{\frac{\id-
\exp(-\nu\sKov)} {\nu\sKov} (A^j -\phi_g^*A^j)
-(A^j_0 - \phi_g^*A^j_0)}\right)
\]
and
\[
t^*_{-(A^j_0 -
\phi_g^*A^j_0)} \exp\left( -\Fop{\frac{\exp(\nu\sKov)-
\id} {\nu\sKov} (A^j -\phi_g^*A^j) -(A^j_0 - \phi_g^*A^j_0)}
\right)
\]
define local automorphisms of $(\Cinf{T^*O_j}[[\nu]],\starNu)$
and $(\Cinf{T^*O_j}[[\nu]],\starEi)$, respectively. But this
implies that $A^j -\phi_g^*A^j$ is closed, hence $B$ is
$G$-invariant. For the case $\kappa \neq 0,1$ we need a more
detailed argument. By invariance of the connection the mapping
$\mathcal S_{\kappa,\phi_g}\Phi_g^*$ is an automorphism of
$\stark$, hence $\mathcal A^j_\kappa \mathcal
S_{\kappa,\phi_g}\Phi_g^* (\mathcal A^j_\kappa)^{-1}$ is a local
automorphism of $\starkB$. Since by assumption $\starkB$ is
$G$-invariant, this yields that
\begin{equation*}
\begin{split}
t^*_{-(A^j_0 - \phi_g^*A^j_0)}&
\exp\left( -\Fop{\frac{\exp(\kappa\nu\sKov)-
\exp((\kappa-1)\nu\sKov)} {\nu\sKov} (A^j -\phi_g^*A^j)
-(A^j_0 - \phi_g^*A^j_0)}\right)\times\\&\times
\exp\left( -\Fop{\frac{\exp(\kappa\nu\sKov)-\id}{\sKov}
(\phi_g^*\alpha -\alpha)}\right)
\end{split}
\end{equation*}
is a local automorphism of $\starkB$. Considering the order zero part
in the formal parameter this implies that $t_{-(A^j_0 - \phi_g^*A^j_0)}$ has
to be a local symplectomorphism with respect to $\omega_{B_0}$.
Thus $A^j_0 - \phi_g^*A^j_0$ is closed.
Factorizing the local automorphism corresponding to $A^j_0 -
\phi_g^*A^j_0$ we now obtain that
\[
\exp\left(-\Fop{\frac{\exp(
\kappa\nu\sKov)- \exp((\kappa-1)\nu\sKov)}{\nu\sKov}
(A^j-A^j_0 -\phi_g^*(A^j-A^j_0)) + \frac{\exp(\kappa\nu
\sKov)-\id}{\sKov} (\phi_g^*\alpha -\alpha)}\right)
\]
is a local automorphism of $\starkB$. In order one of $\nu$ this
means that $- \Fop{A^j_1 - \phi_g^*A^j_1 + \kappa (\phi_g^*
\alpha-\alpha)}$ defines a symplectic vector field with respect
to $\omega_{B_0}$. Therefore $A^j_1 - \phi_g^*A^j_1 +
\kappa (\phi_g^*
\alpha-\alpha)$ is closed. But since $\nabla$ is invariant, $\phi_g^*
\alpha-\alpha$ is closed and so is $A^j_1 - \phi_g^*A^j_1$. Factorizing
again the local automorphism corresponding to the closed
one-form $\nu(A^j_1 - \phi_g^*A^j_1 + \kappa (\phi_g^*
\alpha-\alpha))$, we end up with another local automorphism of
$\starkB$. To lowest order in the formal parameter the exponent
of this automorphism is given by $\Fop{- A^j_2 + \phi_g^*A^j_2 +
\frac{\kappa(\kappa-1)}{2} \sKov(\phi_g^*
\alpha-\alpha)}$. This term again has to define a symplectic vector
field. For $\kappa\neq 0,1$ this is only possible, if
$\sKov(\phi_g^* \alpha-\alpha)=0$, hence $\mathcal
S_{\kappa,\phi_g}$ is an automorphism of $\stark$ and $\stark$ is
$G$-invariant. This means that the mapping $\mathcal A_\kappa^j
\Phi_g^* (\mathcal A_\kappa^j)^{-1}\Phi_{g^{-1}}^*$ from
Eq.\ (\ref{PhigmitAjkappaconjEq}) is a local automorphism of
$\starkB$. Now this local automorphism is an operator of form
$\mathcal A_\kappa$ as given in Eq.\ (\ref{AkappaDefEq}) which in turn is
an automorphism of $\stark$ (and hence of $\starkB$ by Lemma
\ref{starkBlinkrechtsMultLem} ii.)), if and only if $A^j
-\phi_g^*A^j$ is closed. But this implies finally that $\phi_g^*B - B = \d
(\phi_g^*A^j  -A^j)=0$, i.e.\ that $B$ is $G$-invariant.
\end{PROOF}

After having investigated the invariance of the star products
$\stark$ and $\starkB$, we will next consider $G$-equivariant
quantum momentum maps for these star products. More precisely, we
now state one of the key results on phase space reduction of the
star products $\stark$. To this end we henceforth assume $\stark$
to be $G$-invariant. Moreover, we assume the connection $\nabla$ to
be $G$-invariant and that for $\kappa\neq 0,1$ the relation
$\sKov(\phi_g^*\alpha -\alpha)=0$ holds true for all $g\in G$.
\begin{PROPOSITION}\label{starkstrongInvProp}
For every $\kappa\in [0,1]$, the $G$-invariant star product $\stark$
on $(T^*Q,\omega_0)$ is strongly $G$-invariant, i.e.\ the map $J \in
C^1(\mathfrak g,\Cinf{T^*Q})$ with $\Jbold{\xi}:=\JN{\xi}=
\Pol{\xi_\Qind}$ is a $G$-equivariant quantum momentum map for
the lifted Lie group action. In case $J'$, where $J'_0$ is assumed
to be real, is another $G$-equivariant quantum momentum map  for
the same action, then $J'$ is given by $J'(\xi) = \JN{\xi} +
\pair{\Tilde \mu}{\xi}$ with $\Tilde\mu \in {\mathfrak g^*}^G +
\nu {\mathfrak g^*_c}^G[[\nu]]$.
\end{PROPOSITION}
\begin{PROOF}
Put $g= \exp(t\xi)$ in Eq.\ (\ref{RepkappaGeq}),
differentiate the resulting equation with respect to $t$ and evaluate at
$t=0$. Then one obtains $[\Lie_{\xi_\Qind},\repkap{f}]=
\repkap{ -\Fop{\frac{\exp(\kappa\nu\sKov)-\id}{\sKov}\Lie_{\xi_\Qind}
\alpha}f + \Lie_{\xi_\QKotind} f}$. From
the definition of $\varrho_\kappa$ one concludes immediately
that $\repkap{\Pol{\xi_\Qind}} =
-
\nu \Lie_{\xi_\Qind} - \kappa\nu \repkap{\pi^*(\mathsf{div}(\xi_\Qind) +
\alpha(\xi_\Qind))}$, where $\mathsf{div}(\xi_\Qind) = \tr{Y\mapsto
\nabla_Y \xi_\Qind}$ denotes the covariant divergence of the vector
field $\xi_\Qind$. Using the representation property of
$\varrho_\kappa$ this implies
\[
\repkap{-\frac{1}{\nu}\ad_{\stark}(\Pol{\xi_\Qind} + \kappa\nu
\pi^*(\mathsf{div}(\xi_\Qind) + \alpha(\xi_\Qind))) f +
\Fop{\frac{\exp(\kappa\nu\sKov)-\id}{\sKov}\Lie_{\xi_\Qind}
\alpha} f}  = \repkap{\Lie_{\xi_\QKotind} f}.
\]
At this point we need a little technical result.
\begin{SUBLEMMA}
If the star product $\stark$ is $G$-invariant, then the following equality
holds true:
\begin{equation}\label{FopinnerEq}
\kappa \ad_{\stark}(\pi^*(\mathsf{div}(\xi_\Qind) + \alpha(\xi_\Qind)))f =
\Fop{\frac{\exp(\kappa\nu\sKov)-\id}{\sKov}\Lie_{\xi_\Qind}
\alpha}f \quad \text{ for all $f\in \Cinf{T^*Q}$}.
\end{equation}
\end{SUBLEMMA}
\begin{INNERPROOF}
For $\kappa=0$ there is nothing to show, since both sides coincide
with $0$. For $\kappa=1$ we obtain from (\ref{starklinksmultEq})
and (\ref{starkrechtsmultEq}) that
$\ad_{\starEi}(\pi^*(\mathsf{div}(\xi_\Qind) + \alpha(\xi_\Qind)))
=
\Fop{\frac{\exp(\nu \sKov)-\id}{\sKov} \d
(\mathsf{div}(\xi_\Qind) + \alpha(\xi_\Qind))}$.
Now observe that
 $(\Lie_{\xi_\Qind}\nabla)_Y Z =
\Lie_{\xi_\Qind}\nabla_Y Z - \nabla_{\Lie_{\xi_\Qind}Y} Z -
\nabla_Y  \Lie_{\xi_\Qind} Z=0$
by invariance of $\nabla$.
A straightforward computation then shows
$-\tr{R}(\xi_\Qind, Y)= \d \alpha (\xi_\Qind ,Y)=
\d(\mathsf{div}(\xi_\Qind))(Y)$, but from this equation and Cartan's formula
Eq.\ (\ref{FopinnerEq}) is immediate in case $\kappa=1$.
For $\kappa\neq 0,1$ the additional relation $\sKov(\phi_g^* \alpha
-\alpha)=0$ obviously implies $\sKov \Lie_{\xi_\Qind}\alpha =0$,
hence $\Fop{\frac{\exp(\kappa\nu\sKov)-\id}{\sKov}\Lie_{\xi_\Qind}
\alpha}=\kappa\nu\Fop{\Lie_{\xi_\Qind}\alpha} =
\kappa\Fop{\frac{\exp(\kappa\nu\sKov) -\exp((\kappa-1)
\nu \sKov)}{\sKov} \Lie_{\xi_\Qind} \alpha}$. By the argument above
$\Lie_{\xi_\Qind}\alpha = \d (\mathsf{div}(\xi_\Qind) +
\alpha(\xi_\Qind))$, hence
$\Fop{\frac{\exp(\kappa\nu\sKov)-\id}{\sKov}\Lie_{\xi_\Qind}
\alpha}$ equals $\kappa
\ad_{\stark} (\pi^*(\mathsf{div}(\xi_\Qind) + \alpha(\xi_\Qind)))$ due to
the explicit formulas (\ref{starklinksmultEq}) and
(\ref{starkrechtsmultEq}).
\end{INNERPROOF}

By the sublemma we now obtain
$\repkap{-\frac{1}{\nu}\ad_{\stark}(\Pol{\xi_\Qind})f} =
\repkap{\Lie_{\xi_\QKotind}f}$. Recall that
$\varrho_\kappa$ restricted to $\PolFun[[\nu]]$ is injective
and observe that $\Lie_{\xi_\QKotind}$ and $\ad_{\stark}
(\Pol{\xi_\Qind})$ preserve $\PolFun[[\nu]]$. Then we obtain
\[
-\frac{1}{\nu}\ad_{\stark}(\Pol{\xi_\Qind})F = \Lie_{\xi_\QKotind}F
\quad\textrm{for all }F \in \PolFun[[\nu]].
\]
Since $\stark$ is described by bidifferential operators,
this implies $- \Lie_{\xi_\QKotind} =
\frac{1}{\nu}\ad_{\stark}(\Pol{\xi_\Qind})$, hence $\Jbold{\xi} =
\JN{\xi}=\Pol{\xi_\Qind}$ is a quantum Hamiltonian which
is known to be $G$-equivariant. The claim about the ambiguity of
the quantum momentum map is a general result which holds over
arbitrary symplectic manifolds (cf.\ \cite[Prop.\ 6.3]{Xu98}).
\end{PROOF}

For the study of reduction of the invariant star products $\stark$
the $G$-equivariant quantum momentum map given by
$\Jbold{\xi}=\Pol{\xi_\Qind}$ appears to be the preferred one,
since it coincides with the (canonical) $G$-equivariant classical
momentum map. In the sequel we therefore will mostly work with this
momentum map in order to compute the reduced products of $\stark$.

Now we consider the star products $\starkB$ on
$(T^*Q,\omega_{B_0})$ and will derive necessary and sufficient
conditions for the existence of $G$-equivariant quantum momentum
maps for $G$-invariant star products of form $\starkB$. In view of
Theorem \ref{starkBinvThm} we assume additionally to the conditions
which guarantee $\stark$ to be $G$-invariant that $B$ is an element
of $Z^2_\dR (Q)^G[[\nu]]$.
\begin{PROPOSITION}\label{starkBQMMExProp}
Suppose that the star product $\starkB$ on $(T^*Q,\omega_{B_0})$ is
$G$-invariant. Then there is a $G$-equivariant quantum momentum map
for the $G$-action under consideration, if and only if there is an
element $j\in C^1(\mathfrak g, \Cinf{Q})[[\nu]]$ with real-valued
$j_0$ such that
\begin{equation}\label{starkBQMMExEq}
\d \jbold{\xi} = i_{\xi_\Qind}B\quad\textrm{ and }\quad
\phi_g^* \jbold{\xi} = \jbold{\Ad(g^{-1})\xi}
\quad\textrm{for all $\xi\in\mathfrak g$, $g \in G$}.
\end{equation}
In this case we particularly have
\begin{equation}\label{starkBQMMExConsEq}
\jbold{[\xi,\eta]}= B(\xi_\Qind,\eta_\Qind)\quad\textrm{for all }
\xi,\eta \in\mathfrak g.
\end{equation}
Moreover, the map $J\in C^1(\mathfrak g,\Cinf{T^*Q})[[\nu]]$ given
by $\Jbold{\xi}:= \Pol{\xi_\Qind}+ \pi^*\jbold{\xi}$ defines a
$G$-equivariant quantum momentum map, which is unique up to
elements of ${\mathfrak g^*}^G+ \nu{\mathfrak g_c^*}^G[[\nu]]$.
\end{PROPOSITION}
\begin{PROOF}
Consider the following equation over $T^*O_j$:
\[
\Phi_g^* \mathcal A_\kappa^j \Phi_{g^{-1}}^* = t_{-\phi_g^*A^j_0}^*
\exp\left(- \Fop{\frac{\exp(\kappa\nu\sKov)-
\exp((\kappa-1)\nu\sKov)}{\nu\sKov}\phi_g^*A^j - \phi_g^*A^j_0}
\right).
\]
By differentiation at the
neutral element of $G$ we obtain
\[
\Lie_{\xi_\QKotind} = \mathcal
A^j_\kappa
\Lie_{\xi_\QKotind}(\mathcal A^j_\kappa)^{-1} - \Fop{\frac{\exp(\kappa\nu
\sKov)-
\exp((\kappa-1)\nu\sKov)}{\nu\sKov}\Lie_{\xi_\Qind}A^j}.
\]
Now the invariance of $\starkB$ implies that $\stark$ is invariant as well,
hence we conclude from Proposition \ref{starkstrongInvProp} that
$\Lie_{\xi_\QKotind}=
-\frac{1}{\nu}\ad_{\stark}(\Pol{\xi_\Qind})$. By a direct
computation using that $\mathcal A^j_\kappa$ is a local
homomorphism from $\stark$ to $\starkB$ and that $\Pol{\xi_\Qind}$
is a polynomial in the momenta of degree one this implies
that $\mathcal
A^j_\kappa
\Lie_{\xi_\QKotind}(\mathcal A^j_\kappa)^{-1} = -\frac{1}{\nu}
\ad_{\starkB}(\Pol{\xi_\Qind} - \pi^*(A^j(\xi_\Qind)))$. On the other hand
$\Lie_{\xi_\Qind} A^j = i_{\xi_\Qind} B + \d(A^j(\xi_\Qind))$.
By the explicit formula for $\ad_{\starkB} (\pi^*(A^j(\xi_\Qind)))$
we thus obtain
\[
\Lie_{\xi_\QKotind} = -\frac{1}{\nu}\left(\ad_{\starkB}(\Pol{\xi_\Qind})
+ \Fop{\frac{\exp(\kappa\nu \sKov)-
\exp((\kappa-1)\nu\sKov)}{\sKov}i_{\xi_\Qind}B}\right).
\]
From this equation and from the explicit form of the inner
derivations $\ad_{\starkB}(\pi^*\chi)$, $\chi\in \Cinf{Q}[[\nu]]$
it is clear that there is a quantum Hamiltonian for the considered
action, if and only if there is an element $j \in C^1(\mathfrak
g,\Cinf{Q})[[\nu]]$ such that $\d \jbold{\xi}= i_{\xi_\Qind}B$ for
all $\xi \in \mathfrak g$. Observe that the condition necessary for
the solvability of this equation is satisfied, since by invariance
of $\starkB$ we have that $\d i_{\xi_\Qind}B = \Lie_{\xi_\Qind} B
=0$. Evidently, $J$ with $\Jbold{\xi} =
\Pol{\xi_\Qind} + \pi^*\jbold{\xi}$ is additionally $G$-equivariant,
if and only if the second condition in Eq.\ (\ref{starkBQMMExEq})
is satisfied. Finally, differentiating this equation in $g$ at $e$
and using the equality $i_{\xi_\Qind} B = \d \jbold{\xi}$,
it is straightforward to check that Eq.\ (\ref{starkBQMMExConsEq}) holds true.
The claim
about the ambiguity of $J$ is well known, hence the proposition
is proved.
\end{PROOF}

\begin{REMARK}
It is immediate to show that $\starkB$ admits a
quantum Hamiltonian $J$ which additionally satisfies
$\frac{1}{\nu}\ad_{\starkB}(\Jbold{\xi})
\Jbold{\eta}=\Jbold{[\xi,\eta]}$, if and only if
there is an element $j \in C^1(\mathfrak g,\Cinf{Q})[[\nu]]$ such
that $\d \jbold{\xi}= i_{\xi_\Qind} B$ and $B(\xi_\Qind,\eta_\Qind)
=\jbold{[\xi,\eta]}$. Moreover, these conditions determine $j$
up to elements of $\mathfrak g^* + \nu \mathfrak g^*_c[[\nu]]$ which are
invariant with respect to the coadjoint action of $\mathfrak
g$, i.e.\ which vanish on $[\mathfrak g,\mathfrak g]$.
\end{REMARK}
\begin{COROLLARY}
Suppose that the star product $\starkB$ on $(T^*Q,\omega_{B_0})$ is
$G$-invariant and that $J_0$ with $\JN{\xi}= \Pol{\xi_\Qind} +
\pi^*\jN{\xi}$ is a $G$-equivariant classical momentum map.
Then $\starkB$ is strongly $G$-invariant,
if and only if $B_+\in \nu Z^2_\dR(Q)^G[[\nu]]$, where
$B= B_0 + B_+$, is horizontal, i.e.\ if and only if
\begin{equation}\label{starkBStrongInvEq}
i_{\xi_\Qind} B_+ = 0 \quad\textrm{for all }\xi \in \mathfrak g.
\end{equation}
\end{COROLLARY}
\subsection{General Invariant Star Products on $T^*Q$: Relations to
$\starNuB$ and Classification} In this section we consider star
products $\star$ on $(T^*Q,\omega_{B_0})$ which are $G$-invariant
with respect to a lifted group action. Under the assumption that
there is a $G$-invariant torsion free connection $\nabla$ on $Q$ we
will in particular construct for every such $\star$ a
$G$-equivalent star product of form $\starNuB$. Incidently, our
results show that there is a $G$-equivariant quantum momentum map
for an arbitrary $G$-invariant star product in the above sense if
and only if there is a $G$-equivariant quantum momentum map for a
certain star product $\starNuB$. But for these star products we
have already derived criteria for the existence of $G$-equivariant
quantum momentum maps. Thus we obtain necessary and sufficient
conditions for the existence of $G$-equivariant quantum momentum
maps for an arbitrary $G$-invariant star product. Finally, we use
the $G$-equivalence between a $G$-invariant $\star$ and and an
appropriate star product $\starNuB$ to give a classification of
star products on $(T^*Q,\omega_{B_0})$ up to $G$-equivalence.
Actually, our result is a slight refinement of the general
classification results of \cite{BerBieGut98} in the particular case
of a cotangent bundle with lifted $G$-action.

First we need a few results which allow for a comparison of two
different star products $\star_0^{B^{(k)}}$ and $\star_0^{B^{(k+1)}}$ in
case $B^{(k)}=\sum_{l=0}^k \nu^l B_l$ and $B^{(k+1)}=\sum_{l=0}^{k+1}
\nu^l B_l$. Note that the following lemma holds
for arbitrary ordering parameters $\kappa\in [0,1]$ but since we
only need it for $\kappa=0$ we restrict the proof to this particular
case. The proof for general $\kappa$ is an immediate adaption.
\begin{LEMMA}\label{BDifferLem}
For $k\in \mathbb N$ let $B^{(k+1)}=\sum_{l=0}^{k+1} \nu^l B_l$ and
$B^{(k)}=\sum_{l=0}^k \nu^l B_l$ be series of closed two-forms on
$Q$. Then the describing bidifferential operators
$C_r^{\star_0^{B^{(k+1)}}}$ and $C_r^{\star_0^{B^{(k)}}}$ of the
corresponding star products $\star_0^{B^{(k+1)}}$ and
$\star_0^{B^{(k)}}$ on $(T^*Q,\omega_{B_0})$ satisfy
\begin{eqnarray}\label{CBbiskplusersteOrdEq}
\!\!\!\!\!C_r^{\star_0^{B^{(k+1)}}} \!&=& \!C_r^{\star_0^{B^{(k)}}}\quad
\textrm{ for } r = 0,\ldots,k+1\quad\textrm{ and}\\
\label{CBkpluszweiteOrdEq}
\!\!\!\!\!C_{k+2}^{\star_0^{B^{(k+1)}}}(f,f')-
C_{k+2}^{\star_0^{B^{(k+1)}}}(f',f) \!&=& \!
C_{k+2}^{\star_0^{B^{(k)}}}(f,f')-
C_{k+2}^{\star_0^{B^{(k)}}}(f',f) - (\pi^*B_{k+1})(X_f^{B_0},
X_{f'}^{B_0}), \quad
\end{eqnarray}
where $X_f^{B_0}$ denotes the Hamiltonian vector field of $f\in
\Cinf{T^*Q}$ with respect to the symplectic form $\omega_{B_0}$.
\end{LEMMA}
\begin{PROOF}
Let $O$ be an element of a good open cover of $Q$. Over $O$
consider a local potential $A^{(k+1)} = A^{(k)} + \nu^{k+1}A_{k+1}$
of $B^{(k+1)}$, where $A^{(k)}$ is a local potential of $B^{(k)}$
and $A_{k+1}$ is a local potential of $B_{k+1}$. From the very
definition of the star products corresponding to $B^{(k+1)}$ and
$B^{(k)}$ one obtains that $\mathcal S_{k+1}:= \exp \left(-
\Fop{\frac{\id - \exp(-\nu \sKov)}{\nu \sKov}\nu^{k+1} A_{k+1}}
\right)$ defines a local equivalence from $(\Cinf{T^*O}[[\nu]],
\star_0^{B^{(k)}})$ to $(\Cinf{T^*O}[[\nu]],
\star_0^{B^{(k+1)}})$. Expanding $\mathcal S_{k+1}= \id -
\nu^{k+1}\Fop{A_{k+1}} + O(\nu^{k+2})$ and expanding the products
$\star_0^{B^{(k)}}$ and $\star_0^{B^{(k+1)}}$ one obtains
Eqs.\ (\ref{CBbiskplusersteOrdEq}) and (\ref{CBkpluszweiteOrdEq}) by
an immediate computation.
\end{PROOF}

The following results (which are essentially due to
Lichnerowicz \cite[Lemma 1 and 2]{Lic80} and
Bertelson et al.\
\cite[Prop.\ 2.1]{BerBieGut98}) will turn out to be crucial
for our further investigations.
\begin{LEMMA}\label{InvPotLem}
\begin{enumerate}
\item
Suppose that there is a $G$-invariant torsion free connection
$\nabla$ on $Q$. Then every $G$-invariant
differential $\Cinf{T^*Q}$-Hochschild $p$-coboundary $C$ ($p\geq 1$)
which vanishes on constants
is the coboundary of a $G$-invariant differential $p-1$-cochain $c$
vanishing on constants. In case $C(F_1, \ldots, F_p)\in \PolFun$
for all $F_1,\ldots,F_p\in \PolFun$ one can additionally achieve
that $c(F_1, \ldots, F_{p-1})\in \PolFun$ for all
$F_1,\ldots,F_{p-1}\in
\PolFun$.
\item
For every closed $G$-invariant  $p$-form $\Omega$ ($p\geq 1$) on $T^*Q$
there exists a $G$-invariant closed $p$-form $\beta$ on $Q$
and a $G$-invariant $p-1$-form $\Xi$ on $T^*Q$ such that
$i^*\Xi=0$ and
\begin{equation}\label{KotFormHomoBaseFormEq}
\Omega = \d \Xi + \pi^*\beta.
\end{equation}
If $\Omega(X_{F_1}^{B_0},\ldots,X_{F_p}^{B_0})\in \PolFun$ for all
$F_1,\ldots,F_p\in \PolFun$, then $\Xi$ can be chosen such that
$\Xi(X_{F_1}^{B_0},\ldots,X_{F_{p-1}}^{B_0})\in \PolFun$ for all
$F_1,\ldots,F_{p-1}\in \PolFun$.
\end{enumerate}
\end{LEMMA}
\begin{PROOF}
For the proof of i.) recall from \cite[Def.\ 4]{BorNeuWal98} that
every torsion free connection $\nabla$ on $Q$ defines a
torsion free connection $\nabla^\QKotind$ on $T^*Q$ which is
$G$-invariant if the original connection $\nabla$ is invariant.
Moreover, having chosen a torsion free connection $\nabla^\Mind$ on
a manifold $M$ it is well known that every $p$-cochain $C$ on
$\Cinf{M}$ can be uniquely written as $C(f_1,\ldots,f_p)=
C^{I_1;\ldots;I_p}\left(\left(\sKov^\Mind\right)^{|I_1|}
f_1\right)(\partial_{y^{I_1}})\ldots \left(\left(
\sKov^\Mind\right)^{|I_p|} f_p\right)(\partial_{y^{I_p}})$, where
$I_1,\ldots,I_p$ denote multi-indices and the $C^{I_1;\ldots;I_p}$
are components of tensor fields in
$\Ginf{\bigvee^{|I_1|}TM\otimes\ldots\otimes\bigvee^{|I_p|}TM}$ with respect
to local coordinates $y^1,\ldots,y^m$ of $M$. In case $C$ is a
coboundary one can explicitly build a $p-1$-cochain $c$ such that
$\deltaH c =C$
($\deltaH$ denotes the Hochschild differential), where the
tensor fields defining $c$ are given as combinations of those of
$C$. But this implies that in case $C$ and the connection are
invariant $c$ is also invariant (cf.\ \cite[Remark
2.1]{BerBieGut98}) proving the first part of i.). For the proof of
the second part of i.) one first has to observe that the covariant
derivative with respect to the above connection preserves
$\PolFun$ since the Christoffel symbols of $\nabla^\QKotind$ in a
local bundle chart are polynomials in the momenta. Together with
the assumption that $C$ preserves $\PolFun$, which implies that
the components of the corresponding tensor fields are polynomials
in the momenta, this entails that $c$ also preserves $\PolFun$. For
the proof of ii.) consider the closed $p$-form $\beta:= i^* \Omega$
on $Q$ which is evidently $G$-invariant by
$G$-invariance of $\Omega$. Therefore the closed two-form $\Omega -
\pi^* \beta$ is $G$-invariant and $i^*(\Omega -
\pi^* \beta)=0$. Now, consider the homotopy $H : \mathbb R
\times T^*Q \to T^*Q$, $(t,\zeta_x) \mapsto H(t,\zeta_x) := t \zeta_x$.
By means of this homotopy one can explicitly define a $p-1$-form $\Xi$ by
$\Xi(X_1,\ldots, X_{p-1}):=\int_0^1 (H^*(\Omega - \pi^*
\beta))(\partial_t,X_1,\ldots,X_{p-1})\d t$. This
$\Xi$ satisfies $\d \Xi = \Omega - \pi^* \beta$ by the classical
proof of Poincar\'{e}'s lemma. Due to the compatibility of the
above homotopy with the $\Phi_g$ the thus defined $p-1$-form is
$G$-invariant. Since $H(t,i(x))= i(x)$ we also have $i^*\Xi
=0$. From the explicit shape of $\Xi$ it is also obvious that
$\Xi(X_{F_1}^{B_0},\ldots,X_{F_{p-1}}^{B_0})\in \PolFun$ for all
$F_1,\ldots,F_{p-1}\in \PolFun$ in case
$\Omega(X_{F_1}^{B_0},\ldots,X_{F_p}^{B_0})\in \PolFun$ for all
$F_1,\ldots,F_p\in \PolFun$.
\end{PROOF}

Using these technical preparations we can adapt the proof of
\cite[Prop.\ 4.1]{BerBieGut98} to the present situation and obtain:
\begin{PROPOSITION}\label{GEquivaundBExiProp}
Suppose that there is a $G$-invariant torsion free connection
$\nabla$ on $Q$. Then we have:
\begin{enumerate}
\item
For every star product $\star$ on $(T^*Q,\omega_{B_0})$ which is
invariant with respect to the lifted action of a $G$-action on $Q$
there is a formal series $\frac{1}{\nu}B_+ = \sum_{l=1}^\infty
\nu^{l-1}B_l\in Z^2_\dR(Q)^G[[\nu]]$ of $G$-invariant closed
two-forms on $Q$ and a $G$-equivalence transformation $\mathcal T$
from $\star$ to the  $G$-invariant star product
$\starNuB$, where $B = B_0 + B_+$.
\item
In case $\PolFun[[\nu]]$ is a $\star$-subalgebra
one can find $B_+$ and a $G$-equivalence transformation $\mathcal T$
from $\star$ to $\starNuB$ as in i.) such that $\mathcal T F
\in \PolFun[[\nu]]$ for all $F\in \PolFun[[\nu]]$.
\end{enumerate}
\end{PROPOSITION}
\begin{PROOF}
Let $\star$ be an arbitrary $G$-invariant star product on
$(T^*Q,\omega_{B_0})$ as above and consider the $G$-invariant star
product $\star_0^{B^{(0)}}$, where $B^{(0)}=B_0$. In order zero of
the formal parameter, the star products $\star$ and
$\star_0^{B^{(0)}}$ trivially coincide.
Since both are star products with respect to
$\omega_{B_0}$ the anti-symmetric part of $C_1^{\star_0^{B^{(0)}}}-
C_1^\star$ vanishes. Therefore Lemma \ref{InvPotLem} i.) implies
that there is a $G$-invariant $1$-cochain $c_1$ with
$C_1^{\star_0^{B^{(0)}}}- C_1^\star = \deltaH c_1$. Now put
$T^{(0)} := \id - \nu c_1$, which is clearly $G$-invariant.
Define another $G$-invariant star product by
$\star^{(0)} := T^{(0)}\star$, i.e.\ let $T^{(0)}(f
\star f') = (T^{(0)} f) \star^{(0)} (T^{(0)}f')$. Then an easy
computation shows that $C_r^{\star_0^{B^{(0)}}} =
C_r^{\star^{(0)}}$ for $r=0,1$. By associativity of these
two star products one obtains that the anti-symmetric part of
$C_2^{\star_0^{B^{(0)}}} - C_2^{\star^{(0)}}$ defines a
$G$-invariant closed two-form $\Omega_1$ on $T^*Q$ via
\[
C_2^{\star_0^{B^{(0)}}}(f,f')- C_2^{\star_0^{B^{(0)}}}(f',f) -
C_2^{\star^{(0)}}(f,f') + C_2^{\star^{(0)}}(f',f) =
\Omega_1\left(X_f^{B_0}, X_{f'}^{B_0}\right).
\]
Again from Lemma \ref{InvPotLem} ii.) we get that $B_1:= i^*
\Omega_1$ is a $G$-invariant closed two-form on $Q$ and that there is a
$G$-invariant one-form $\Xi_1$ on $T^*Q$ such that $\Omega_1 = \d
\Xi_1 + \pi^* B_1$. Then we consider the $G$-invariant star product
$\star_0^{B^{(1)}}$ with $B^{(1)}= B^{(0)} + \nu B_1$
and the $G$-invariant star product
$\widetilde{\star}^{(0)}:=\widetilde{T}^{(0)}\star^{(0)}$, where
$\widetilde{T}^{(0)} f:= f + \nu \Xi_1(X_f^{B_0})$.
According to Eq.\ (\ref{CBbiskplusersteOrdEq}) we have
$C_r^{\star_0^{B^{(1)}}} = C_r^{\star_0^{B^{(0)}}}$ for $r=0,1$.
Now it is straightforward to check that
$C_r^{\widetilde{\star}^{(0)}}=C_r^{\star_0^{B^{(1)}}}$ for $r=0,1$
and that the anti-symmetric part of
$C_2^{\star_0^{B^{(1)}}}-C_2^{\widetilde{\star}^{(0)}}$ vanishes
due to the definition of $\widetilde{\star}^{(0)}$ and Eq.\
(\ref{CBkpluszweiteOrdEq}). But then Lemma \ref{InvPotLem}
i.) yields
$C_2^{\star_0^{B^{(1)}}}-C_2^{\widetilde{\star}^{(0)}}=
\deltaH c_2$ with a $G$-invariant $1$-cochain $c_2$. Putting
$\star^{(1)}:= T^{(1)} \widetilde{\star}^{(0)}$ with
$T^{(1)}:= \id - \nu^2 c_2$ we then obtain $C_r^{\star_0^{B^{(1)}}}=
C_r^{\star^{(1)}}$ for $r=0,1,2$. Proceeding inductively we
thus can find $G$-invariant operators $T^{(l)} = \id -
\nu^{l+1} c_{l+1}$ for $l=0,\ldots, k$, $G$-invariant
$\widetilde{T}^{(m)}$ for $m=0,\ldots, k-1$ with
$\widetilde{T}^{(m)}f = f + \nu^{m+1} \Xi_{m+1}(X_f^{B_0})$ and
$G$-invariant closed two-forms $B_1,\ldots, B_k$
on $Q$  such that
$C_r^{\star_0^{B^{(k)}}} = C_r^{\star^{(k)}}$ for $r=0,\ldots, k+1$.
Hereby, $B^{(k)}= \sum_{l=0}^k\nu^l B_l$ and
$\star^{(k)} = \mathcal T^{(k)} \star$ with $\mathcal T^{(k)} :=
T^{(k)} \widetilde{T}^{(k-1)} T^{(k-1)}\ldots \widetilde{T}^{(0)}
T^{(0)}$. For $k\to \infty$ we thus obtain a well-defined
$G$-invariant formal series of differential operators $\mathcal T:=
\mathcal T^{(\infty)}$ and a formal series $B=
\sum_{l=0}^\infty\nu^l B_l$ of $G$-invariant closed two-forms on
$Q$ such that the $G$-invariant star product $\mathcal T \star$
coincides with $\starNuB$. Hence $\mathcal T$ is a
$G$-equivalence from $\star$ to $\starNuB$, and $B_+$ is given by
$B_+= B - B_0$ proving i.). For the proof of ii.) one just has to
observe that $\PolFun[[\nu]]$ is a $\star_0^{B^{(0)}}$-subalgebra
and that by assumption $\PolFun[[\nu]]$ is a $\star$-subalgebra.
Using Lemma \ref{InvPotLem} this implies that in every step of the
above construction one can achieve that $T^{(l)}$ and
$\widetilde{T}^{(m)}$ map elements of $\PolFun[[\nu]]$ to elements
of $\PolFun[[\nu]]$. To verify this check  by induction
that $\PolFun[[\nu]]$ is both a $\star^{(k)}$-subalgebra as
well as a $\widetilde{\star}^{(k)}$-subalgebra since
$\PolFun[[\nu]]$ is a $\star_0^{B^{(l)}}$-subalgebra for all
occuring $B^{(l)}$.
\end{PROOF}

Note that even after having fixed a $G$-invariant torsion free connection
on $Q$ one cannot use the construction of $\frac{1}{\nu}B_+$ in the
proof of the preceding proposition to define a map from the
space of $G$-invariant star products on $(T^*Q,\omega_{B_0})$ to
the space of formal series of closed $G$-invariant two-forms on
$Q$. This fact is caused by the freedom of choice in the
equivalence transformations which in fact affects the explicit
form of $\frac{1}{\nu}B_+$. For instance, in the definition
of $T^{(0)}$ we could have replaced $c_1$ by $c_1 + \Lie_X$, where
$X$ is a $G$-invariant vector field on $T^*Q$. An easy
computation then shows that this gives rise to a modification of $B_1$
by the additional term $-\d i^*(i_X\omega_{B_0})$.
Later on, we will show that the above construction
nevertheless induces a bijection between the $G$-equivalence classes of
$G$-invariant star products and formal series of elements of the space
$H^2_{\mbox{\rm\tiny dR},\Gind}(Q)= Z^2_{\mbox{\rm\tiny dR}}(Q)^G /
\d (\Ginf{T^*Q}^G)$ of second degree cohomology classes of
$G$-invariant de Rham cohomology.
Moreover, this bijection will actually turn out
to be independent of the chosen connection, hence is canonical
(cf.\ \cite[Thm.\ 4.1]{BerBieGut98} for an analogous statement on
general symplectic manifolds).

As an immediate corollary Proposition \ref{GEquivaundBExiProp}
implies:
\begin{COROLLARY}\label{GenQMMExiCor}
Under the assumptions of the proposition the following holds true:
\begin{enumerate}
\item
There is a $G$-equivariant quantum momentum map for a $G$-invariant
star product $\star$ on $(T^*Q,\omega_{B_0})$, if and only if there
is a $G$-equivariant quantum momentum map for the star product
$\starNuB$, where $B=B_0 + B_+$ denotes a formal series of closed
$G$-invariant two-forms on $Q$ as in Proposition
\ref{GEquivaundBExiProp} i.).
\item
If $\PolFun[[\nu]]$ is in addition a $\star$-subalgebra, then every
$G$-equivariant quantum momentum map $J$ for $\star$ satisfies
$\Jbold{\xi} \in \PolFun[[\nu]]$ for all $\xi\in \mathfrak g$.
\end{enumerate}
\end{COROLLARY}
\begin{PROOF}
For the proof of i.) consider a $G$-equivariant quantum momentum
map $J$ for $\star$ and a $G$-equivalence $\mathcal T$ to the star
product $\starNuB$. Then $J^B$ with $J^B(\xi):= \mathcal T
\Jbold{\xi}$ clearly defines a $G$-equivariant quantum momentum map
for $\starNuB$. Vice versa, every $G$-equivariant quantum momentum
map $J^B$ for $\starNuB$ defines a $G$-equivariant quantum momentum
map $J$ for $\star$ via $\Jbold{\xi}:=
\mathcal T^{-1}J^B(\xi)$. For the proof of ii.) apply
Proposition \ref{GEquivaundBExiProp} to show that $\mathcal T$ can
be chosen to preserve $\PolFun[[\nu]]$. Given a $G$-equivariant
quantum momentum map $J$ for $\star$ we then get one for $\starNuB$
by $J^B(\xi):= \mathcal T \Jbold{\xi}$. But from Proposition
\ref{starkBQMMExProp} we have that $J^B$ is of form $J^B(\xi) =
\Pol{\xi_\Qind} +\pi^* \jbold{\xi}\in
\PolFun[[\nu]]$, where $j\in C^1(\mathfrak g,\Cinf{Q})[[\nu]]$
satisfies the conditions in Eq.\ (\ref{starkBQMMExEq}). This
implies in particular that $\Jbold{\xi} = \mathcal
T^{-1}J^B(\xi)\in \PolFun[[\nu]]$. Since any other $G$-equivariant
quantum momentum map $J'$ for $\star$ differs from $J$ by an
element of ${\mathfrak g^*}^G + \nu {\mathfrak g^*_c}^G[[\nu]]$ we
obtain $J'(\xi) \in \PolFun[[\nu]]$ for every $G$-equivariant $J'$
which is a quantum momentum map for $\star$.
\end{PROOF}

In view of the second part of the above corollary it now becomes clear that
one of the assumptions we made for our reduction scheme -- namely that
$\Jbold{\xi} \in \PolFun[[\nu]]$ -- is in fact not an additional
assumption but a consequence of the assumption that
$\PolFun[[\nu]]$ is a $\star$-subalgebra.

For the purposes of the following section, where we will compute
the characteristic class of a reduced star product the results
achieved so far would completely suffice. But with a little more
effort we can give a classification of the $G$-invariant
star products on cotangent bundles up to $G$-equivalence, a result
which is of independent interest. To this end we show in a first step
the following proposition. Its proof is rather technical but yields
the key tools for the main results of the last part of this section.
\begin{PROPOSITION}\label{BBStrichGhomProp}
Let $\starNuB$ resp.\ ${\star'}_0^{B'}$ denote a $G$-invariant star
products on $(T^*Q,\omega_{B_0})$ which is obtained from a
$G$-invariant torsion free connection $\nabla$ resp.\ $\nabla'$ and
a $G$-invariant formal series of closed two-forms $B$ resp.\ $B'$
on $Q$ starting with $B_0$. Then $\starNuB$ and ${\star'}_0^{B'}$
are $G$-equivalent, if and only if $\frac{1}{\nu}B_+ =\frac{1}{\nu}
(B-B_0)$ and $\frac{1}{\nu}B'_+=
\frac{1}{\nu}(B'-B_0)$ are $G$-cohomologous, i.e.\ if and only if
there is a $G$-invariant formal series of one-forms $\beta$ on $Q$
such that $\frac{1}{\nu} B_+ = \frac{1}{\nu}B'_+ + \d \beta$.
\end{PROPOSITION}
\begin{PROOF}
Let us first assume that $B_+$ and $B'_+$ are $G$-cohomologous. Then
we want to prove that $\starNuB$ and ${\star'}_0^{B'}$ are $G$-equivalent.
To this end we need the following result about standard ordered star
products obtained from torsion free connections.
\begin{SUBLEMMA}
Let $\nabla$ and $\nabla'$ denote two torsion free connections on
$Q$ and $\varrho_0$ resp.\ $\varrho'_0$ the corresponding standard
ordered representation of the star product $\starNu$ resp.\
$\star'_0 $ on $(T^*Q,\omega_0)$. Then there is a uniquely
determined formal series $\mathcal S$ of differential operators on
$\Cinf{T^*Q}$ such that $\mathcal S F \in \PolFun[[\nu]]$ for all
$F\in\PolFun[[\nu]]$ and
\begin{equation}\label{repStrichrepEq}
\varrho'_0 (f) = \repNu{\mathcal S f}\quad\textrm{for all }f\in
\Cinf{T^*Q}[[\nu]].
\end{equation}
This implies that $\mathcal S$ is an equivalence transformation from
$(\Cinf{T^*Q}[[\nu]],\star'_0)$ to $(\Cinf{T^*Q}[[\nu]],\starNu)$.
Moreover, $\mathcal S$ satisfies $\mathcal S \pi^*\chi
= \pi^*\chi$ for all $\chi \in \Cinf{Q}[[\nu]]$. Consequently, one has
\begin{equation}\label{FopSConjEq}
\frac{1}{\nu}\mathcal S \ad_{\star'_0} (\pi^*\chi) \mathcal S^{-1}
= \frac{1}{\nu}\ad_{\starNu} (\pi^*\chi).
\end{equation}
\end{SUBLEMMA}
\begin{INNERPROOF}
First recall that the operators of symmetric covariant derivation
$\sKov$ and $\sKov '$ satisfy $\sKov ' = \sKov - \d x^i \vee \d x^j
\vee i_s(S(\partial_{x^i},\partial_{x^j}))$, where the symmetric
tensor field $S$ is given by $\nabla'_{X} Y = \nabla_X Y + S(X,Y)$.
Now it is easy to see that $\varrho'_0 (F)$ lies in the image of
$\varrho_0$ for all $F\in \PolFun[[\nu]]$. By injectivity of the restriction
of $\varrho_0$ to $\PolFun[[\nu]]$ one can
define a map $\mathcal S :\PolFun[[\nu]]\to \PolFun[[\nu]]$ by
$\mathcal S F := {\varrho_0}^{-1} (\varrho'_0(F))$.
Using the explicit form of the standard ordered representations it is
immediate to check that this map is given by a formal series of
differential operators on $\PolFun$ and that this series starts
with $\id$. By definition, $\mathcal S$
satisfies Eq.\ (\ref{repStrichrepEq}) on polynomial functions in
the momenta and $\mathcal S(F \star'_0 F')
= \mathcal S F \starNu \mathcal S F'$ for all $F,F'\in
\PolFun[[\nu]]$. This implies that $\mathcal S$, which
clearly extends uniquely to a mapping on $\Cinf{T^*Q}[[\nu]]$,
is an equivalence from $\star'_0$ to $\starNu$ and satisfies
Eq.\ (\ref{repStrichrepEq}). Uniqueness of $\mathcal S$
is again a direct consequence of $\varrho_0$ being injective when
restricted to $\PolFun[[\nu]]$. For the proof of the further
properties of $\mathcal S$ observe first that $\varrho_0$ satisfies
$\repNu{\pi^*\chi F}\psi = \chi \repNu{F}\psi$ for all
$\chi,\psi\in \Cinf{Q}[[\nu]]$ and that an analogous relation holds
for $\varrho'_0$. Using the definition of $\mathcal S$ this yields
that $\mathcal S (\pi^*\chi F) = \pi^*\chi \mathcal S F$ for all $F
\in \PolFun[[\nu]]$, hence $\mathcal S$ commutes with
left-multiplications by formal functions pulled-back from $Q$. In
particular, we obtain by setting $F=1$ that $\mathcal
S \pi^*\chi = \pi^*\chi$. From this relation Eq.\ (\ref{FopSConjEq}) is
immediate, since $\mathcal S$ is an equivalence from $\star'_0$ to
$\starNu$.
\end{INNERPROOF}

Using $G$-equivariance of $\varrho'_0$ and $\varrho_0$
the preceding sublemma entails that $\varrho'_0(f)=
\repNu{\mathcal S f} = \repNu{\Phi_g^* \mathcal S
\Phi_{g^{-1}}^* f} $ for all $f \in \Cinf{T^*Q}[[\nu]]$.
Clearly, $\Phi_g^* \mathcal S \Phi_{g^{-1}}^*$ preserves
$\PolFun[[\nu]]$. But since $\mathcal S$ is the uniquely determined
map which satisfies Eq.\ (\ref{repStrichrepEq}) and preserves
$\PolFun[[\nu]]$ we conclude that $\mathcal S =
\Phi_g^* \mathcal S \Phi_{g^{-1}}^*$, i.e.\ $\mathcal S$ is a
$G$-equivalence from $\star'_0$ to $\starNu$. We claim that if
$\frac{1}{\nu}B_+ = \frac{1}{\nu} B'_+ + \d \beta$ with $\beta\in
\Ginf{T^*Q}^G[[\nu]]$ one can use $\mathcal S$ to build a
$G$-equivalence from ${\star'}_0^{B'}$ to $\starNuB$. To this end
consider ${A'}^j = A_0^j + {A'}_+^j$ with local potentials $A_0^j$
of $B_0$ and ${A'}_+^j$ of $B'_+$ over some $O_j$ which is assumed to be
an element of a good open cover of $Q$. Then $A^j=A_0^j + {A'}_+^j + \nu \beta$
is a local potential of $B$. Composition of $\mathcal S$
with the local isomorphisms $({\mathcal A'}_0^j)^{-1}$ and
$\mathcal A_0^j$ defined in Eq.\ (\ref{AkappaDefEq}) gives rise to
a local isomorphism $\mathcal T_j:= \mathcal A_0^j
\mathcal S ({\mathcal A'}_0^j)^{-1}$ from
$(\Cinf{T^*O_j}[[\nu]],{\star'}_0^{B'})$ to $(\Cinf{T^*O_j}[[\nu]],
\starNuB)$.
We now have to show that these local isomorphisms actually glue together
to a globally defined $G$-equivalence. To verify this
recall from \cite[Thm.\ 3.4]{BorNeuPflWal03} that $({\mathcal
A'}_0^j)^{-1} {\mathcal A'}_0^i f = \Tilde f(1)$ for every $f\in
\Cinf{T^*(O_j\cap O_i)}[[\nu]]$, where $\Tilde f$ is given by the
unique solution of the differential equation $\frac{\d}{\d t}
\Tilde f (t) = \frac{1}{\nu} \ad_{\star'_0} (\pi^*a_{ji})
\Tilde f(t)$ with $\Tilde f(0)=f$. Here, $a_{ji}\in \Cinf{O_j
\cap O_i}[[\nu]]$ satisfies $\d a_{ji} = {A'}^j - {A'}^i$ over $O_j
\cap O_i$. Since $\mathcal S$ is an equivalence from
$\star'_0$ to $\starNu$ and $\mathcal S \pi^*\chi = \pi^*\chi$, this
implies that $\mathcal S ({\mathcal A'}_0^j)^{-1} {\mathcal A'}_0^i
\mathcal S^{-1} f = \Tilde{\Tilde{f}}(1)$, where
$\Tilde{\Tilde{f}}$ solves $\frac{\d}{\d t}
\Tilde{\Tilde{f}} (t) = \frac{1}{\nu} \ad_{\starNu} (\pi^*a_{ji})
\Tilde{\Tilde{f}}(t)$ with $\Tilde{\Tilde{f}}(0)=f$. But from the
choice of the local potenials $A^j$ and ${A'}^j$ we obtain  $\d
a_{ji} = {A}^j - {A}^i$ over $O_j \cap O_i$. Therefore, $\mathcal
S ({\mathcal A'}_0^j)^{-1} {\mathcal A'}_0^i \mathcal S^{-1}$ is
equal to $(\mathcal A_0^j)^{-1} \mathcal A_0^i$, i.e.\ over
$\Cinf{T^*(O_j\cap O_i)}[[\nu]]$ we have
$\mathcal T_j {\mathcal T_i}^{-1} = \id$.
But this entails that we can define a global isomorphism
$\mathcal T$ from $(\Cinf{T^*Q}[[\nu]],{\star'}_0^{B'})$ to
$(\Cinf{T^*Q}[[\nu]],\starNuB)$ by $(\mathcal T f)|_{T^*O_j}:=
\mathcal T_j f|_{T^*O_j}$. From the fact that ${A'}^j$ and $A^j$
coincide at zeroth order in $\nu$ it is obvious that this
isomorphism is in fact an equivalence transformation. It remains to
show that $\mathcal T$ is $G$-invariant. Let us fix $g\in G$. Using
that $\mathcal S$ is $G$-invariant we obtain $\Phi_g^*\mathcal
T_j\Phi_{g^{-1}}^*= \Phi_g^* \mathcal A_0^j\Phi_{g^{-1}}^*
\mathcal S \Phi_g^*({\mathcal A'}_0^j)^{-1}\Phi_{g^{-1}}^*$ over
$\phi_{g^{-1}}(O_j)$. Now consider an index $i$ such that
$\phi_{g^{-1}}(O_j)\cap O_i\neq \emptyset$. Then we claim that
$\Phi_g^*\mathcal T_j\Phi_{g^{-1}}^*{\mathcal T_i}^{-1}=\id$ over
$\phi_{g^{-1}}(O_j)\cap O_i$. To prove this we may assume without loss of
generality that $\phi_{g^{-1}}(O_j)\cap O_i$ is contractible
since in case it were not contractible we could cover it by open
contractible subsets and use the following argument for each
element of this covering. By $G$-invariance of $B'$ we
have $\d (\phi_g^* {A'}^j-{A'}^i)=0$. Hence there  exist
formal functions $b_{ji}$ over $\phi_{g^{-1}}(O_j)\cap O_i$
such that $\phi_g^* {A'}^j-{A'}^i = \d b_{ji}$. Then
$\Phi_g^*{{\mathcal A'}_0^j}^{-1} \Phi_{g^{-1}}^*{{\mathcal A'}_0^i}$
turns out to be the local automorphism of $\star'_0$ generated by
$\frac{1}{\nu}\ad_{\star'_0}(\pi^* b_{ji})$. Like in the argument which showed
that $\mathcal T$ is well-defined one concludes that  $\mathcal S
\Phi_g^*({\mathcal A'}_0^j)^{-1} \Phi_{g^{-1}}^*{{\mathcal A'}_0^i}
\mathcal S^{-1}$ is the local automorphism of $\starNu$ generated
by $\frac{1}{\nu}\ad_{\starNu}(\pi^* b_{ji})$. On the other hand
$\Phi_g^*({\mathcal A}_0^j)^{-1} \Phi_{g^{-1}}^*
{\mathcal A}_0^i$ coincides with this local automorphism, since by
$G$-invariance of $\beta$ and the choice of the local potentials
$A^j$ and ${A'}^j$ the equation $\phi_g^*A^j - A^i =\phi_g^*{A'}^j
- {A'}^i = \d b_{ji}$ is valid. But this implies that
$\Phi_g^*\mathcal T_j \Phi_{g^{-1}}^* {\mathcal T_i}^{-1} = \id$
over $\phi_{g^{-1}}(O_j)\cap O_i$. Hence $\mathcal T$ is
$G$-invariant. So we have shown that $\starNuB$ and
${\star'}_0^{B'}$ are $G$-equivalent, if $\frac{1}{\nu}B_+$ and
$\frac{1}{\nu}B'_+$ are $G$-cohomologous. For the proof of the
converse statement assume that $\starNuB$ and ${\star'}_0^{B'}$ are
$G$-equivalent and that $l\geq 1$ is the smallest index such that
$B_l$ and $B'_l$ are not $G$-cohomologous. As usual we have set hereby $B
= B_0 + \sum_{r=1}^\infty\nu^r B_r$ and $B' = B_0 +
\sum_{r=1}^\infty\nu^r B'_r$. Now consider $B'':= B_0 +
\sum_{r=0}^{l-1} \nu^r B_r + \sum_{r=l}^\infty \nu^r B'_r$. Then
$B''_+$ is $G$-cohomologous to $B'_+$. Consequently,
we know from above that ${\star'}_0^{B'}$ is $G$-equivalent to
$\star_0^{B''}$. But this implies that $\starNuB$ and $\star_0^{B''}$ are
also $G$-equivalent. It is now immediate to deduce from Lemma
\ref{BDifferLem} that the describing bidifferential operators
$C_r^{\starNuB}$ and $C_r^{\star_0^{B''}}$ coincide for $r
=0,\ldots, l$ and that
\[
C_{l+1}^{\starNuB}(f,f') - C_{l+1}^{\starNuB}(f',f) -
C_{l+1}^{\star_0^{B''}}(f,f') + C_{l+1}^{\star_0^{B''}}(f',f) =
(\pi^*B'_l - \pi^*B_l)(X_f^{B_0},X_{f'}^{B_0}).
\]
But then the $G$-equivalence of $\starNuB$ and $\star_0^{B''}$
implies that $\pi^*B'_l - \pi^*B_l$ is $G$-exact (cf.\ \cite[Thm.\
2.1]{BerBieGut98}), i.e.\ there is a $G$-invariant one-form $\Xi$
on $T^*Q$ such that $\pi^*B'_l - \pi^*B_l= \d \Xi$. Thus $B'_l
= \d i^*\Xi + B_l$, where $i^*\Xi$ is $G$-invariant.
Hence $B_l$ and $B'_l$ are $G$-cohomologous, which is a
contradiction, proving the other direction of the statement of the
proposition.
\end{PROOF}

Using the construction in the proof of Proposition
\ref{GEquivaundBExiProp} and the preceding proposition we can state
the following classification result:
\begin{THEOREM}\label{GenGInvClassThm}
\begin{enumerate}
\item
To every star product $\star$ on $(T^*Q,\omega_{B_0})$ which is
invariant with respect to the lifted action of a $G$-action on $Q$
one can assign a formal series in the second $G$-invariant de Rham
cohomology of $Q$ by
\begin{equation}\label{starToClassDefEq}
c_\Gind : \star \mapsto c_\Gind (\star):=\frac{1}{\nu}[B_+]_\Gind
\in H^2_{\mbox{\rm\tiny dR},\Gind}(Q)[[\nu]],
\end{equation}
where $\frac{1}{\nu}B_+$ denotes a formal series of
$G$-invariant closed two-forms on $Q$ as in Proposition
\ref{GEquivaundBExiProp}.
\item
The map $c_\Gind$ in Eq.\ (\ref{starToClassDefEq}) is
independent of the chosen $G$-invariant torsion free connection on
$Q$ which was used to define $\frac{1}{\nu}B_+$. Moreover,
$c_\Gind$ induces by
$[\star]_\Gind \mapsto c_\Gind (\star)$ a bijection
between the set of $G$-equivalence classes of $G$-invariant star
products as in i.) and
$H^2_{\mbox{\rm\tiny dR},\Gind}(Q)[[\nu]]$.
\end{enumerate}
\end{THEOREM}
\begin{PROOF}
For the proof of i.) we just have to verify that $c_\Gind$ is
well-defined. To this end consider two $G$-equivalence
transformations $\mathcal T,\mathcal T'$ and two formal series
$B,B'$ of closed $G$-invariant two-forms on $Q$ starting with $B_0$
such that $\mathcal T\star = \starNuB$ and $\mathcal T'
\star=\star_0^{B'}$ as in Proposition \ref{GEquivaundBExiProp}.
Then $\starNuB$ and $\star_0^{B'}$ are $G$-equivalent, clearly.
Proposition \ref{BBStrichGhomProp} implies that $\frac{1}{\nu}B_+$
and $\frac{1}{\nu}B'_+$ are $G$-cohomologous, therefore
$c_\Gind$ is well-defined. For the proof of ii.) consider
$G$-equivalences $\mathcal T,\mathcal T'$ and two formal series
$B,B'$ of closed $G$-invariant two-forms on $Q$ starting with $B_0$
such that $\mathcal T\star = \starNuB$ and $\mathcal T'
\star={\star'}_0^{B'}$, where $\starNuB$ and ${\star'}_0^{B'}$ are
obtained from different connections $\nabla$ and $\nabla'$. Then
$\starNuB$ and ${\star'}_0^{B'}$ are $G$-equivalent.
Proposition \ref{BBStrichGhomProp} implies again that $\frac{1}{\nu}B_+$
and $\frac{1}{\nu}B'_+$ are $G$-cohomologous which shows that $c_\Gind$
is independent of the connection used to construct
$\frac{1}{\nu}B_+$. Furthermore, observe that for $G$-equivalent
star products $\star$ and $\star'$ one has $c_\Gind (\star) = c_\Gind
(\star')$, since there exists a $G$-equivalence from $\star$ to
$\starNuB$, hence we obtain a $G$-equivalence from $\star'$ to $\starNuB$.
This implies that $c_\Gind$ induces a mapping from the set of
$G$-equivalence classes of $G$-invariant star products as in i.)
to $H^2_{\mbox{\rm\tiny dR},\Gind}(Q)[[\nu]]$ as given above.
To prove surjectivity of this
map consider $\starNuB$, where $B = B_0 +
\nu \beta$ and $\beta\in Z^2_\dR(Q)^G[[\nu]]$ is an arbitrary
formal series of closed $G$-invariant two-forms on $Q$. By
definition of $c_\Gind$ and choosing  $\id$ as $G$-equivalence
we obtain  $c_\Gind (\starNuB) = [\beta]_\Gind$. Since $\beta$ is
arbitrary this proves surjectivity. To prove injectivity let
$\star,\star'$ be star products with $c_\Gind(\star) =
c_\Gind(\star')$. By Proposition
\ref{BBStrichGhomProp} the corresponding star products $\starNuB$
and $\star_0^{B'}$ are $G$-equivalent which implies that $\star$ and
$\star'$ are $G$-equivalent.
\end{PROOF}
\section{The Characteristic Class of the Reduced Star Products
$\starred{J,\mu}$}
\label{CharClassSec}
In this section we want to compute the characteristic class of the
reduced star products $\starred{J,\mu}$ in order to clarify how
the equivalence classes of these products depend on the initially chosen
parameters of the reduction scheme. To this end
we proceed in two steps. Under the general assumption of a proper
and free $G$-action
on $Q$ we first consider a star product
$\starNuB$ on $(T^*Q,\omega_{B_0})$ constructed from a
$G$-invariant torsion free connection $\nabla$ on $Q$
and a formal series $B$ of $G$-invariant closed two-forms on $Q$.
Recall that due to the properness of the $G$-action such a $G$-invariant
connection exists on $Q$ by Palais' Theorem and that the resulting $\starNuB$
is $G$-invariant by the results of the preceding section.
Additionally, we assume that $j \in C^1(\mathfrak
g,\Cinf{Q})[[\nu]]$ satisfies the conditions given in Eq.\
(\ref{starkBQMMExEq}) so that we can use $J^B$ with $J^B(\xi)=
\Pol{\xi_\Qind} + \pi^*\jbold{\xi}$ as $G$-equivariant quantum
momentum map. Then it is possible to compute the characteristic
class of the reduced star product ${\starNuB}^{J^B,\mu}$ explicitly
in terms of $B, j, \mu$ and the connection on $p:Q\to
\cc{Q}$. In a second step we use the relation between a
star product $\star$ that satisfies all necessary assumptions for
the applicability of our reduction procedure and a $G$-equivalent
star product $\starNuB$ to relate the characteristic class of
${\star}^{J,\mu}$ to the one of ${\starNuB}^{J^B,\mu}$.

Let us now provide a few results needed in the sequel for the computation
of characteristic class of star products.
For more details we refer the reader to \cite{GutRaw99,Neu02} which treats
the case of arbitrary symplectic manifolds and to \cite{BorNeuPflWal03}, where
the special case of cotangent bundles is considered.

Recall that the characteristic class of a
star product $\star$ on $(M,\omega)$ is an element of
$\frac{[\omega]}{\nu} + H^2_{\mbox{\tiny dR}}(M)[[\nu]]$.
For its computation one first has to find
local derivations of $(\Cinf{O_j}[[\nu]],\star)$,
so-called local $\nu$-Euler derivations, where $\{O_j\}_{j\in I}$ is a good
open cover of $M$. These local $\nu$-Euler derivations are of form
\begin{equation}\label{genEulDerEq}
E_j  = \nu \partial_\nu + \Lie_{\xi_j} + \sum_{r=1}^\infty \nu^r
D_{j,r},
\end{equation}
where $\xi_j$ is a local conformally symplectic vector field (i.e.\
$\Lie_{\xi_j}\omega|_{O_j}=\omega|_{O_j}$), and the $D_{j,r}$ are
locally defined differential operators on $\Cinf{O_j}$. With the
help of these the characteristic class  can be
determined except for the part of order zero in the formal
parameter. For the computation of that term one additionally needs
an explicit expression for the anti-symmetric part
$C_2^-(f,f')=\frac{1}{2} \left(C_2(f,f')- C_2(f',f)\right)$ of the
bidifferential operator  describing the considered star product
in the second order of the formal parameter.
More explicitly, to determine the characteristic class from the $E_j$
one considers $E_i - E_j$ over $O_i\cap O_j$, which
is a quasi-inner derivation, i.e.\ there are local formal functions
$d_{ij}\in \Cinf{O_i\cap O_j}[[\nu]]$ such that $E_i -
E_j=\frac{1}{\nu}\ad_\star(d_{ij})$. Now, whenever $O_i \cap O_j
\cap O_k\neq\emptyset$ the sums $d_{ijk} = d_{jk} - d_{ik}
+ d_{ij}$ lie in $\mathbb C[[\nu]]$ and define a $2$-cocycle whose
\v{C}ech class $[d_{ijk}]\in H^2(M,\mathbb C)[[\nu]]$ is
independent of the choices made. The corresponding class
$d(\star) \in H^2_{\mbox{\rm\tiny dR}}(M)[[\nu]]$ is called
Deligne's derivation-related class of $\star$. In addition, let
${C_2^-}^\sharp$ denote the image of $C_2^-$ under the projection
onto the second component of the decomposition $H^2_{\mbox{\rm\tiny
Chev, nc}}(\Cinf{M},\Cinf{M})= \mathbb C \oplus H^2_{\mbox{\rm\tiny
dR}}(M)$ which describes the second cohomology group of the null-on-constants
Chevalley cohomology of $(\Cinf{M},\{\,\,,\,\,\})$,
taken with respect to the adjoint representation.
Together, $d(\star)$ and ${C_2^-}^\sharp$ define
the characteristic class $c(\star)$ of $\star$ by
\begin{equation}\label{CharKlassAllgDefEq}
c(\star)_0 = -2 {C_2^-}^\sharp\quad\textrm{ and }\quad \partial_\nu
c(\star) = \frac{1}{\nu^2} d(\star).
\end{equation}
The so-defined element of $\frac{[\omega]}{\nu} +
H^2_{\mbox{\rm\tiny dR}}(M)[[\nu]]$ classifies the equivalence
classes of star products on a symplectic manifold $(M,\omega)$ in a
functorial way (cf.\ \cite[Thm.\ 6.4]{GutRaw99}).

In the following $\{O_i\}_{i\in I}$ denotes a $G$-invariant good
open cover of $Q$ which projects via $p$ to a good open cover
$\{\cc{O_i}\}_{i\in I}=\{p(O_i)\}_{i\in I}$ of $\cc{Q}$. Such a
cover exists due to the fact that the action on $Q$ is proper. Our
first goal is to define local $\nu$-Euler derivations of
${\starNuB}^{J^B,\mu}$ on every $\cc{O_i}$ using certain local
$\nu$-Euler derivations of $\starNuB$ on $O_i$. Unfortunately, an
arbitrary $\nu$-Euler derivation of $\starNuB$ cannot be projected
down to such a derivation of ${\starNuB}^{J^B,\mu}$ in general,
since such derivations usually neither preserve $\PolFun[[\nu]]$
nor are $G$-invariant and even then do not preserve the ideal of
$G$-invariant formal functions generated by the $G$-equivariant
quantum momentum map. Therefore, we have to find appropriately
modified $\nu$-Euler derivations of $\starNuB$, where we let us
lead by intuition rather than by a deductive procedure. Actually,
the relation between $b_0$ and $B_0$ in the lowest order of the
characteristic classes of $\starNuB^{J^B,\mu}$ and $\starNuB$
suggests that a similar relation might also hold in  higher orders.
In analogy to classical phase space reduction we therefore consider
the following formal series of closed two-forms on $Q$:
\begin{equation}\label{ConsiderBGammaEq}
B + \d \Gamma_{\Check j - \mu},
\end{equation}
where $\Gamma$ has been defined by Eq.\ (\ref{GamDefEq}) and has been
extended by $\mathbb C[[\nu]]$-linearity. A straightforward argument
which is completely analogous to the computation in the proof of
Theorem \ref{BNullRedThm} now shows by $G$-equivariance of $j$
that $B+\d \Gamma_{\Check j - \mu}$ is a formal series of
$G$-invariant closed horizontal two-forms on $Q$. Hence
there exists a uniquely defined formal series $b$ of closed two-forms
on $\cc{Q}$ such that
\begin{equation}\label{BGammaDefbEq}
B + \d \Gamma_{\Check j - \mu} = p^*b.
\end{equation}
Now we choose local potentials $a^i$ of $b$ over $\cc{O_i}$, i.e.\
$\d a^i = b|_{\cc{O_i}}$. For $\cc{O_i}\cap\cc{O_j}$ we choose
local formal functions $f_{ij}$ with
$\d f_{ij}|_{\cc{O_i}\cap\cc{O_j}}
= (a^i - a^j)|_{\cc{O_i}\cap\cc{O_j}}$. Furthermore, we consider
the local formal one-forms $A^i$ on $O_i$ defined by $A^i := p^*a^i
-\Gamma_{\Check j - \mu}$. Then $A^i - A^j = \d p^*f_{ij}$
holds true on $O_i\cap O_j$ by construction. Using these $A^i$ we now consider
the local isomorphisms $\mathcal A_0^i:(\Cinf{T^*
O_i}[[\nu]],\starNu) \to (\Cinf{T^* O_i}[[\nu]],\starNuB)$
from Eq.\ (\ref{AkappaDefEq}) and claim the following:
\begin{PROPOSITION}
With notations from above the following holds true:
\begin{enumerate}
\item
The mappings $\mathcal E_i : \Cinf{T^* O_i}[[\nu]] \to \Cinf{T^*
O_i}[[\nu]]$ defined by
\begin{equation}\label{NuEulerobenDefEq}
\mathcal E_i := \mathcal A_0^i \mathcal H (\mathcal A_0^i)^{-1}
\end{equation}
are $G$-invariant local $\nu$-Euler derivations of $\starNuB$ which
preserve $\mathcal P(O_i)[[\nu]]$ and $\mathcal P(O_i)[[\nu]]\cap
I_{\starNuB,\mu}$.
\item
The mappings $\mathsf E_i : \hor(\mathcal P(O_i)^G)[[\nu]] \to
\hor(\mathcal P(O_i)^G)[[\nu]]$ defined by
\begin{equation}
\mathsf E_i F := \hor_{\Check j_0-\mu_0}\left(\frac{\id}{\id -\nu
\triangle_{\mu,\starNuB}} \mathcal E_i F
\right),\quad F \in \hor(\mathcal P(O_i)^G)[[\nu]]
\end{equation}
are local derivations of $(\hor(\mathcal P(O_i)^G)[[\nu]],
{\bullet^B_0}^{J^B,\mu})$. The corresponding mappings
$E_i:=l^{-1}\circ\mathsf E_i \circ l$ are local $\nu$-Euler
derivations of $(\mathcal P(\cc{O_i})[[\nu]],{\starNuB}^{J^B,\mu})$
which uniquely extend to such derivations of
$(\Cinf{T^*\cc{O_i}}[[\nu]], {\starNuB}^{J^B,\mu})$.
\item
On $\Cinf{T^*(\cc{O_i}\cap\cc{O_j})}[[\nu]]$ one has
\begin{equation}\label{NuEulDiffEq}
E_i - E_j=
\frac{1}{\nu}\ad_{{\starNuB}^{J^B,\mu}} ((\nu\partial_\nu - \id)
\cc{\pi}^*f_{ij}).
\end{equation}
This implies that the characteristic class $c({\starNuB}^{J^B,\mu})$ is
given by
\begin{equation}\label{CharKlaOhneNullEq}
c({\starNuB}^{J^B,\mu}) = \frac{1}{\nu}[\cc{\pi}^*b] -
[\cc{\pi}^*b_1] + c({\starNuB}^{J^B,\mu})_0.
\end{equation}
\end{enumerate}
\end{PROPOSITION}
\begin{PROOF}
For the proof of i.) first observe that by definition of $A^i$ one
has $\d A^i = B|_{O_i}$. Therefore, the mappings $\mathcal A_0^i$
are in fact local isomorphisms from $(\Cinf{T^*
O_i}[[\nu]],\starNu)$ to $(\Cinf{T^* O_i}[[\nu]],\starNuB)$. By
their explicit form it is obvious that they preserve $\mathcal P
(O_i)[[\nu]]$. Moreover, they also turn out to be $G$-invariant due
to the $G$-invariance of the $A^i$. Since $\mathcal H =
\Lie_{\xi_0} + \nu\partial_\nu$ is evidently $G$-invariant and
preserves $\mathcal P (O_i)[[\nu]]$ these properties hold for
$\mathcal E_i$, as well. In addition, $\mathcal E_i$ is a local
$\nu$-Euler derivation of $\starNuB$, as $\mathcal H$ is a global
$\nu$-Euler derivation of $\starNu$. It remains to show that
$\mathcal E_i$ preserves $\mathcal P (O_i)[[\nu]]\cap
I_{\starNuB,\mu}$. But this follows from a straightforward proof of
$\mathcal E_i J^B(\xi)|_{T^*O_i} = J^B(\xi)|_{T^*O_i} -
\pair{\mu}{\xi} +
\nu\partial_\nu \pair{\mu}{\xi}$ which uses the explicit shape of $J^B$
and $\mathcal A_0^i$. Using i.) it is rather evident that $\mathsf
E_i$ defines a local derivation of ${\bullet_0^B}^{J^B,\mu}$ on
$\hor(\mathcal P(O_i)^G)[[\nu]]$, and we only have to show that
$E_i$ is of form provided in Eq.\ (\ref{genEulDerEq}). To this end
recall from \cite[Lemma 4.4]{BorNeuPflWal03} that $\mathcal E_i =
\mathcal H +
\Fop{\frac{\id- \exp (-\nu \sKov)}{\nu \sKov}(\nu\partial_\nu -
\id)A^i}$ which directly implies that $E_i=l^{-1} \circ \mathsf E_i
\circ l$ has form $\nu \partial_\nu + \sum_{r=0}^\infty
\nu^r D_{i,r}$, where the $D_{i,r}$ are locally
defined differential operators. Hence, these
mappings uniquely extend to $\Cinf{T^*\cc{O_i}}[[\nu]]$, since
differential operators are completely determined by their values on
polynomial functions in the momenta. We only have to show that
$D_{i,0}=\Lie_{\cc{\xi}_i}$ with a locally defined vector field
$\cc{\xi}_i\in \Ginf{T(T^*\cc{O_i})}$ satisfying $\Lie_{\cc{\xi}_i}
\omega_{b_0}  = \omega_{b_0}$. But this follows from an easy
computation expanding the exponent in the above given expression
for $\mathcal E_i$. For the proof of iii.) this expression
again together with $A^i - A^j = \d p^*f_{ij}$ entails
\begin{eqnarray*}
(\mathsf E_i -\mathsf E_j) F &=& \hor_{\Check j_0 -\mu_0}\left(
\frac{\id}{\id -\nu \triangle_{\mu,\starNuB}}
\frac{1}{\nu}\Fop{(\id - \exp (-\nu \sKov))(\nu\partial_\nu - \id )
p^*f_{ij} } F \right)\\ &=&
\frac{1}{\nu}\ad_{{\bullet_0^B}^{J^B,\mu}}(\pi^* p^*(\nu\partial_\nu
-\id)f_{ij}) F,
\end{eqnarray*}
where the second equality follows from Eqs.\
(\ref{starklinksmultEq}) and (\ref{starkrechtsmultEq}) together
with Lemma \ref{starkBlinkrechtsMultLem} i.). Conjugation of
this result by $l^{-1}$ yields Eq.\ (\ref{NuEulDiffEq}). From the
very definition of Deligne's derivation related class and the
definition of the local functions $f_{ij}$ we thus obtain
$d({\starNuB}^{J^B,\mu})= (\nu \partial_\nu
- \id)[\cc{\pi}^*b]$. By definition of the characteristic class this
implies Eq.\ (\ref{CharKlaOhneNullEq}).
\end{PROOF}

To determine the missing part $c({\starNuB}^{J^B,\mu})_0$ of the
characteristic class of ${\starNuB}^{J^B,\mu}$ one has to compute
the anti-symmetric part of the bidifferential operator describing
${\starNuB}^{J^B,\mu}$ in the second order of the formal parameter.
As this is a rather cumbersome but nevertheless important
computation we only give here the important intermediate steps
and omit details of the proof.
\begin{LEMMA}\label{CharKlaNullLem}
\begin{enumerate}
\item
Writing $f\starNuB f' =  f f' + \nu C_1^{\starNuB}(f,f') +
\nu^2 C_2^{\starNuB}(f,f')+ O(\nu^3)$ for $f,f'\in \Cinf{T^*Q}$ we have:
\begin{eqnarray}\label{CEinsNullBEq}
C_1^{\starNuB}(f,f') &=& \frac{1}{2}\left(
\{f,f'\}_{B_0} - \Delta_0 (f f') + (\Delta_0 f)f' + f(\Delta_0 f')
\right)\\ \label{CZweiNullBEq}
C_2^{\starNuB}(f,f') - C_2^{\starNuB}(f',f) &=& -\frac{1}{2}\left(
\Delta_0 \{f,f'\}_{B_0} - \{\Delta_0 f,f'\}_{B_0} -
\{f,\Delta_0 f'\}_{B_0}
\right)\\
\nonumber
& &- (\pi^*(B_1 - \frac{1}{2}\tr{R}))(X_f^{B_0},X_{f'}^{B_0}),
\end{eqnarray}
where $\tr{R}$ denotes the trace of the curvature tensor of
$\nabla$ and $\Delta_0$ the differential
operator defined in Eq.\ (\ref{DeltaDefEq}).
\item
For $s,t \in \Ginf{\bigvee T \cc{Q}}$ the anti-symmetric part of
the bidifferential operator $C_2^{{\starNuB}^{J^B,\mu}}$ which describes
the star product ${\starNuB}^{J^B,\mu}$ on $\PolFunRed[[\nu]]$ is
given by
\begin{eqnarray}
\lefteqn{C_2^{{\starNuB}^{J^B,\mu}}(\PolRed{s},\PolRed{t}) -
C_2^{{\starNuB}^{J^B,\mu}}(\PolRed{t},\PolRed{s})}\\
&\!\!\!\!=&\!\!\! l^{-1}
\!\left(\!\hor_{\Check j_0 -\mu_0}\! \left(\!
\triangle_{\mu,\starNuB}\!\left\{ \Pol{s^\hor},\Pol{t^\hor}
\right\}_{B_0} \!+C_2^{\starNuB}\!\left(\Pol{s^\hor},\Pol{t^\hor}\right) -
C_2^{\starNuB}\!\left(\Pol{t^\hor},\Pol{s^\hor}\right)
\!\right)\!\right)\nonumber\\
&\!\!\!\!=&\!\!\! - \frac{1}{2}\left(
\delta_0\{\PolRed{s},\PolRed{t}\}_{b_0} -
\{\delta_0\PolRed{s},\PolRed{t}\}_{b_0} -
\{\PolRed{s},\delta_0\PolRed{t}\}_{b_0}
\right)\\
\nonumber
&\!\!\!\!&\!\!\! + \cc{\pi}^*\!\left(- b_1 - \frac{1}{2} r +
\frac{1}{2}
\tau_\lambda\right) (X_{\PolRed{s}}^{b_0}, X_{\PolRed{t}}^{b_0}).
\end{eqnarray}
Hereby, $r$ and $\tau_\lambda$ are the unique closed two-forms on
$\cc{Q}$ determined by $p^*r = - \tr{R} + \d \Gamma_{\Check d}$
and $p^*\tau_\lambda = \tr{\ad(\lambda)}$, where
$d\in C^1(\mathfrak g, \Cinf{Q})$ is defined by $d(\xi):=
\mathsf{div}(\xi_\Qind)$. Moreover, the mapping $\delta_0 :
\PolFunRed\to \PolFunRed$ is given by $\delta_0 (\PolRed{s}):=
l^{-1} (\hor_{\Check j_0 - \mu_0} (\Delta_0 \Pol{s^\hor}))$.
\end{enumerate}
\end{LEMMA}
\begin{PROOF}
For the proof of i.) first recall from \cite[Thm.\ 10]{BorNeuWal98}
that the star product $\starWe$ is of Weyl type. Hence
the describing bidifferential operators at order one
and two of the formal parameter satisfy
$C_1^{\starWe}(f,f')=\frac{1}{2}\{f,f'\}_0$ and
$C_2^{\starWe}(f,f')-C_2^{\starWe}(f',f)=0$. Using these two
relations together with the equalities $(N_{1/2} f)\starNu (N_{1/2}f
') = N_{1/2}(f\starWe f')$ and $\mathcal A^i_0 (f \starNu
f'|_{T^*O_i})= (\mathcal A^i_0 f|_{T^*O_i}) \starNuB (\mathcal
A^i_0 f'|_{T^*O_i})$ a straightforward computation which involves
expansion of $N_{1/2}$, $\mathcal A^i_0$ and the products
$\starWe$, $\starNu$, and $\starNuB$ up to the second order in $\nu$
then shows Eqs.\ (\ref{CEinsNullBEq}) and (\ref{CZweiNullBEq}).
The first equality stated in ii.) is obtained by an immediate computation using
the definition of ${\starNuB}^{J^B,\mu}$. In contrast, the
proof of the second equality turns out to be more involved but
requires nothing more than a consequent application of the definitions.
Last, we provide the argument showing that $r$ and
$\tau_\lambda$ are well-defined.
For $\tau_\lambda$ this is well known by Chern-Weil theory, and
$[\tau_\lambda]$ is a characteristic class of the principal $G$-bundle
$p: Q \to \cc{Q}$.
To prove that $r$ is well-defined observe that
$-i_{\xi_\Qind}\tr{R}=\d d (\xi)$ and $\phi_g^*d (\xi) = d(
\Ad(g^{-1})\xi)$ for all $\xi\in \mathfrak g$ and
repeat the argument used for Theorem
\ref{BNullRedThm} showing that $b_0$ is well-defined.
\end{PROOF}

By definition of the zeroth order of the characteristic
class and the above result one directly obtains
\begin{equation}
c({\starNuB}^{J^B,\mu})_0 = \left[\cc{\pi}^*\left(b_1 +
\frac{r}{2}-
\frac{\tau_\lambda}{2}\right)\right].
\end{equation}
At this point one might expect that the characteristic
class of the reduced star product depends on the chosen
$G$-invariant connection $\nabla$ and
that the geometry of the principal bundle enters the characteristic
class $c({\starNuB}^{J^B,\mu})$ via $[\tau_\lambda]$. As we will
show in the next lemma none of these dependencies occurs.
\begin{LEMMA}\label{CNullBestimmLem}
\begin{enumerate}
\item
Let $\nabla$ and $\nabla'$ be two $G$-invariant torsion free connections
on $Q$. Then the corresponding closed two-forms $r$ and $r'$
constructed as in Lemma \ref{CharKlaNullLem} ii.) are cohomologous.
\item
With notations from Lemma \ref{CharKlaNullLem} ii.) we have
\begin{equation}
r - \tau_\lambda = - \tr{\cc{R}} - \d w.
\end{equation}
Hereby, $\cc{R}$ denotes the curvature of the torsion free
connection $\cc{\nabla}$ on $\cc{Q}$ defined by $T p\,
\nabla_{t^\hor} u^\hor = \cc{\nabla}_t u \circ p$ for $t,u \in
\Ginf{T\cc{Q}}$. Moreover, $w\in \Ginf{T^*\cc{Q}}$ is defined by
$p^*w = \mathrm{H}(W) = W - \sum_{i=1}^{\dim{(G)}}
W({e_i}_\Qind)\Gamma_{e^i}$ with $W\in \Ginf{T^*Q}$ given by $W(X)=
\sum_{l=1}^{\dim{(G)}}\Gamma_{e^l} (\nabla_X{e_l}_\Qind)$.
Consequently, the cohomology classes of $r$ and $\tau_\lambda$
coincide: $[r -
\tau_\lambda]=[0]$.
\end{enumerate}
\end{LEMMA}
\begin{PROOF}
For the proof of i.) write $\nabla'_X Y = \nabla_X Y + S(X,Y)$
with a symmetric $G$-invariant tensor field $S\in
\Ginf{\bigvee^2T^*Q\otimes TQ}$. Obeserve that the $G$-invariant
one-form $\tr{S}$ defined by $\tr{S}(X) := \d
x^i(S(X,\partial_{x^i}))$ satisfies $\tr{R'} = \tr{R} + \d
(\tr{S})$ and $\mathsf{div}'(X)=
\mathsf{div}(X)+\tr{S}(X)$. Defining $\sigma \in C^1(\mathfrak g,
\Cinf{Q})$ by $\sigma(\xi):= \tr{S}(\xi_\Qind)$ it is easy to see that
there exists a unique one-form $s$ on $\cc{Q}$ which satisfies $p^*s =
\tr{S} - \Gamma_{\Check \sigma}$. From the above identities
relating the traces of the curvature tensors and the covariant
divergences it is evident that $r' = r - \d s$. To prove ii.)
first verify that $\cc{\nabla}$ actually defines a torsion
free connection on $\cc{Q}$ and that $w$ is well-defined.
Then the proof consists of an easy computation showing that
\[
(- \tr{R} + \d \Gamma_{\Check d})(t^\hor,u^\hor) -
\tr{\ad(\lambda(t^\hor,u^\hor))} =
p^*(-\tr{\cc{R}}(t,u) - \d w (t,u)),
\]
This computation mainly relies on splitting the definition of $\tr{R}$ into
horizontal and vertical part. But since the
trace of the curvature tensor of a torsion free connection is an
exact two-form (cf.\ \cite[Lemma 16]{BorNeuWal98}) this implies
$[r] = [\tau_\lambda]$.
\end{PROOF}

Collecting our results we have shown:
\begin{THEOREM}\label{CharKlasNuBThM}
The characteristic class of the star product ${\starNuB}^{J^B,\mu}$
on $(T^*\cc{Q},\omega_{b_0})$ is given by
\begin{equation}
c({\starNuB}^{J^B,\mu}) = \frac{1}{\nu}[\cc{\pi}^* b],
\end{equation}
where $J^B(\xi)= \Pol{\xi_\Qind} + \pi^*\jbold{\xi}$ and $b\in
Z^2_\dR(\cc{Q})[[\nu]]$ is determined by $p^*b
= B + \d \Gamma_{\Check j - \mu}$.
\end{THEOREM}

Now we consider an arbitrary star product $\star$ satisfying the
assumptions for our reduction procedure and
present the main result of this section.
\begin{THEOREM}\label{GeneralCharKlassThm}
Let $\star$ be a $G$-invariant star product on
$(T^*Q,\omega_{B_0})$ such that $\PolFun[[\nu]]$ is a
$\star$-subalgebra. Let $J$ denote a $G$-equivariant quantum
momentum map for $\star$ and  $\mu \in {\mathfrak g^*}^G+
\nu {\mathfrak g_c^*}^G[[\nu]]$.
\begin{enumerate}
\item
Assume that $\mathcal T$ is a $G$-equivalence from $\star$ to
$\starNuB$ as in Proposition \ref{GEquivaundBExiProp} ii.), where
$B=B_0 + B_+$ is an appropriate formal series of $G$-invariant
closed two-forms. Then the star product $\starred{J,\mu}$ is
equivalent to ${\starNuB}^{J^B,\mu}$ with $G$-equivariant quantum
momentum map $J^B$ given by $J^B(\xi) := \mathcal T \Jbold{\xi}$.
\item
With notations from i.), the characteristic class of the star
product $\starred{J,\mu}$ is given by
\begin{equation}
c(\starred{J,\mu}) = c({\starNuB}^{J^B,\mu})
=\frac{1}{\nu}[\cc{\pi}^*b],
\end{equation}
where $b$ denotes the uniquely determined formal series of
two-forms on $\cc{Q}$ such that $p^*b = B + \d \Gamma_{\Check j -
\mu}$ and where $j\in C^1(\mathfrak g,\Cinf{Q})[[\nu]]$ is given by
$\jbold{\xi}=i^*\mathcal T \Jbold{\xi}$.
\end{enumerate}
\end{THEOREM}
\begin{PROOF}
i.) is a direct consequence of Proposition \ref{IsoAutoDerRedProp}
i.) and iv.) since a $G$-equivalence $\mathcal T$ exists due to
Proposition \ref{GEquivaundBExiProp} and since $J^B$ is a
$G$-equivariant quantum momentum map for $\starNuB$ by Corollary
\ref{GenQMMExiCor}. By i.) the equality $c(\starred{J,\mu}) =
c({\starNuB}^{J^B,\mu})$ follows immediately. By Theorem
\ref{CharKlasNuBThM} we obtain the claim about the explicit form of
the characteristic class of $\starred{J,\mu}$. We only have to
verify that $j$ is actually given by $\jbold{\xi}= i^*\mathcal T
\Jbold{\xi}$, but this is obvious since $J^B$ satisfies
$J^B(\xi) = \Pol{\xi_\Qind} + \pi^*\jbold{\xi}$ by
Proposition \ref{starkBQMMExProp}.
\end{PROOF}

\begin{COROLLARY}
Let $\gamma$ and $\gamma'$ denote two connection one-forms on
$p:Q\to \cc{Q}$ and let $\starred{J,\mu}$ resp.\
$(\starred{J,\mu})'$ be the corresponding reduced star product on
$(T^*\cc{Q},\omega_{b_0})$ resp.\ $(T^*\cc{Q},\omega_{b'_0})$
obtained by reduction of the star product $\star$ on
$(T^*Q,\omega_{B_0})$. Then the characteristic classes of
$\starred{J,\mu}$ and $(\starred{J,\mu})'$ coincide and there is an
isomorphism from $(\Cinf{T^*\cc{Q}}[[\nu]],\starred{J,\mu})$ to
$(\Cinf{T^*\cc{Q}}[[\nu]],(\starred{J,\mu})')$. Moreover, the
corresponding star products $\star_{\Psi_{\mu_0}}^{J,\mu}$ and
$(\starred{J,\mu})'_{\Psi'_{\mu_0}}$ on
$((T^*Q)_{\mu_0},\omega_{\mu_0})$ (cf.\ Remark \ref{PullBackRem})
are equivalent.
\end{COROLLARY}
\begin{PROOF}
Let $\Gamma'_{\Check j_0 - \mu_0}$ be the one-form defined by
$\pair{\Check j_0 - \mu_0}{\gamma'}$. Then one observes first that the
translation $t_{\Gamma_{\Check j_0 - \mu_0}-\Gamma'_{\Check j_0 -
\mu_0}}$ along the fibre which maps the zero level set
$(T^*Q)^0 =\{\zeta_x \in T^*Q\,|\,
\zeta_x(\xi_\Qind (x)) =0\:\text{ for all $\xi \in\mathfrak g$}\}$ to itself
clearly passes to the quotient defining a diffeomorphism of
$T^*\cc{Q}$. Moreover, it is easy to see that this diffeomorphism
consists of a translation $t_{\beta_0}$ along the fibres on
$T^*\cc{Q}$, where $\beta_0$ is the unique one-form on $\cc{Q}$
such that $p^*\beta_0 = \Gamma_{\Check j_0 -
\mu_0}-\Gamma'_{\Check j_0 - \mu_0}$. By definition of $b_0$ and
$b'_0$, where $b'_0$ is defined as in Theorem \ref{BNullRedThm}
using $\Gamma'_{\Check j_0 - \mu_0}$ instead of $\Gamma_{\Check j_0
- \mu_0}$, we have $b'_0 = b_0 - \d \beta_0$. Hence, this
diffeomorphism is in fact a symplectomorphism from
$(T^*\cc{Q},\omega_{b'_0})$ to $(T^*\cc{Q},\omega_{b_0})$. Let
$\beta$ be the unique formal series of one-forms on $\cc{Q}$
such that $p^*\beta = \Gamma_{\Check j -
\mu}-\Gamma'_{\Check j - \mu}$. Then we analogously find that $b' = b -
\d \beta$, therefore $c((\starred{J,\mu})') = \frac{1}{\nu}
[\cc{\pi}^* b'] = \frac{1}{\nu} [\cc{\pi}^* b] =
c(\starred{J,\mu})$. Now we consider the star product
$(\starred{J,\mu})'':= t_{\beta_0}^* \starred{J,\mu}$
on $(T^*\cc{Q},\omega_{b'_0})$, which has characteristic class
$c((\starred{J,\mu})'') = t_{\beta_0}^* c (\starred{J,\mu}) =
\frac{1}{\nu} [t_{\beta_0}^* \cc{\pi}^* b] = \frac{1}{\nu}
[\cc{\pi}^* b] = c((\starred{J,\mu})')$. Therefore it is equivalent
to $(\starred{J,\mu})'$. But then the composition of an equivalence
transformation $T$ from $(\starred{J,\mu})''$ to
$(\starred{J,\mu})'$ with $t_{\beta_0}^*$ yields an isomorphism
from $(\Cinf{T^*\cc{Q}}[[\nu]], \starred{J,\mu})$ to
$(\Cinf{T^*\cc{Q}}[[\nu]], (\starred{J,\mu})')$. One finally concludes
that the corresponding star products on $((T^*Q)_{\mu_0},\omega_{\mu_0})$
are equivalent, since the characteristic class is natural with respect
to diffeomorphisms (cf.\ \cite[Thm.\ 6.4]{GutRaw99}) and
since $t_{\beta_0} = \Psi_{\mu_0}\circ (\Psi'_{\mu_0})^{-1}$
by the above result, where
$\Psi'_{\mu_0}$ is defined as in Theorem \ref{BNullRedThm} using
$\gamma'$ instead of $\gamma$.
\end{PROOF}

Similarly, we also find the dependence of the characteristic class
on different choices of the $G$-equivariant quantum momentum map
and different choices of the momentum value:
\begin{COROLLARY}
\begin{enumerate}
\item
Let $\starred{J,\mu}$ and $\starred{J',\mu}$ denote the star
products on $(T^*\cc{Q}, \omega_{b_0})$ and
$(T^*\cc{Q},\omega_{b'_0})$ obtained from two possibly different
$G$-equivariant quantum momentum maps $J$ and $J'$ satisfying $J -
J' = \Tilde{\mu}
\in {\mathfrak g^*}^G + \nu {\mathfrak g^*_c}^G [[\nu]]$. Then the
characteristic classes fulfill
\begin{equation}
c (\starred{J,\mu}) - c(\starred{J',\mu}) =
\frac{1}{\nu}[\cc{\pi}^* \Tilde{\mu}_\lambda],
\end{equation}
where $\Tilde{\mu}_\lambda \in Z^2_\dR(\cc{Q})[[\nu]]$ is
determined by $p^*\Tilde{\mu}_\lambda =
\pair{\Tilde{\mu}}{\lambda}$.
\item
Let $\starred{J,\mu}$ and $\starred{J,\mu'}$ denote the star
products on $(T^*\cc{Q}, \omega_{b_0})$ and
$(T^*\cc{Q},\omega_{b'_0})$ obtained from two possibly
different momentum values $\mu$ and $\mu'$
and let $\Tilde{\mu} = \mu - \mu' \in {\mathfrak g^*}^G + \nu {\mathfrak
g^*_c}^G[[\nu]]$. Then the characteristic classes satisfy
\begin{equation}
c (\starred{J,\mu}) - c(\starred{J,\mu'}) =
- \frac{1}{\nu}[\cc{\pi}^* \Tilde{\mu}_\lambda],
\end{equation}
where $\Tilde{\mu}_\lambda \in Z^2_\dR(\cc{Q})[[\nu]]$ is
determined by $p^*\Tilde{\mu}_\lambda =
\pair{\Tilde{\mu}}{\lambda}$.
\end{enumerate}
\end{COROLLARY}
\begin{PROOF}
Both claims are direct consequences of Theorem
\ref{GeneralCharKlassThm}. For the proof of i.) one just has to
observe that for a $G$-equivalence $\mathcal T$ from $\star$ to
some $\starNuB$ we have $\jbold{\xi} - j'(\xi)= i^*(\mathcal T
\Jbold{\xi} - \mathcal TJ'(\xi)) =\pair{\Tilde{\mu}}{\xi}$,
whereas the second claim is obvious.
\end{PROOF}

Finally, we are able to recover a relation between the
characteristic class of the original star product $\star$ on
$(T^*Q,\omega_{B_0})$ and the characteristic class of the reduced
star product $\star_{\Psi_{\mu_0}}^{J,\mu}$ (cf.\ Remark
\ref{PullBackRem}) on $((T^*Q)_{\mu_0},\omega_{\mu_0})$ which
already has been observed by M. Bordemann \cite{Bor03} (cf.\ also
\cite{Bor04}) for arbitrary symplectic manifolds.
\begin{COROLLARY}
Denote by $i_{\mu_0}$ the inclusion of $\Check J_0^{-1}(\mu_0)$
into $T^*Q$ and by $\pi_{\mu_0}$ the projection from $\Check
J_0^{-1}(\mu_0)$ to $(T^*Q)_{\mu_0}$. Then the characteristic classes of
$\star$ and $\star_{\Psi_{\mu_0}}^{J,\mu}$ satisfy
\begin{equation}
i_{\mu_0}^*c(\star) = \pi_{\mu_0}^* c
(\star_{\Psi_{\mu_0}}^{J,\mu}).
\end{equation}
\end{COROLLARY}
\begin{PROOF}
By naturality of characteristic classes with respect to
diffeomorphisms we obtain $c (\star_{\Psi_{\mu_0}}^{J,\mu}) =
\Psi_{\mu_0}^* c (\starred{J,\mu})$. The commutative diagram
(\ref{ComClassRedDiag}) and Theorem \ref{GeneralCharKlassThm} then
entail
\[
\pi_{\mu_0}^* c (\star_{\Psi_{\mu_0}}^{J,\mu}) =
\frac{1}{\nu}[\pi_{\mu_0}^*\Psi_{\mu_0}^*\cc{\pi}^*b] =
\frac{1}{\nu}[i_{\mu_0}^*t_{\Gamma_{\Check j_0 - \mu_0}}^*
\pi^*(B +\d \Gamma_{\Check j - \mu})] =
i_{\mu_0}^*c(\star),
\]
where the last equality follows from $\pi \circ t_{\Gamma_{\Check
j_0 - \mu_0}}=\pi$ and the fact that $c(\star) =
\frac{1}{\nu}[\pi^* B]$.
\end{PROOF}

To conclude this section we discuss some special cases of
Theorem \ref{GeneralCharKlassThm} which clarify the result and allow
for some interesting observations.

First consider the star product $\starNu$ and the $G$-equivariant
quantum momentum map $\JN{\xi} =
\Pol{\xi_\Qind}$. Then $c({\starNu}^{J_0,\mu}) = - \frac{1}{\nu}
[\cc{\pi}^*\mu_\lambda]$, where $\mu_\lambda$ is the unique formal
series of closed two-forms on $\cc{Q}$ defined by $p^*\mu_\lambda =
\pair{\mu}{\lambda}$. By Chern-Weil theory one knows that
$[\mu_\lambda]$ is independent of the chosen connection and defines
a formal series with values in the first characteristic classes of
the principal $G$-bundle which only depends on $\mu$. Now choose
$\mu_0$ in order to fix the symplectic form on $T^*\cc{Q}$ to be
$\cc{\omega}_0 - \cc{\pi}^* {\mu_0}_\lambda$. In general, one then
cannot obtain star products on $(T^*\cc{Q},\cc{\omega}_0 -
\cc{\pi}^* {\mu_0}_\lambda)$ in every possible characteristic class
by reduction of the star product $\starNu$ on the original phase space.
The reason for this is that not every de Rham cohomology class in
$H^2_\dR(\cc{Q})$ is equal to a characteristic class, which is a
direct consequence of the fact that the Chern-Weil
homomorphism is not surjective, in general. Moreover, by our result it is
clear that in general `quantization does not commute with
reduction'. Consider for example $\mu_0=0$. Then the star product
$\starNu^{J_0,\mu_+}$ is equivalent to an intrinsically defined
star product $\cc{\star}_0$ (starting
from a torsion free connection $\cc{\nabla}$ on $\cc{Q}$),
if and only if $[{\mu_+}_\lambda]=
[0]$. Furthermore, we will see in the following section that there
are even more conditions which have to be fulfilled in order to
achieve that $\starNu^{J_0,\mu_+}$  equals $\cc{\star}_0$.

Starting with the star products $\starNuB$ on $(T^*Q,\omega_{B_0})$
one can actually get representatives for every characteristic class
of star products on $(T^*\cc{Q},\omega_{b_0})$ by varying $B_+$ and
the corresponding mappings $j_+$. This follows from the
observation that $[b_+]$, where $b=b_0 + b_+$ with
fixed $b_0$, corresponds to the formal series
of $G$-equivariant cohomology classes defined by the pair $(B_+,
\Check j_+ - \mu_+)$ and the well-known fact that the
$G$-equivariant cohomology of $Q$ (which is by definition the
cohomology of the complex of basic differential forms on $Q$) is
isomorphic to the de Rham cohomology of the quotient $\cc{Q}=Q/G$
(cf.\ \cite{GuiSte99}).
\section{Applications and Examples}
\label{AppExaSec}
\subsection{Reduction of $\stark$ and $\starkB$}
\label{ExaRedSubSec}
In this section we apply the reduction scheme developped
in Section \ref{RedSec} to several concrete examples of star
products on $(T^*Q,\omega_{B_0})$. As we want to identify the resulting reduced
products on $(T^*\cc{Q}, \omega_{b_0})$ with naturally defined star
products we will often make use of the fact that certain star
products are determined by their representations.
Hence we will construct representations of the reduced star
products, which is also of independant interest because of
some strong relations to the results of \cite{Emm93a,Emm93b,Got86}, where
the quantization is formulated in terms of representations alone
neglecting the algebra of observables.

As a first step towards the concrete computation of reduced star products
derived from $\stark$ and $\starkB$, we are going to establish a
relation between the reduced star products obtained from
different momentum values $\mu \in {\mathfrak g^*}^G +
\nu{\mathfrak g^*_c}^G[[\nu]]$. This will allow us to restrict
our further considerations to the case $\mu=0$. In the following we
assume that the connection used to define $\stark$ and $\starkB$ is
$G$-invariant; since the group action on $Q$ is assumed to be
proper such a connection always exists. Moreover , we assume that
$B \in Z^2_\dR (Q)^G[[\nu]]$. In case $\kappa \neq 0,1$ we finally
assume that $\sKov(\phi_g^* \alpha -\alpha)=0$. Note that one can
even achieve $\phi_g^*\alpha = \alpha$ using a different
$G$-invariant volume density for the definition of $\alpha$, namely
the Riemannian volume corresponding to an invariant Riemannian
metric on $Q$; if $\nabla$ is the pertinent Levi Civita connection
one even has $\alpha=0$. For reduction of $\stark$ we always use
the canonical $G$-equivariant classical momentum map with
$\JN{\xi}=\Pol{\xi_\Qind}$ as quantum momentum map $J^0$. By Remark
\ref{MomentValueRem} this causes no loss of generality. For
reduction of the products $\starkB$ we use a quantum momentum map
of the form $J^B(\xi)  = \Pol{\xi_\Qind} + \pi^* \jbold{\xi}$,
where $\d
\jbold{\xi} = i_{\xi_\Qind}B$ and $\phi_g^* \jbold{\xi} =
\jbold{\Ad(g^{-1})\xi}$.

In order to relate the reduced star products $\stark^{J^0,\mu}$ and
$\starkB^{J^B,\mu}$ for different momentum values with each another
we are going to construct local isomorphisms between them. For $\mu
\in {\mathfrak g^*}^G + \nu {\mathfrak g_c^*}^G[[\nu]]$ consider
the formal series of one-forms $\Gamma_\mu = \pair{\mu}{\gamma}\in
\Ginf{T^*Q}^G[[\nu]]$. Clearly, there is a uniquely defined
$b_\mu \in Z^2_\dR(\cc{Q})[[\nu]]$ such that $p^*b_\mu = - \d
\Gamma_\mu = - \pair{\mu}{\lambda}$. Now let $\{O_i\}_{i \in I}$
denote a good open cover of $Q$ such that $\{\cc{O_i}\}_{i\in I}$
with $\cc{O_i} = p (O_i)$ is a good open cover of $\cc{Q}$. Over every
$\cc{O_i}$ choose a local potential $a_\mu^i$ of $b_\mu$ which means
$b_\mu |_{\cc{O_i}} =
\d a_\mu^i$, and consider the formal series of locally defined
one-forms $A_\mu^i := \Gamma_\mu + p^* a_\mu^i\in
\Ginf{T^*O_i}[[\nu]]$. Then the $A_\mu^i$ turn out to be closed, hence the
operator $\mathcal A_{\mu,\kappa}^i$ defined by
Eq.\ (\ref{AkappaDefEq}) using $A_\mu^i$ instead of $A$ is a
local automorphism of $\stark$ and also of
$\starkB$ by Lemma \ref{starkBlinkrechtsMultLem} ii.).
Furthermore, $\mathcal A_{\mu,\kappa}^i$ is
$G$-invariant due to the invariance of the connection and the
invariance of $A_\mu^i$. By its form it is clear
that $\mathcal A_{\mu,\kappa}^i$  preserves $\mathcal
P(O_i)[[\nu]]$. Moreover, an easy computation shows that $ \mathcal
A_{\mu,\kappa}^i
\Pol{\xi_\Qind}|_{T^*O_i} = \Pol{\xi_\Qind}|_{T^*O_i} -
\pair{\mu}{\xi}$ and $\mathcal A_{\mu,\kappa}^i
(\Pol{\xi_\Qind} + \pi^*\jbold{\xi})|_{T^*O_i} = (\Pol{\xi_\Qind} +
\pi^*\jbold{\xi})|_{T^*O_i} - \pair{\mu}{\xi}$
for all $\xi \in \mathfrak g$. By these observations and
Proposition \ref{IsoAutoDerRedProp} ii.) we obtain:
\begin{LEMMA}\label{ImpulsVerglLem}
With notations from above, the mapping $\mathsf
A_{\mu,\kappa}^i: \hor (\mathcal P (O_i)^G)[[\nu]] \to  \hor
(\mathcal P (O_i)^G)[[\nu]]$ defined by
\begin{equation}
\mathsf A_{\mu,\kappa}^i F := \hor_{-\mu_0} \left(
\frac{\id}{\id - \nu \triangle_{\mu,\stark}}
\mathcal A_{\mu,\kappa}^i F \right), \quad
F \in \hor (\mathcal P(O_i)^G)[[\nu]],
\end{equation}
yields a local isomorphism from $(\hor (\mathcal P
(O_i)^G)[[\nu]],{\bullet_\kappa}^{J^0,0})$ to $(\hor (\mathcal P
(O_i)^G)[[\nu]],{\bullet_\kappa}^{J^0,\mu})$ which fulfills
\begin{equation}
\mathsf A_{\mu,\kappa}^i F = \mathcal A_{\mu,\kappa}^i F, \quad
 \text{ for all $F \in \hor (\mathcal P(O_i)^G)[[\nu]]$}.
\end{equation}
Replacing $\stark$ by $\starkB$, $J^0$ by $J^B$, and $-\mu_0$ by
$\Check j_0-\mu_0$ in the definition of $\mathsf A_{\mu,\kappa}^i$,
one obtains a local isomorphism from $(\hor (\mathcal P
(O_i)^G)[[\nu]], {\bullet_\kappa^B}^{J^B,0})$ to $(\hor (\mathcal P
(O_i)^G)[[\nu]],{\bullet_\kappa^B}^{J^B,\mu})$ which coincides with
$\mathsf A_{\mu,\kappa}^i$. Moreover, the corresponding operator
$A_{\mu,\kappa}^i:= l^{-1}\circ \mathsf A_{\mu,\kappa}^i\circ l$ is
a local isomorphism from $(\mathcal P
(\cc{O_i})[[\nu]],{\stark}^{J^0,0})$ to $(\mathcal P
(\cc{O_i})[[\nu]],{\stark}^{J^0,\mu})$ and from $(\mathcal P
(\cc{O_i})[[\nu]],{\starkB}^{J^B,0})$ to $(\mathcal P
(\cc{O_i})[[\nu]],{\starkB}^{J^B,\mu})$. This operator
extends uniquely to a local isomorphism, also denoted by $A_{\mu,\kappa}^i$,
from $(\Cinf{T^*\cc{O_i}}[[\nu]], {\stark}^{J^0,0})$ to
$(\Cinf{T^*\cc{O_i}}[[\nu]], {\stark}^{J^0,\mu})$ and from
$(\Cinf{T^*\cc{O_i}}[[\nu]],{\starkB}^{J^B,0})$ to
$(\Cinf{T^*\cc{O_i}}[[\nu]],{\starkB}^{J^B,\mu})$. Together with
the product ${\stark}^{J^0,0}$ resp.\ ${\starkB}^{J^B,0}$ the
local isomorphisms $\{A_{\mu,\kappa}^i\}_{i\in I}$ completely
determine the product ${\stark}^{J^0,\mu}$ resp.\
${\starkB}^{J^B,\mu}$.
\end{LEMMA}
\begin{PROOF}
First we observe that translations along the fibre and the
operators defined by Eq.\ (\ref{FOpDefEq}) preserve $\hor
(\PolFun^G)[[\nu]]$. But this implies that $\mathcal
A_{\mu,\kappa}^i$ preserves $\hor (\mathcal P(O_i)^G)[[\nu]]$,
therefore the contributions involving $\triangle_{\mu,\stark}$
vanish. Observing that $\hor_{-\mu_0}F' =F'$ for all $F'\in \hor
(\mathcal P(O_i)^G)[[\nu]]$ this implies $\mathsf A_{\mu,\kappa}^i
F =
\mathcal A_{\mu,\kappa}^i F$ for  $F\in \hor (\mathcal
P(O_i)^G)[[\nu]]$. A similar argument shows that the operator
analogous to $\mathsf A_{\mu,\kappa}^i$ which is obtained by
replacing $\stark$ by $\starkB$, $J^0$ by $J^B$, and $-\mu_0$ by
$\Check j_0-\mu_0$ also coincides with $\mathcal A_{\mu,\kappa}^i$.
Last, it remains to show that the product ${\stark}^{J^0,\mu}$
resp.\ ${\starkB}^{J^B,\mu}$ is completely determined by
${\stark}^{J^0,0}$ resp.\ ${\starkB}^{J^B,0}$ and the local
isomorphisms $\{\mathsf A_{\mu,\kappa}^i\}_{i\in I}$. In order to
check this one just has to observe that for $O_i\cap O_k
\neq \emptyset$ the compositions $(\mathsf A_{\mu,\kappa,}^k)^{-1}
\mathsf A_{\mu,\kappa}^i$ define automorphisms of $(\hor
(\mathcal P (O_i\cap O_k)^G)[[\nu]],{\bullet_\kappa}^{J^0,0})$
resp.\ $(\hor (\mathcal P (O_i\cap O_k)^G)[[\nu]],
{\bullet_\kappa^B}^{J^B,0} )$. Therefore, ${\stark}^{J^0,\mu}$
resp.\ ${\starkB}^{J^B,\mu}$ is globally defined by
$f{\stark}^{J^0,\mu} f'|_{T^*\cc{O_i}} = A_{\mu,\kappa}^i (
(A_{\mu,\kappa}^i)^{-1} f|_{T^*\cc{O_i}} {\stark}^{J^0,0}
(A_{\mu,\kappa}^i)^{-1} f'|_{T^*\cc{O_i}})$ resp.\ $f
{\starkB}^{J^B,\mu} f'|_{T^*\cc{O_i}}
= A_{\mu,\kappa}^i ( ( A_{\mu,\kappa}^i)^{-1} f|_{T^*\cc{O_i}}
{\starkB}^{J^B,0}( A_{\mu,\kappa}^i)^{-1}f'|_{T^*\cc{O_i}})$, where
$f,f' \in \Cinf{T^*\cc{Q}}[[\nu]]$.
\end{PROOF}

Clearly, the choice of $J^0$ and $J^B$ in the above lemma causes no
loss of generality since $\stark^{J', \mu}= \stark^{J^0, \mu -
\Tilde{\mu}}$, if $J'(\xi)= J^0(\xi)+ \pair{\Tilde{\mu}}{\xi}$,
and $\starkB^{J'',\mu}=\starkB^{J^B,\mu - \Tilde{\Tilde{\mu}}}$, if
$J''(\xi)= J^B(\xi)+ \pair{\Tilde{\Tilde{\mu}}}{\xi}$ (cf.\ Remark
\ref{MomentValueRem}). Therefore, the above result allows us to
compute the star products obtained from $\stark$ and $\starkB$ by
our reduction scheme in case we have at least determined one
reduced star product explicitly for a special choice of the quantum
momentum map and a special choice of the momentum value.

Now let us consider the standard ordered star product $\starNu$ and
the reduced product $\starNu^{J^0,0}$ more closely. Our goal is to
find out whether the reduced star product $\starNu^{J^0,0}$ is
again a standard ordered star product corresponding to a certain
torsion free connection $\cc{\nabla}$ on $\cc{Q}$.
\begin{LEMMA}\label{repNuRedLem}
\begin{enumerate}
\item
Let us assign to every $f\in \PolFunRed[[\nu]]$ a formal
series $\Tilde{\varrho}_0 (f)$
of differential operators on $\Cinf{\cc{Q}}$ by the following relation:
\begin{equation}\label{repNuRedDefEq}
p^* \Tilde{\varrho}_0(f) \chi =
\repNu{l(f)} p^* \chi, \quad \chi \in \Cinf{\cc{Q}}.
\end{equation}
Then
$\Tilde{\varrho}_0$ defines a representation of
$(\PolFunRed[[\nu]],\starNu^{J^0,0})$ on $\Cinf{\cc{Q}}[[\nu]]$ by
$\mathbb C[[\nu]]$-linear extension.
\item
If $\Tilde{\varrho}_0$ coincides with the standard ordered
representation $\cc{\varrho}_0$ with respect to some torsion free
connection on $\cc{Q}$, then $\starNu^{J^0,0}$ coincides with the
standard ordered star product $\cc{\star}_0$ corresponding to the
connection $\cc{\nabla}$ determined by
\begin{equation}\label{RedConnDefEq}
T p \, \nabla_{s^\hor} t^\hor = \cc{\nabla}_s t \circ p,
\quad s,t \in \Ginf{T\cc{Q}}.
\end{equation}
\item
If $\starNu^{J^0,0}$ coincides with the standard ordered star
product $\cc{\star}_0$ on $(T^*\cc{Q},\cc{\omega}_0)$ corresponding
to some torsion free connection on $\cc{Q}$, then
$\Tilde{\varrho}_0$ coincides with the standard ordered
representation $\cc{\varrho}_0$ with respect to the connection
$\cc{\nabla}$ determined by Eq.\ (\ref{RedConnDefEq}).
\end{enumerate}
\end{LEMMA}
\begin{PROOF}
By equivariance of $\varrho_0$ (see Eq.\ (\ref{repkapalmostequiEq}))
and by the fact that $l(f)\in
\hor(\PolFun^G)[[\nu]]$ the right-hand side of
Eq.\ (\ref{repNuRedDefEq}) is a $G$-invariant element of
$\Cinf{Q}[[\nu]]$, therefore $\Tilde{\varrho}_0(f)
\chi\in \Cinf{\cc{Q}}[[\nu]]$ is well-defined by this
equation, indeed. By the form of $\repNu{l(f)}$ one concludes that
$\Tilde{\varrho}_0(f)$ is a formal series of differential
operators, hence it can be extended to $\Cinf{\cc{Q}}[[\nu]]$ by
$\mathbb C[[\nu]]$-linearity. To prove that $\Tilde{\varrho}_0$ is
a representation of $(\PolFunRed[[\nu]],\starNu^{J^0,0})$ let us
first note that $\repNu{F}p^*\chi=0$ for all $F \in I^{\mathcal
P}_{0,\starNu}$, since $\varrho_0$ is a representation of $\starNu$
and since $\repNu{J^0(\xi)} p^*\chi
= - \nu \Lie_{\xi_\Qind} p^*\chi =0$ for all
$\xi\in \mathfrak g$.
Using this observation and the definition of $\starNu^{J^0,0}$
one immediately verifies that $\Tilde{\varrho}_0$ is a reperesentation.
This proves i.). Let us consider ii.).
In case $\Tilde{\varrho}_0$ coincides with the standard ordered
representation $\cc{\varrho}_0$ with respect to some torsion free
connection on $\cc{Q}$ the star product $\starNu^{J^0,0}$
coincides with the corresponding standard ordered star product
$\cc{\star}_0$, since the representation completely determines the
star product (cf.\ Section \ref{StarProdConsSubSec}). Moreover,
considering $\Tilde{\varrho}_0 (\PolRed{s\vee t})\chi =
\cc{\varrho}_0 (\PolRed{s \vee t})\chi$ for $s,t\in \Ginf{T\cc{Q}}$
one finds that the torsion free connection used to define
$\cc{\varrho}_0$ is uniquely determined and is given by
$\cc{\nabla}$ as in Eq.\ (\ref{RedConnDefEq}). For the proof
of iii.) assume that $\starNu^{J^0,0}=\cc{\star}_0$, where
$\cc{\star}_0$ is obtained from some torsion free connection on
$\cc{Q}$. Using the definition of $\Tilde{\varrho}_0$ it is immediate
to verify that $\Tilde{\varrho}_0(f)=
\cc{\varrho}_0(f)$ for all $f \in \mathcal P^0 (\cc{Q})
\oplus \mathcal P^1 (\cc{Q})$; note that the torsion
free connection used to define $\cc{\varrho}_0$ is of no
importance, hereby. Let us now assume that $\Tilde{\varrho}_0(f)=
\cc{\varrho}_0(f)$ for all $f \in \bigoplus_{k=0}^r
\mathcal P^k (\cc{Q})$ and consider $\PolRed{x_1\vee
\ldots \vee x_{r+1}}$ with $x_j\in \Ginf{T\cc{Q}}$. By Lemma
\ref{RedHomoCondLem} we know that $\starNu^{J^0,0}$ is a
homogeneous star product. Therefore
$\PolRed{x_1}\starNu^{J^0,0}\ldots
\starNu^{J^0,0} \PolRed{x_{r+1}} = \PolRed{x_1\vee
\ldots \vee x_{r+1}} + \sum_{l=1}^{r+1} \nu^l f_l =
\PolRed{x_1}\cc{\star}_0\ldots
\cc{\star}_0 \PolRed{x_{r+1}}$, where $f_l\in
\mathcal P^{r+1-l}(\cc{Q})$. Furthermore, we have
$\Tilde{\varrho}_0 (\PolRed{x_1}\starNu^{J^0,0}\ldots
\starNu^{J^0,0} \PolRed{x_{r+1}}) = \cc{\varrho}_0
(\PolRed{x_1}\cc{\star}_0\ldots \cc{\star}_0 \PolRed{x_{r+1}})$ by
the representation properties and the fact that $\Tilde{\varrho}_0
(\PolRed{x_j})= \cc{\varrho}_0 (\PolRed{x_j})$. Using the above
expression for $\PolRed{x_1}\starNu^{J^0,0}\ldots
\starNu^{J^0,0} \PolRed{x_{r+1}}$, this equation implies that
$\cc{\varrho}_0(\PolRed{x_1\vee \ldots \vee x_{r+1}})
=\Tilde{\varrho}_0(\PolRed{x_1\vee \ldots \vee x_{r+1}})$.
By an induction argument we then conclude that $\Tilde{\varrho}_0$ coincides
with $\cc{\varrho}_0$ on $\PolFunRed$, hence
$\Tilde{\varrho}_0=\cc{\varrho}_0$. Like for ii.) it then follows
that the connection used to define $\cc{\varrho}_0$ is given by
$\cc{\nabla}$ as in Eq.\ (\ref{RedConnDefEq}).
\end{PROOF}

Unfortunately, it is in general not true that
$\Tilde{\varrho}_0$ equals $\cc{\varrho}_0$ by the following
equality:
\[
p^*(\Tilde{\varrho}_0 -\cc{\varrho}_0)(\PolRed{x_1 \vee x_2\vee
x_3})\chi
=\frac{(-\nu)^3}{3} \!\sum_{\sigma \in S_3}\!\sum_{i=1}^{\dim{(G)}}
\Gamma_{e^i}(\nabla_{x^\hor_{\sigma(1)}}x^\hor_{\sigma(2)})
(\d p^*\chi)(\nabla_{x^\hor_{\sigma(3)}}{e_i}_\Qind).
\]
In fact, the analysis of the condition $\Tilde{\varrho}_0 =
\cc{\varrho}_0$ turns out to be rather involved, but at least we
can give two conditions which guarantee that this equality holds
true and then $\starNu^{J^0,0}= \cc{\star}_0$.
\begin{LEMMA}\label{NullRedCommNotBedLem}
If the connection $\nabla$ satisfies either
\begin{equation}\label{NablaCondEq1}
\nabla_{X}V \in \Ginf{VQ}^G \quad\textrm{for all $X\in
\Ginf{TQ}^G$, $V \in \Ginf{VQ}^G$},
\end{equation}
or
\begin{equation}\label{NablaCondEq2}
\nabla_{x^\hor}y^\hor \in \Ginf{HQ}^G\quad\textrm{for all }
x,y \in \Ginf{T\cc{Q}},
\end{equation}
the representation $\Tilde{\varrho}_0$ coincides with the standard
ordered representation $\cc{\varrho}_0$ with respect to the torsion
free connection $\cc{\nabla}$ defined by Eq.\ (\ref{RedConnDefEq}).
\end{LEMMA}
\begin{PROOF}
To show that $\Tilde{\varrho}_0$ actually coincides with the
standard ordered representation defined by means of $\cc{\nabla}$,
if one of the conditions (\ref{NablaCondEq1}),
(\ref{NablaCondEq2}) is satisfied, note first that it suffices
to prove the equality of the representations on elements
of $\PolFunRed$ of form $\PolRed{x_1\vee \ldots \vee x_k}$,
where $x_1,\ldots,x_k \in \Ginf{T\cc{Q}}$. Using the definition
of $\cc{\nabla}$ it is straightforward to check that
$(\nabla^k p^*\chi)(x_1^\hor,\ldots
,x_k^\hor)= p^*((\cc{\nabla}^k \chi)(x_1,\ldots,x_k))$ holds true for
for the $k$-fold covariant derivative, if
one of the assumptions (\ref{NablaCondEq1}),
(\ref{NablaCondEq2}) is satisfied. But from this
and the definition of $\varrho_0$ we immediately obtain
\[
\repNu{l (\PolRed{x_1\vee \ldots \vee
x_k})}p^*\chi = \repNu{\Pol{x_1^\hor\vee \ldots \vee
x_k^\hor}}p^*\chi = \frac{(-\nu)^k}{k!}p^*(i_s(x_1)\ldots
i_s(x_k)\cc{\sKov}^k\chi).
\]
This implies that $\Tilde{\varrho}_0(\PolRed{x_1\vee \ldots
\vee x_k})\chi= \frac{(-\nu)^k}{k!}i_s(x_1)\ldots
i_s(x_k)\cc{\sKov}^k\chi$. The last expression now is the standard
ordered representation of $\PolRed{x_1\vee \ldots \vee x_k}$ with
respect to the connection $\cc{\nabla}$. This proves the lemma. Finally,
let us note that if Eq.\ (\ref{NablaCondEq1}) is
satisfied, we moreover have  $(\nabla^k
p^*\chi)(Y_1,\ldots,Y_k)=0$ for $Y_1,\ldots,Y_k\in
\Ginf{TQ}^G$ in case at least one $Y_i$ is vertical.
\end{PROOF}

As a direct consequence of the lemma we obtain:
\begin{PROPOSITION}\label{starNuRedNullProp}
Let $\starNu$ be the standard ordered star product on
$(T^*Q,\omega_0)$ obtained from a $G$-invariant torsion free
connection $\nabla$ on $Q$ which satisfies one of the conditions
(\ref{NablaCondEq1}), (\ref{NablaCondEq2}). Then the reduced star
product $\starNu^{J^0,0}$ on $(T^*\cc{Q},\cc{\omega}_0)$ coincides
with the standard ordered star product $\cc{\star}_0$ corresponding
to the connection $\cc{\nabla}$ on $\cc{Q}$ defined by Eq.\
(\ref{RedConnDefEq}).
\end{PROPOSITION}
\begin{PROOF}
The claim follows immediately from Lemma \ref{repNuRedLem} and Lemma
\ref{NullRedCommNotBedLem}.
\end{PROOF}

In other words the preceding proposition just
means that, using appropriate connections, standard ordered quantization
commutes with reduction. Symbolically we have the following
commutative diagram:
\begin{equation}
\begin{CD}
(\Cinf{T^*Q},\{\,\,,\,\,\}_0)@>\mathcal
Q_0(\nabla)>>(\Cinf{T^*Q}[[\nu]],\starNu)\\ @VV \mathcal
R(J^0,0,\cdot)V @VV\mathcal R(J^0,0,\starNu)V\\
(\Cinf{T^*\cc{Q}},\{\,\,,\,\,\}_0)@>\mathcal
Q_0(\cc{\nabla})>>(\Cinf{T^*\cc{Q}}[[\nu]],\starNu^{J^0,0}=
\cc{\star}_0).
\end{CD}
\end{equation}
Hereby, $\mathcal Q_0(\nabla)$ resp.\ $\mathcal Q_0(\cc{\nabla})$
denotes the standard ordered quantization using the respective
connection and $\mathcal R(J^0,0,\cdot)$ resp.\ $\mathcal
R(J^0,0,\starNu)$ denotes classical resp.\ quantum reduction using
the indicated momentum map, momentum value and associative product.
Note that the condition expressed by  Eq.\ (\ref{NablaCondEq2}) is
rather restrictive, since it particularly implies, by using that
$\nabla$ is torsion free, that the horizontal distribution has to
be integrable, hence the principal connection corresponding to
$\gamma$ has to be flat. In contrast, the next lemma shows that
given a $G$-invariant torsion free connection $\nabla$, we can
always find another $G$-invariant torsion free connection
$\Hat{\nabla}$ which satisfies Eq.\ (\ref{NablaCondEq1}) and even
induces the same connection $\cc{\nabla}$ on $\cc{Q}$.
\begin{LEMMA}\label{nablaadaptLem1}
Let $\nabla$ denote a torsion free $G$-invariant connection on $Q$.
Define $\Hat{\nabla}$  by
\begin{equation}\label{nablaadaptEq1}
\begin{array}{rcl}
\Hat{\nabla}_H H' &:=& \nabla_H H',\\
\Hat{\nabla}_H V &:=& \nabla_H V - \mathrm{H} (\nabla_H V),
\end{array}
\qquad
\begin{array}{rcl}
\Hat{\nabla}_V H&:=& \nabla_V H- \mathrm{H} (\nabla_H V),\\
\Hat{\nabla}_V V'&:=& \nabla_V V'- \mathrm{H} (\nabla_V V'),
\end{array}
\end{equation}
where $H,H'\in \Ginf{HQ}$ and $V,V'\in \Ginf{VQ}$. Then $\Hat{\nabla}$ is
a torsion free $G$-invariant connection on $Q$ such that the induced
connection on $\cc{Q}$ coincides with the one induced by $\nabla$
and such that $\Hat{\nabla}_X V\in\Ginf{VQ}$ for all $X\in
\Ginf{TQ}, V\in \Ginf{VQ}$. Any other connection
$\Hat{\Hat{\nabla}}$ satisfying these conditions as well
is of form $\Hat{\Hat{\nabla}}_X Y =
\Hat{\nabla}_X Y + S(X,Y)$ with $S\in
\Ginf{\bigvee^2 T^*Q\otimes V Q}^G$.
\end{LEMMA}
\begin{PROOF}
Since $TQ=HQ\oplus VQ$, the vector field $\Hat{\nabla}_X
Y$ is well-defined for all $X,Y\in \Ginf{TQ}$ by Eq.\
(\ref{nablaadaptEq1}). The claim now follows by a
straightforward argument using the definition of $\Hat{\nabla}$.
\end{PROOF}

Analogously, if the principal connection
corresponding to $\gamma$ is flat, there exists a $G$-invariant torsion
free connection on $Q$ satisfying condition (\ref{NablaCondEq2}).
\begin{LEMMA}\label{nablaadaptLem2}
Let $\nabla$ denote a torsion free $G$-invariant connection on $Q$
and assume that the principal connection corresponding to $\gamma$
is flat. Define $\Check{\nabla}$ by
\begin{equation}
  \Check{\nabla}_H H':=\mathrm{H}(\nabla_H H'), \qquad
  \Check{\nabla}_Y Y':=\nabla_Y Y',
\end{equation}
where $H,H'\in \Ginf{HQ}$ and at least one of the
vector fields $Y,Y'\in \Ginf{TQ}$ is vertical. Then
$\Check{\nabla}$ is a torsion free $G$-invariant connection on $Q$
such that the induced connection on $\cc{Q}$ coincides with the one
induced by $\nabla$ and such that $\Check{\nabla}_H H'\in \Ginf{HQ}$
for all $H,H'\in \Ginf{HQ}$. Any other connection
$\Check{\Check{\nabla}}$ satisfying these conditions as well
is of form $\Check{\Check{\nabla}}_X Y =
\Check{\nabla}_X Y + S(X,Y)$ with $S\in
\Ginf{\bigvee^2 T^*Q\otimes  TQ}^G$ satisfying $S(H,H')=0$ for all
$H,H'\in \Ginf{HQ}$.
\end{LEMMA}
\begin{PROOF}
Again the proof is straightforward; the only crucial point to
observe is that the flatness of the principal connection implies
$\Check{\nabla}_H H' - \Check{\nabla}_{H'} H= [H,H']$.
Moreover, the condition that $\Check{\Check{\nabla}}$ induces the
same connection like $\Check{\nabla}$ and that it satisfies
$\Check{\Check{\nabla}}_{H}H'\in \Ginf{HQ}$ for all $H,H'\in
\Ginf{HQ}$ entails that $S(H,H')$ has to be both vertical and horizontal.
Thus $S(H,H')$ has to vanish for all $H,H'\in \Ginf{HQ}$.
\end{PROOF}

To conclude our study of the reduced star product
$\starNu^{J^0,0}$ let us mention that one can also use
the relation between the standard ordered representation
$\varrho_0$ and a symbolic calculus for pseudo-differential
operators on $\Cinf{Q}$ (see \cite[Sect.\ 6]{BorNeuWal99}
and \cite[Sect.\ 10]{BorNeuPflWal03}) to obtain the above `reduction
commutes with quantization' result. For details about
reduction of star products in this framework we
refer the interested reader to the thesis \cite{Kow01}.

Now we consider the reduction of the products $\stark$ with
$\kappa\neq 0$. These investigations will turn out to be slightly
more involved. First of all let us recall that the reduction scheme
introduced in Section \ref{RedSec} works for all the star products
$\stark$ under the general assumption that the connection $\nabla$
is $G$-invariant and that $\sKov(\phi_g^*\alpha -\alpha)=0$ holds
true additionally, if $\kappa\neq 0,1$. We will now show that
without additional assumptions the reduced star products
$\stark^{J^0,0}$ for $\kappa\neq 0$ are in general not even
equivalent to $\starNu^{J^0,0}$. This destroys the expectation that
$\stark^{J^0,0}$ could coincide with some naturally defined star
product $\cc{\star}_\kappa$ (cf.\ Section \ref{CompareSubSec} for a
concrete example in case $\kappa=1/2$).
\begin{LEMMA}
For all $\kappa \in [0,1]$ the $G$-invariant star product $\stark$
is $G$-equivalent to $\star_0^{B_\kappa}$ with $B_\kappa = \kappa
\nu \tr{R}$. Hence, $\stark$ is $G$-equivalent to $\starNu$, if and only
if $[B_\kappa]_\Gind=[0]_\Gind$. Consequently, the characteristic
class of $\stark^{J^0,0}$ is given by
\begin{equation}
c(\stark^{J^0,0}) = -\kappa[\cc{\pi}^*r] =-
\kappa[\cc{\pi}^*\tau_\lambda],
\end{equation}
where we have used the notation of Lemma \ref{CNullBestimmLem}.
\end{LEMMA}
\begin{PROOF}
Consider the operator defined by Eq.\
(\ref{AkappaDefEq}) using $A_\kappa=-\kappa\nu \alpha$ instead of
$A$ and denote it by $\mathcal A_0^{B_\kappa}$. Clearly,
$\mathcal A_0^{B_\kappa}$ is an equivalence transformation from
$\starNu$ to $\star_0^{B_\kappa}$, where $B_\kappa = -\kappa\nu \d
\alpha = \kappa\nu \tr{R}$. Therefore, $\mathcal A_0^{B_\kappa}
\Nk$ defines an equivalence transformation from $\stark$ to
$\star_0^{B_\kappa}$. Using Eq.\ (\ref{NkFactEq}), a
straightforward computation shows $\mathcal A_0^{B_\kappa}
\Nk = \exp \left(\Fop{\sum_{r=2}^\infty
\frac{1}{r!}(\kappa^r - \kappa) (-\nu)^r\sKov^{r-1}\alpha}\right)
\exp (-\kappa \nu \Delta_0)$. For $\kappa=0$ and $\kappa=1$ we have
$\mathcal A_0^{B_0} N_0= \id$ and $\mathcal A_0^{B_1} N_1 = \exp
(-\nu \Delta_0)$, which are both $G$-invariant operators. For
$\kappa \neq 0,1$ the operator $\mathcal A_0^{B_\kappa} \Nk$ is
$G$-invariant as well due to the additional condition
$\sKov(\phi_g^*\alpha -\alpha)=0$. This proves that $\stark$ is
$G$-equivalent to $\star_0^{B_\kappa}$. Thus, by
Proposition \ref{BBStrichGhomProp} the star product $\stark$ is
$G$-equivalent to $\starNu$, if and only if $[B_\kappa]_\Gind
=[0]_\Gind$. The result for the characteristic class of $\stark^{J_0,0}$
is an immediate consequence of Theorem \ref{GeneralCharKlassThm} and the
observation that the equations
$j_\kappa(\xi) = i^*\mathcal A_0^{B_\kappa} \Nk J^0(\xi) = -\kappa
\nu \mathsf{div}(\xi_\Qind)$ and $p^*b_\kappa = B_\kappa + \d
\Gamma_{\Check j_\kappa}$ with $b_\kappa=- \kappa \nu r$ hold true.
\end{PROOF}

Clearly, one could now change the momentum value to the one defined by
$\pair{\mu}{\xi} = - \kappa\nu\tr{\ad(\xi)}$. Thus one could achieve
$c(\stark^{J^0, \mu})=[0]$, but the result would not be a naturally
defined star product $\cc{\star}_\kappa$ for $\kappa\neq 0$.
Instead, we remain at the choice of $0$ momentum value and
merely adjust the parameters entering the construction of $\stark$
in order to obtain a star product in the
desired equivalence class. Henceforth, we thus assume that the
volume density $\upsilon$ is $G$-invariant. Consequently, the
one-form $\alpha$ defined in Eq.\ (\ref{alphaDefEq}) is also
$G$-invariant. Then $\Nk$ is a $G$-equivalence
from $\stark$ to $\starNu$, imlying that $\stark^{J^0,0}$ is
equivalent to $\starNu^{J'_\kappa,0}$, where $J'_\kappa(\xi)=
\Nk J^0(\xi)= J^0(\xi) - \kappa \nu \pi^*(\mathsf{div}(\xi_\Qind) +
\alpha(\xi_\Qind))$. But by $G$-invariance of
$\upsilon$ we have $\Lie_{\xi_\Qind}\upsilon =0$, which
entails  $\mathsf{div}(\xi_\Qind) + \alpha(\xi_\Qind)=0$  by
definition of $\alpha$. Hence
$\starNu^{J'_\kappa,0}=\starNu^{J^0,0}$. Analogously
to Lemma \ref{repNuRedLem} we now obtain:
\begin{LEMMA}\label{repkapRedLem}
Assume that $\kappa\neq 0$ and that
the volume density $\upsilon$ is $G$-invariant. Then assign to
every $f\in \PolFunRed[[\nu]]$ a formal series $\Tilde{\varrho}_\kappa(f)$
of differential operators on $\Cinf{\cc{Q}}$ by
\begin{equation}\label{repkapRedDefEq}
p^* \Tilde{\varrho}_\kappa(f) \chi =
\repkap{l(f)} p^* \chi\quad\textrm{for all }\chi
\in \Cinf{\cc{Q}} .
\end{equation}
Then $\Tilde{\varrho}_\kappa$ gives rise to a representation of
$(\PolFunRed[[\nu]],\stark^{J^0,0})$ on $\Cinf{\cc{Q}}[[\nu]]$ by
$\mathbb C[[\nu]]$-linear extension.
\end{LEMMA}
\begin{PROOF}
The proof of the claim is along the lines of the proof of Lemma
\ref{repNuRedLem} i.). The only additional relation which should be noted
for the proof of the representation property is the equality $\repkap{J^0(\xi)}
=\repNu{\Nk J^0(\xi)}= - \nu \Lie_{\xi_\Qind}$ which holds by
$G$-invariance of $\upsilon$.
\end{PROOF}

In order to interpret certain star products $\stark^{J^0,0}$ as
naturally defined star products $\cc{\star}_\kappa$ we need a
further condition which guarantees that the volume density
$\upsilon$ induces a volume density $\cc{\upsilon}$ on $\cc{Q}$.
The function $\upsilon ({e_1}_\Qind,\ldots, {e_{\dim{(G)}}}_\Qind,
x_1^\hor, \ldots, x_{n-\dim{(G)}}^\hor)$ turns out to be
$G$-invariant, if the group $G$ is unimodular. Unimodularity hereby
means that $|\det(\Ad(g))|=1$ for all $g\in G$, whence
$\tr{\ad(\xi)}=0$ for all $\xi\in \mathfrak g$. With the additional
assumption of $G$ to be unimodular we can define a volume density
$\cc{\upsilon}$ on $\cc{Q}$ by
\begin{equation}\label{RedDensDefEq}
p^*(\cc{\upsilon}(x_1,\ldots, x_{n-\dim{(G)}})) =\upsilon
({e_1}_\Qind,\ldots, {e_{\dim{(G)}}}_\Qind, x_1^\hor, \ldots,
x_{n-\dim{(G)}}^\hor).
\end{equation}
Evidently, the so-defined volume density $\cc{\upsilon}$ depends on
the chosen basis $\{e_i\}_{1\leq i\leq \dim{(G)}}$ of $\mathfrak
g$. But the choice of a different basis $\{e'_i\}_{1\leq i\leq
\dim{(G)}}$ yields a volume density $\cc{\upsilon}'=\mathsf a
\cc{\upsilon}$ with $\mathsf a \in \mathbb R^+$. Therefore, the
one-forms $\cc{\alpha},\cc{\alpha}'$ defined by
$\cc{\nabla}_x\cc{\upsilon} =\cc{\alpha}(x)\cc{\upsilon}$ and
$\cc{\nabla}_x\cc{\upsilon}' =\cc{\alpha}'(x)\cc{\upsilon}'$ with
$x\in \Ginf{T\cc{Q}}$ coincide. Hence, the corresponding
$\kappa$-ordered star products also coincide and are independent of
the above choice of a basis of $\mathfrak g$.

\begin{LEMMA}\label{repkRedLem}
Assume that $\kappa\neq 0$, that the volume density $\upsilon$ is
$G$-invariant and that $G$ is unimodular.
Then $\Tilde{\varrho}_\kappa$ coincides with the $\kappa$-ordered representation
$\cc{\varrho}_\kappa$ induced by the connection $\cc{\nabla}$
defined by Eq.\ (\ref{RedConnDefEq}) and the volume density
$\cc{\upsilon}$ defined by Eq.\ (\ref{RedDensDefEq}), if and only if
$\stark^{J^0,0}$ coincides with the $\kappa$-ordered star product
$\cc{\star}_\kappa$ corresponding to $\cc{\nabla}$ and the volume
density $\cc{\upsilon}$.
\end{LEMMA}
\begin{PROOF}
First we note that
$\Tilde{\varrho}_\kappa (f) = \cc{\varrho}_\kappa(f)$ for all $f
\in \mathcal P^0 (\cc{Q}) \oplus \mathcal P^1 (\cc{Q})$,
where $\cc{\varrho}_\kappa$ is defined using the connection
$\cc{\nabla}$ and the one-form $\cc{\alpha}\in \Ginf{T^*\cc{Q}}$
determined by $p^*\cc{\alpha} = \mathrm{H}(\alpha + W)$.
Here, the one-form $W$ is given by $W(X) =
\sum_{i=1}^{\dim{(G)}}\Gamma_{e_i} (\nabla_X {e_i}_\Qind)$
(cf.\ Lemma \ref{CNullBestimmLem}). In order to
interpret $\cc{\alpha}$ as the one-form defined by $\cc{\nabla}$
and some volume density $\cc{\upsilon}$ we must necessarily have
$\d \cc{\alpha} = - \tr{\cc{R}}$. Using the results and notation of
Lemma \ref{CNullBestimmLem} it is easy to compute that $\d
\cc{\alpha} = \tau_\lambda - \tr{\cc{R}}$, since $r = \d
(\cc{\alpha} - w)$. Actually, this relation suggests to assume that
the Lie group $G$ is unimodular, since then $\tau_\lambda=0$. In
this case, one finds $\cc{\nabla}_x
\cc{\upsilon} = \cc{\alpha}(x)\cc{\upsilon}$ for all $x\in
\Ginf{T\cc{Q}}$. Furthermore, it turns out that
$\tilde{\varrho}_\kappa (f)= \cc{\varrho}_\kappa (f)$ also holds
for all $f\in\mathcal P^2 (\cc{Q})$ without any further conditions
on the connection $\nabla$. With these observations, the proof
of the claim is completely analogous to the one of Lemma
\ref{repNuRedLem} ii.) and iii.).
\end{PROOF}

Like in the case $\kappa=0$, the operator $\Tilde{\varrho}_\kappa(f)$
coincides with $\cc{\varrho}_\kappa(f)$ without any further conditions on
$\nabla$, if $f \in \bigoplus_{k = 0}^2
\mathcal P^k(\cc{Q})$,
but for polynomials in the momenta of higher degree one does
not have $\Tilde{\varrho}_\kappa(f)=\cc{\varrho}_\kappa(f)$, in
general. Fortunately, the conditions imposed on $\nabla$
in case $\kappa=0$ turn out to be also guarantee
that the reduced star product $\stark^{J^0,0}$ coincides with
$\cc{\star}_\kappa$ defined by $\cc{\nabla}$ and $\cc{\upsilon}$
resp.\ $\cc{\alpha}$.
\begin{PROPOSITION}\label{starkRedNullProp}
For $\kappa\neq 0$ let $\stark$ be the $\kappa$-ordered star
product on $(T^*Q,\omega_0)$ obtained from a $G$-invariant volume
density $\upsilon$ and a $G$-invariant torsion free connection
$\nabla$ on $Q$ which satisfies one of the conditions
(\ref{NablaCondEq1}), (\ref{NablaCondEq2}). If $G$ is
unimodular, then the reduced star product $\stark^{J^0,0}$ on
$(T^*\cc{Q},\cc{\omega}_0)$ coincides with the $\kappa$-ordered
star product $\cc{\star}_\kappa$ which corresponds to the connection
$\cc{\nabla}$ on $\cc{Q}$ defined by Eq.\ (\ref{RedConnDefEq}) and to
the volume density $\cc{\upsilon}$ determined by Eq.\
(\ref{RedDensDefEq}).
\end{PROPOSITION}
\begin{PROOF}
By Lemma \ref{repkRedLem}, we only have to prove that each
of the conditions (\ref{NablaCondEq1}), (\ref{NablaCondEq2})
implies $\Tilde{\varrho}_\kappa = \cc{\varrho}_\kappa$. First,
we consider the case, where condition (\ref{NablaCondEq2})
is satisfied. We have to determine $\repkap{l(\PolRed{x_1\vee\ldots
\vee x_k})} p^* \chi = \repNu{\Nk \Pol{x^\hor_1\vee\ldots
\vee x^\hor_k}} p^* \chi$ for the proof of the claim.
Using the explicit form of the
operator $\Delta$ it is easy to find (cf.\ \cite[Eq.\
(1.3)]{BorNeuPflWal03}) that
\begin{eqnarray*}
\Delta \Pol{x^\hor_1\vee\ldots \vee x^\hor_k} &=& \Pol{
\sum_{l=1}^k x^\hor_1\vee\ldots\vee x^\hor_{l-1}
(\mathsf{div}(x_l^\hor)+ \alpha(x_l^\hor)) x_{l+1}^\hor\vee
\ldots\vee x^\hor_k}\\
& & +\Pol{
\sum_{\stackrel{l,j=1}{\scriptscriptstyle j\neq l}}^k
\nabla_{x_j^\hor} x_l^\hor \vee x^\hor_1\vee
\stackrel{l}{\Hat{\ldots}}\, \stackrel{j}{\Hat{\ldots}}
\vee x^\hor_k},
\end{eqnarray*}
where $\stackrel{j}{\Hat{\ldots}}$ denotes omission of the $j^{\rm
th}$ term. Using the definition of $\cc{\nabla}$ and $\cc{\alpha}$,
it is straightforward to show that $\mathsf{div}(x_l^\hor)+
\alpha(x_l^\hor) = p^*(\cc{\mathsf{div}}(x_l) +\cc{\alpha}(x_l))$.
Furthermore, the assumption about the connection $\nabla$ implies
$\nabla_{x_j^\hor} x_l^\hor = (\cc{\nabla}_{x_j}x_l)^\hor$. Putting
these formulas together we get $\Delta (l(\PolRed{x_1\vee\ldots
\vee x_k})) = l(\cc{\Delta}(\PolRed{x_1\vee\ldots \vee x_k}))$,
where the differential operator $\cc{\Delta}$ on $\Cinf{T^*\cc{Q}}$
is defined completely analogously to $\Delta$ using the connection
$\cc{\nabla}$ and the one-form $\cc{\alpha}$. Defining $\cc{\Nk}:=
\exp(-\kappa\nu \cc{\Delta})$ we obtain by induction that $\Nk l(f)
= l (\cc{\Nk} f)$ for all $f\in \PolFunRed[[\nu]]$. Using Lemma
\ref{NullRedCommNotBedLem} we then get
$p^*\Tilde{\varrho}_\kappa(f)\chi= \repNu{\Nk l(f)}p^*\chi=
\repNu{l(\cc{\Nk} f)}p^*\chi = p^*\cc{\varrho}_0
(\cc{\Nk}f)\chi=p^*\cc{\varrho}_\kappa(f)\chi$, which says that
$\Tilde{\varrho}_\kappa$ coincides with $\cc{\varrho}_\kappa$. Now
we consider the case, where condition
(\ref{NablaCondEq1}) is satisfied. As a first step we will to show
that $\hor (\Delta F) = l(\cc{\Delta}(l^{-1}(\hor(F))))$ for all
$F\in \PolFun^G$. To this end consider $F = \Pol{V_1\vee\ldots \vee
V_r\vee x_1^\hor\vee \ldots \vee x_k^\hor}$ with $V_1,\ldots,V_r\in
\Ginf{VQ}^G$. After application of $\Delta$ to $F$ several types of
terms appear according to the above formula. The terms involving
$\mathsf{div}(V_i) + \alpha(V_i)$ vanish, since $\mathsf{div}(V_i) +
\alpha(V_i)=0$ by  unimodularity of $G$ and the fact that
$\mathsf{div}(\xi_\Qind) + \alpha(\xi_\Qind)=0$. Moreover, as in
the first part of the proof we have $\mathsf{div}(x_j^\hor) +
\alpha(x_j^\hor)= p^*(\cc{\mathsf{div}}(x_j) + \cc{\alpha}(x_j))$.
From the second sum in the above formula four types of terms arise,
namely those involving $\nabla_{V_i}V_j$, $\nabla_{x_i^\hor}V_j$,
$\nabla_{V_j}x_i^\hor$, and $\nabla_{x_i^\hor} x_j^\hor$. By the
assumption on the connection it is evident that $\nabla_{V_i}V_j$,
$\nabla_{x_i^\hor}V_j$, and $\nabla_{V_j}x_i^\hor$ are all
vertical, hence these contributions vanish after projection to the
total horizontal part. Projecting to the horizontal part we may
also replace $\nabla_{x_i^\hor} x_j^\hor$ by $(\cc{\nabla}_{x_i}
x_j)^\hor$. Combining these results we get $\hor(\Delta F)=0$, if
$r\geq 1$, and $\hor(\Delta F) =
l(\cc{\Delta}(\PolRed{x_1\vee\ldots \vee x_k}))$, if $r=0$. But
this implies $\hor (\Delta F) = l(\cc{\Delta}(l^{-1}(\hor(F))))$
for all $F\in \PolFun^G$, since these are sums of terms of form
$\Pol{V_1\vee\ldots \vee V_r\vee x_1^\hor\vee \ldots \vee
x_k^\hor}$. By induction, this implies that $\hor(\Delta^k F) =
l(\cc{\Delta}^k(l^{-1}(\hor(F))))$ for all $k\in \mathbb N$.
Finally, one has to observe that $\repNu{\Pol{Y_1\vee\ldots\vee
Y_r}}p^*\chi=0$ for $Y_1,\ldots,Y_r\in \Ginf{TQ}^G$ in case at
least one of the $Y_i$ is vertical (cf.\ proof of Lemma
\ref{NullRedCommNotBedLem}). But then one finds
$p^*\Tilde{\varrho}_\kappa(f)\chi =
\sum_{k=0}^\infty \frac{(-\kappa\nu)^k}{k!}
\repNu{\Delta^k l(f)}p^*\chi =
\sum_{k=0}^\infty \frac{(-\kappa\nu)^k}{k!}
\repNu{\hor(\Delta^k l(f))}p^*\chi=
\sum_{k=0}^\infty \frac{(-\kappa\nu)^k}{k!}
\repNu{l(\cc{\Delta}^k f)}p^*\chi =
p^*(\cc{\varrho}_0(\cc{\Nk}f)\chi) = p^*\cc{\varrho}_\kappa(f)\chi$.
This shows that $\Tilde{\varrho}_\kappa =\cc{\varrho}_\kappa$.
\end{PROOF}

Symbolically, the above proposition can be
expressed by a commutative diagramm analogous to the one
for $\kappa=0$. More precisely, for $\kappa\in (0,1]$ and the appropriate
connections and volume densities the $\kappa$-ordered
quantization commutes with reduction, if $G$ is unimodular:
\begin{equation}
\begin{CD}
(\Cinf{T^*Q},\{\,\,,\,\,\}_0)@>\mathcal
Q_\kappa(\nabla,\upsilon)>>(\Cinf{T^*Q}[[\nu]],\stark)\\ @VV
\mathcal R(J^0,0,\cdot)V @VV\mathcal R(J^0,0,\stark)V\\
(\Cinf{T^*\cc{Q}},\{\,\,,\,\,\}_0)@>\mathcal
Q_\kappa(\cc{\nabla},\cc{\upsilon})>>(\Cinf{T^*\cc{Q}}[[\nu]],
\stark^{J^0,0}=
\cc{\star}_\kappa).
\end{CD}
\end{equation}
Hereby, $\mathcal Q_\kappa(\nabla,\upsilon)$ resp.\ $\mathcal
Q_\kappa(\cc{\nabla},\cc{\upsilon})$ denotes the $\kappa$-ordered
quantization using the respective connection and the respective
volume density, and $\mathcal R(J^0,0,\cdot)$ resp.\ $\mathcal
R(J^0,0,\stark)$ denotes classical resp.\ quantum reduction using
the indicated momentum map, momentum value and associative product.
\begin{REMARK}
At this point let us mention that the result of
Proposition \ref{starkRedNullProp} does not contradict \cite[Thm.\
4]{Emm93a}, where it was shown that in the special case of a certain
$G$-invariant Riemannian connection $\nabla^{\mathrm g}$ the
horizontal distribution has to be integrable in order to achieve
$\Tilde{\varrho}_{1/2} = \cc{\varrho}_{1/2}$. In contrast, our
result shows that the, in general, weaker condition
(\ref{NablaCondEq1}) suffices for the equality
$\Tilde{\varrho}_\kappa = \cc{\varrho}_\kappa$ to hold. The reason for this
peculiarity is that in the special case, where the $G$-invariant
Riemannian metric $\mathrm{g}$ has been chosen such that
$\mathrm{g}(V,H)=0$ for all $V\in \Ginf{VQ}$ and all
$H\in\Ginf{HQ}$, condition (\ref{NablaCondEq1}) implies that
(\ref{NablaCondEq2}) has to be satisfied, too. To check this, one
observes $\mathrm{g}(\nabla^{\mathrm g}_{x^\hor} V,
y^\hor) = -
\mathrm{g}(V,\nabla^{\mathrm g}_{x^\hor} y^\hor)$. But this term has to
vanish for all $V\in \Ginf{VQ}^G$ due to Eq.\ (\ref{NablaCondEq1})
and the orthogonality of vertical and horizontal vector fields.
This implies that $\nabla^{\mathrm g}_{x^\hor} y^\hor \in \Ginf{HQ}^G$
for all $x,y\in
\Ginf{T\cc{Q}}$, i.e.\ (\ref{NablaCondEq2}) holds true.
\end{REMARK}

After the investigation of the reduced star products obtained from
$\stark$, we finally consider the reduced star products
$\starkB^{J^B,0}$ obtained from $\starkB$. Having identified the
products $\stark^{J^0,0}$ we now want to relate these products to
$\starkB^{J^B,0}$ by means of local isomorphisms. In fact, the main
steps for the construction of these isomorphisms have already been
achieved in Section \ref{CharClassSec}, where we have considered
the case $\kappa =0$. Like in Eq.\ (\ref{BGammaDefbEq}) let us now
define $b_{B,\Check j}\in Z^2_\dR(\cc{Q})[[\nu]]$ by
$p^*b_{B,\Check j} = B + \d
\Gamma_{\Check j}$. Denote by $\{O_i\}_{i\in I}$ a good open
cover of $Q$ as in Section \ref{CharClassSec} and by
$\{\cc{O}_i\}_{i\in I}$ the corresponding good open cover of
$\cc{Q}$. Then we obtain formal local one-forms $A_{B,\Check j}^i
= p^* a_{B,\Check j}^i - \Gamma_{\Check j}$ on the $O_i$, where
$a_{B,\Check j}^i$ denotes a local potential of $b_{B,\Check j}$ on
$\cc{O}_i$. These formal local one-forms induce local isomorphisms $\mathcal
A^i_{B,\Check j,\kappa}
: (\Cinf{T^*O_i}[[\nu]],\stark)\to (\Cinf{T^*O_i}[[\nu]],\starkB)$
as defined by Eq.\ (\ref{AkappaDefEq}).
With these preparations we now get:
\begin{LEMMA}
With notations from above, the mapping $\mathsf A^i_{B,\Check
j,\kappa}: \hor(\mathcal P(O_i)^G)[[\nu]]\to
\hor(\mathcal P(O_i)^G)[[\nu]]$ defined by
\begin{equation}
  \mathsf A^i_{B,\Check j,\kappa} F:=
  \hor_{\Check{j_0}}\left(\frac{\id}{\id - \nu \triangle_{0,\starkB}}
  \mathcal A^i_{B,\Check j,\kappa}  F\right), \quad
  F\in \hor(\mathcal P(O_i)^G)[[\nu]],
\end{equation}
yields a local isomorphism from
$(\hor(\mathcal P(O_i)^G)[[\nu]],{\bullet_\kappa}^{J^0,0})$ to
$(\hor(\mathcal P(O_i)^G)[[\nu]],{\bullet^B_\kappa}^{J^B,0})$
which fulfills
\begin{equation}
 \mathsf A^i_{B,\Check j,\kappa} F = \mathcal A^i_{B,\Check
 j,\kappa} F\quad \text{ for all $F\in \hor(\mathcal P(O_i)^G)[[\nu]]$}.
\end{equation}
The induced mapping $A_{B,\Check j,\kappa}^i :=
l^{-1}\circ\mathsf A^i_{B,\Check j,\kappa} \circ l$ is a local
isomorphism from $(\mathcal P(\cc{O}_i)[[\nu]],{\stark}^{J^0,0})$
to $(\mathcal P (\cc{O}_i)[[\nu]],{\starkB}^{J^B,0})$. It
extends uniquely to a local isomorphism from
$(\Cinf{T^*\cc{O}_i}[[\nu]],{\stark}^{J^0,0})$ to
$(\Cinf{T^*\cc{O}_i}[[\nu]],{\starkB}^{J^B,0})$ which will be
denoted by $A^i_{B,\Check j,\kappa}$ as well.
\end{LEMMA}
\begin{PROOF}
The claim is evident by Proposition \ref{IsoAutoDerRedProp}
ii.), since $\mathcal A^i_{B,\Check j,\kappa}
\Pol{\xi_\Qind}|_{T^*O_i}= (\Pol{\xi_\Qind} +
\pi^*\jbold{\xi})|_{T^*O_i}$ and since
$\mathcal A^i_{B,\Check j,\kappa}$ preserves $\hor(\mathcal
P(O_i)^G)[[\nu]]$.
\end{PROOF}

In the following lemma we give sufficient conditions which allow for
a concrete computation of the above local isomorphisms
$\{A_{B,\Check j,\kappa}^i\}_{i\in I}$.
\begin{LEMMA}\label{AkapRedBestimmLem}
\begin{enumerate}
\item
Assume that $\nabla$ satisfies condition
(\ref{NablaCondEq1}) and that $B\in Z^2_\dR(Q)^G[[\nu]]$ is horizontal
so that we can choose $j =0$; this means in particular that
$\starkB$ is strongly $G$-invariant. Then the local isomorphism
$A_{B,0,\kappa}^i
: (\Cinf{T^*\cc{O}_i}[[\nu]],{\stark}^{J^0,0}) \to
(\Cinf{T^*\cc{O}_i}[[\nu]],{\starkB}^{J^0,0})$ is given by
\begin{equation}\label{AkapparedEq1}
A_{B,0,\kappa}^i = t^*_{-(a_{B,0}^i)_0}\exp\left(
-\Fop{\frac{\exp(\kappa\nu \cc{\sKov})- \exp((\kappa -1)\nu
\cc{\sKov})}{\nu\cc{\sKov}}a_{B,0}^i - (a_{B,0}^i)_0}
\right),
\end{equation}
where $\cc{\sKov}$ denotes the operator of symmetric covariant
derivation with respect to $\cc{\nabla}$.
\item
Assume that the $G$-invariant torsion free connection $\nabla$
satisfies condition (\ref{NablaCondEq2}). Then the local isomorphism
$A_{B,\Check j,\kappa}^i
: (\Cinf{T^*\cc{O}_i}[[\nu]],{\stark}^{J^0,0})\to
(\Cinf{T^*\cc{O}_i}[[\nu]],{\starkB}^{J^B,0})$ is given by
\begin{equation}\label{AkapparedEq2}
A_{B,\Check j,\kappa}^i = t^*_{-(a_{B,\Check j}^i)_0}\exp\left(
-\Fop{\frac{\exp(\kappa\nu \cc{\sKov})- \exp((\kappa -1)\nu
\cc{\sKov})}{\nu\cc{\sKov}}a_{B,\Check j}^i - (a_{B,\Check j}^i)_0}
\right).
\end{equation}
\end{enumerate}
\end{LEMMA}
\begin{PROOF}
In order to determine $A_{B,\Check j,\kappa}^i$ explicitly, it suffices
to compute $A_{B,\Check j,\kappa}^i
\PolRed{x_1\vee \ldots \vee x_k}$ for $x_1,\ldots,x_k \in
\Ginf{T\cc{Q}}$ and arbitrary $k\in \mathbb N\setminus \{0\}$,
since for $\chi\in \Cinf{\cc{Q}}$ the equation $A_{B,\Check
j,\kappa}^i\cc{\pi}^*\chi = \cc{\pi}^*\chi$ is evident. But from
the definition of $\mathsf F$ this means that we have to evaluate
terms of form $(\sKov^{k-1}(p^* a_{B,\Check j}^i -
\Gamma_{\Check j})) (x_1^\hor,\ldots,x_k^\hor)$. Let us first
consider the term involving $p^*a_{B,\Check j}^i$. To this end observe that
$p^*a_{B,\Check j}^i$ is a formal series of sums consisting of terms
of form $p^*(\chi \d \chi') = (p^*\chi) (\sKov p^*\chi')$ with $\chi,\chi'
\in \Cinf{\cc{O}_i}$. Hence it suffices to determine
$(\sKov^{k-1}((p^*\chi) (\sKov p^*\chi')))(
x_1^\hor,\ldots,x_k^\hor)$. Now recall that $\sKov$ is a derivation and
observe that the result is a sum of terms we have already
computed in the proof of Lemma \ref{NullRedCommNotBedLem}. Since
$\cc{\sKov}$ is a derivation as well, it is now straightforward
to compute that $(\sKov^{k-1}p^*(\chi \d
\chi'))(x_1^\hor,\ldots,x_k^\hor) =
p^*( (\cc{\sKov}^{k-1}(\chi \d \chi'))(x_1,\ldots,x_k))$, if one of
the conditions (\ref{NablaCondEq1}), (\ref{NablaCondEq2}) is
satisfied. If $B$ is horizontal and we may choose $j=0$, this
already implies Eq.\ (\ref{AkapparedEq1}). Now we turn to consider
the term $(\sKov^{k-1} \Gamma_{\Check j})
(x_1^\hor,\ldots,x_k^\hor)$ in case $j$ is arbitrary and $\nabla$
satisfies condition (\ref{NablaCondEq2}). But in this case a
straightforward induction argument shows that $(\sKov^{k-1}
\Gamma_{\Check j}) (x_1^\hor,\ldots,x_k^\hor)=0$, since
$\Gamma_{\Check j}$ vanishes on horizontal vector fields. Therefore,
$A^i_{B,\Check j,\kappa}$ then is given by Eq.\ (\ref{AkapparedEq2}).
\end{PROOF}

Using the result of the preceding lemma we can now easily identify
the star products $\starkB^{J^B,0}$ with naturally defined star
products on $(T^*\cc{Q},\omega_{b_0})$, if certain conditions are
satisfied. Moreover, observing that the local isomorphisms
$A^i_{\mu,\kappa}$ relating the reduced star products obtained from
different momentum values that were constructed in Lemma
\ref{ImpulsVerglLem} are of the same shape as the operators
$A^i_{B,\Check j,\kappa}$ we can -- slightly modifying the above
proof -- obtain the main result of this section.

\begin{THEOREM}
Let $\stark$ and $\starkB$ be the star products obtained from a
$G$-invariant torsion free connection $\nabla$ on $Q$ and $B\in
Z^2_\dR(Q)^G[[\nu]]$. For $\kappa \neq 0$, we moreover assume that
the volume density $\upsilon$ used to define $\alpha$ is
$G$-invariant and that the Lie group $G$ is unimodular.
\begin{enumerate}
\item
If $\nabla$ satisfies condition (\ref{NablaCondEq1}) and $B$ is
horizontal, then the star product $\starkB^{J^0,0}$ on
$(T^*\cc{Q},\omega_{(b_B)_0})$ coincides with the naturally defined
star product $\cc{\star}_\kappa^{b_B}$ obtained from the connection
$\cc{\nabla}$, the volume density $\cc{\upsilon}$, and $b_B\in
Z^2_\dR(\cc{Q})[[\nu]]$ defined by $p^*b_B = B$. In particular, the
star product $\stark^{J^0,0}$ on $(T^*\cc{Q},\cc{\omega}_0)$
coincides with $\cc{\star}_\kappa$.
\item
If $\nabla$ satisfies condition (\ref{NablaCondEq2}), then the star
product $\starkB^{J^B,\mu}$ on $(T^*\cc{Q},\omega_{b_0})$ coincides
with the naturally defined star product $\cc{\star}_\kappa^b$
obtained from the connection $\cc{\nabla}$, the volume density
$\cc{\upsilon}$, and $b\in Z^2_\dR(\cc{Q})[[\nu]]$ defined by $p^*b
= B+ \d \Gamma_{\Check j - \mu}$. In
particular, the star product $\stark^{J^0,\mu}$ on
$(T^*\cc{Q},\omega_{(b_\mu)_0})$ coincides with
$\cc{\star}_\kappa^{b_\mu}$, where $b_\mu
\in Z^2_\dR(\cc{Q})[[\nu]]$ is defined by $p^*b_\mu
=-\d \Gamma_\mu$.
\end{enumerate}
\end{THEOREM}
\begin{PROOF}
The claim about $\stark^{J^0,0}$ in i.) is just a restatement
of Proposition \ref{starNuRedNullProp} and Proposition
\ref{starkRedNullProp}. Moreover, Lemma \ref{AkapRedBestimmLem}
shows that the local isomorphisms from $\stark^{J^0,0}$ to
$\starkB^{J^0,0}$ are of the same form as the operators $\mathcal
A_\kappa$ of Eq.\ (\ref{AkappaDefEq}) using the
connection $\cc{\nabla}$ and local potentials of the formal series
$b_B=b_{B,0}$ of closed two-forms on $\cc{Q}$ defined by $p^*b_B =
B$. But this implies that $\starkB^{J^0,0}$ coincides with
$\cc{\star}_\kappa^{b_B}$ which is evidently a star product with
respect to $\omega_{(b_B)_0}$. For the proof of ii.) we observe
that the composition of the operators $A_{\mu,\kappa}^i$ and
$A_{B,\Check j,\kappa}^i$, which is a local isomorphism from
$\stark^{J^0,0}$ to $\starkB^{J^B,\mu}$, is given by
$t^*_{-(a_{B,\Check j}^i+ a_\mu^i)_0}\exp\left(\!
-\Fop{\frac{\exp(\kappa\nu \cc{\sKov})- \exp((\kappa -1)\nu
\cc{\sKov})}{\nu\cc{\sKov}}(a_{B,\Check j}^i + a_\mu^i) -
(a_{B,\Check j}^i+ a_\mu^i)_0}
\!\right)$ and that $\d(a_{B,\Check j}^i + a_\mu^i) =
b_{B,\Check j} + b_\mu = b$, where $p^*b = B + \d \Gamma_{\Check j
-\mu}$. But this implies that $\starkB^{J^B,\mu}$ coincides with
$\cc{\star}_\kappa^b$ using that
$\stark^{J^0,0}=\cc{\star}_\kappa$. Finally, the local isomorphism
$A_{\mu,\kappa}^i$ from $\stark^{J^0,0}$ to $\stark^{J^0,\mu}$ is
given by $t^*_{- (a_\mu^i)_0}\exp\left(
-\Fop{\frac{\exp(\kappa\nu \cc{\sKov})- \exp((\kappa -1)\nu
\cc{\sKov})}{\nu\cc{\sKov}} a_\mu^i -
(a_\mu^i)_0} \right)$. Hence $\stark^{J^0,\mu}$ equals
$\cc{\star}_\kappa^{b_\mu}$, where $p^*b_\mu = -\d \Gamma_\mu$.
\end{PROOF}
\subsection{Comparison to Existing Results}\label{CompareSubSec}
In this section we establish some relations between our
construction of reduced star products and known concepts of
reduction in deformation quantization. Additionally, we
consider the more specific example $T^*S^{n-1}$ which has been
discussed in the literature.
\subsubsection*{Remarks on Fedosov's Method}
In order to relate our results to the investigations of Fedosov
\cite{Fed98}, we will assume in this paragraph that the group
acting on $Q$ is compact. Then we consider the usual Fedosov star
product $\starF$ on $(T^*Q,\omega_{B_0})$ associated to a
$G$-invariant torsion free symplectic connection $\nabla^\QKotind$
on $T^*Q$ with $c(\starF)= \frac{1}{\nu}[\pi^*B_0]$. In order to
achieve that $\PolFun[[\nu]]$ is a $\starF$-subalgebra, we restrict
our considerations to connections whose Christoffel symbols in a
bundle chart are polynomials in the momenta (cf.\ \cite[Appx.\
A]{BorNeuWal98}). Under these general assumptions it is clear that
$\starF$ is $G$-invariant and even strongly $G$-invariant, i.e.\ we
can use $\JN{\xi} = \Pol{\xi_\Qind} +
\pi^*\jN{\xi}$ as $G$-equivariant quantum momentum map for
reduction. Moreover, Proposition \ref{GEquivaundBExiProp} tells us that
$\starF$ is $G$-equivalent to some $G$-invariant star product
$\starNuB$ with $B=B_0+ B_+$, where $B_+ \in \nu
Z^2_\dR(Q)^G[[\nu]]$. But then the characteristic
classes of $\starNuB$ and $\starF$ coincide, therefore $B_+$ has
to be exact. Since the group acting on $Q$ is compact, one can therefore
find a $G$-invariant formal potential $A_+$ for $B_+$.
Hence $\starNuB$ is $G$-equivalent to $\star_0^{B_0}$ by
Proposition \ref{BBStrichGhomProp}. Together with Theorem
\ref{GeneralCharKlassThm} this implies that the characteristic
class of $\starF^{J_0,0}$ is given by
\begin{equation}
c(\starF^{J_0,0}) = \frac{1}{\nu}[\cc{\pi}^*b_0] +
\frac{1}{\nu} [\cc{\pi}^*\Tilde{\mu}_\lambda],
\end{equation}
where $p^*b_0 = B_0 + \d \Gamma_{\Check j_0}$ and $\Tilde{\mu} \in
\nu \mathfrak g^*_c[[\nu]]$ is defined by $\pair{\Tilde{\mu}}{\xi}:
= i^*\mathcal T \JN{\xi}-\jN{\xi}$ with a $G$-equivalence
$\mathcal T$ from $\starF$ to $\star_0^{B_0}$. Here again,
$\Tilde{\mu}_\lambda\in \nu Z^2_\dR(\cc{Q})[[\nu]]$ is determined
by $p^*\Tilde{\mu}_\lambda = \pair{\Tilde{\mu}}{\lambda}$. Thus,
the reduced star product $\starF^{J_0,0}$ is in general not
equivalent to a canonical star product on the reduced phase space
with characteristic class $\frac{1}{\nu}[\cc{\pi}^*b_0]$. In view
of this fact -- similar to the situation for the star products
$\stark$ -- the `reduction commutes with quantization' theorem
proved by Fedosov appears to be a consequence of the appropriate
choice of the original star product on the large phase space rather
than a general principle in deformation quantization. Finally, let
us emphasize that the picture changes in case the Lie algebra
$\mathfrak g$ is semi-simple, since then there are no non-zero
elements of $\mathfrak g^*$ vanishing on $[\mathfrak g,\mathfrak
g]$ which implies that in this case $\Tilde{\mu}=0$. Another
example, where the term $[\cc{\pi}^*\Tilde{\mu}_\lambda]$ obviously
vanishes, is the case of a trivial principal $G$-bundle $p:Q\to
\cc{Q}$, where all the characteristic classes of the bundle are
zero.
\subsubsection*{The Star Product on $T^*S^{n-1}$ \`{a} la
Bayen et al.} In this paragraph we consider a very concrete
example of reduction. The reduced phase space is the
cotangent bundle of the $n-1$-sphere, which is obtained from
$T^*(\mathbb R^n\setminus \{0\})$ by classical Marsden-Weinstein reduction.
Using our reduction method for star products we will show that the
reduced star product obtained from the Weyl-Moyal star product
$\starWe$ on $T^*(\mathbb R^n\setminus \{0\})$  coincides with the
deformation quantization obtained by Bayen et al.\ in \cite{BayFla78}.

Let $G$ be the group $\mathbb R^+$ of positive real numbers with
the usual multiplication as composition, and consider the action of
$G$ on $\mathbb R^n\setminus \{0\}$ given by $\phi_g (x)= g x$ for
$g \in \mathbb R^+$, $x \in \mathbb R^n\setminus \{0\}$. Then the
quotient $(\mathbb R^n\setminus \{0\})/\mathbb R^+$ is isomorphic
to the sphere $S^{n-1} \subset \mathbb R^n\setminus\{0\}$. Then,
the corresponding projection $p : \mathbb R^n\setminus \{0\} \to
S^{n-1}$ is given by $p(x) =
\frac{1}{|x|}x$, where $|x|$ denotes the Euclidean length of
$x\in\mathbb R^n\setminus \{0\}$. For $\xi
\in \mathfrak g =
\mathbb R$, the generating vector field is explicitly given by
$\xi_{\mbox{\tiny $\mathbb R^n\setminus \{0\}$}}(x) = \xi x^i
\partial_{x^i}$. Hence, the (canonical) $G$-equivariant
classical momentum map reads $\JN{\xi}(q,p) =
\xi q^i p_i$. It is straightforward to check that $\gamma \in
\Ginf{T^*(\mathbb R^n\setminus \{0\})}$ with
$\gamma(x) = \frac{x^i}{|x|^2} \d x^i$ is a connection one-form
for the principal $\mathbb R^+$-bundle under consideration.
Clearly, $\gamma$ is closed, whence the connection is flat. In order
to compute the horizontal lift of a vector field $t\in
\Ginf{TS^{n-1}}$, it turns out to be convenient to describe $t$ by a
smooth mapping $\Tilde{t}: S^{n-1}\to \mathbb R^n$ which satisfies
$\big\langle\Tilde{t}\big(\frac{1}{|x|}x\big)\big| x \big\rangle=
0$ for all $x \in \mathbb R^n\setminus \{0\}$, where $\langle
{}\cdot{}|{}\cdot{}\rangle$ denotes the Euclidean inner product on
$\mathbb R^n$. With this mapping, the horizontal lift $t^\hor\in
\Ginf{T(\mathbb R^n\setminus \{0\})}$ is easily computed and given
by $t^\hor (x) = |x| \Tilde{t}\big(\frac{1}{|x|}x\big)$. But from
this the following relation is immediate:
\begin{equation}\label{SnLiftEq}
(l (\PolRed{t}))(q,p) = \Pol{t^\hor}(q,p)=
\Big\langle \Tilde{t}\Big(\frac{1}{|q|}q\Big) \Big|
|q|p- \frac{\langle q | p \rangle }{|q|}q \Big\rangle
= (\Pi^*\PolRed{t})(q,p),
\end{equation}
where $\Pi : T^*(\mathbb R^n\setminus\{0\})\to T^*S^{n-1}$,
$(q,p) \mapsto
\left(\frac{1}{|q|}q , |q|p- \frac{\langle q | p \rangle}{|q|}q\right)$
denotes the canonical projection onto $T^*S^{n-1}$. Note hereby,
that $T^*S^{n-1}$ is naturally embedded in $T^*(\mathbb
R^n\setminus\{0\})$ as the submanifold defined by the constraints
$|q|^2 =1$ and $\langle q | p \rangle =0$. Clearly, the formula in
Eq.\ (\ref{SnLiftEq}) also holds for all $t
\in \Ginf{\bigvee TS^{n-1}}$, since $l$ and $\Pi^*$ are homomorphisms
with respect to pointwise multiplication. The inverse $l^{-1}$
is just given by restriction of $F\in \hor (\mathcal P
(\mathbb R^n\setminus\{0\})^G)$ to $T^*S^{n-1}$, i.e.\ $l^{-1}=
I^*$, where $I: T^*S^{n-1}\to T^*(\mathbb R^n\setminus\{0\})$
denotes the embedding of $T^*S^{n-1}$ into $T^*(\mathbb
R^n\setminus\{0\})$.

After these preparations recall that the Weyl-Moyal star
product $\starWe$ on $T^*(\mathbb R^n\setminus \{0\})$ can be written as
\begin{equation}\label{WeylMoyalEq}
f \starWe f' = m \circ \exp
\left(\frac{\nu}{2}\left(\partial_{q^i}\otimes \partial_{p_i} -
\partial_{p_i} \otimes \partial_{q^i}\right)
\right)(f \otimes f'),\quad
\text{$f,f'\in \Cinf{T^*(\mathbb R^n\setminus\{0\})}[[\nu]]$},
\end{equation}
where $m$ is defined by $m (f\otimes f') := f f'$. Given $\mu
\in \mathbb R + \nu \mathbb C[[\nu]]$, we now want to compute the
reduced star product $\starWe^{J_0,\mu}$. According to
Theorem \ref{StarProdRedThm} and the above results for $l$ and
$l^{-1}$ it is given by
\begin{equation}\label{StarProdSnEq}
\PolRed{s} \starWe^{J_0,\mu} \PolRed{t} = \left.
\hor_{-\mu_0}\left(\frac{\id}{\id - \nu \triangle_{\mu,\starWe}}(
(\Pi^*\PolRed{s}) \starWe (\Pi^*\PolRed{t}))
\right) \right|_{T^*S^{n-1}}.
\end{equation}
The differential operator $\Fop{\gamma}$ with $(\Fop{\gamma}f)(q,p)
= \frac{q^i}{|q|^2} (\partial_{p_i}f)(q,p)$
turns out to be a derivation of $\starWe$, which is easily verified
using the explicit formula for $\starWe$ and the fact that $\gamma
$ is closed. Hence $(\Pi^*\PolRed{s})
\starWe (\Pi^*\PolRed{t})$ lies in $\hor(\mathcal P(\mathbb
R^n\setminus\{0\})^G)[[\nu]]$, and Eq.\ (\ref{StarProdSnEq}) simplifies to
\begin{equation}\label{BayFlaSnEq}
\PolRed{s} \starWe^{J_0,\mu} \PolRed{t} = \left.
(\Pi^*\PolRed{s}) \starWe (\Pi^*\PolRed{t}) \right|_{T^*S^{n-1}}=
I^*((\Pi^*\PolRed{s}) \starWe (\Pi^*\PolRed{t})),
\end{equation}
since $\triangle_{\mu,\starWe} F =0$ and $\hor_{-\mu_0}F
=F$ for every invariant totally horizontal polynomial function $F$.
Thus, the resulting star product $\starWe^{J_0}:=\starWe^{J_0,\mu}$
does not depend on the momentum value $\mu$ and therefore not on
the particular choice of the $G$-equivariant classical momentum map
(cf.\ Remark \ref{MomentValueRem}). This means in particular that
$\starWe^{J_0}$ is a star product with respect to the canonical
symplectic form on $T^* S^{n-1}$. Moreover, the above expression in
Eq.\ (\ref{BayFlaSnEq}) coincides with the star product constructed
by Bayen et al.\ in \cite{BayFla78}. There, using a slightly
different approach, $T^* S^{n-1}$ has been regarded as the quotient
$T^*(\mathbb R^n\setminus\{0\})/G'$, where the two-dimensional
non-Abelian group $G'=\mathbb R^+ \ltimes
\mathbb R$ is the semi-direct product of $(\mathbb R^+,\cdot)$ and
$(\mathbb R,+)$ over $\sigma:\mathbb R^+ \to \mathsf{Aut}(\mathbb
R)$, $\sigma(g)h = g^{-1} h$, which acts on $T^*(\mathbb
R^n\setminus\{0\})$ by $\Phi'_{(g,h)}(q,p)= (g q, g^{-1}p+ h q)$.
Finally, let us mention that our construction of reduced star
products can be applied to all star products $\stark$, but for
$\kappa\neq 1/2$ the concrete computation turns out to be much more
involved. This is caused by the fact that $\hor(\mathcal P(\mathbb
R^n\setminus\{0\})^G)[[\nu]]$ is not a $\stark$-subalgebra unless
$\kappa = 1/2$.  In particular, additional contributions arising
from $\frac{\id}{\id - \nu \triangle_{\mu,\stark}}$ have to be
taken into account, if $\kappa\neq 1/2$, which makes a concrete
computation of $\stark^{J_0,\mu}$ more difficult. To conclude the
discussion of this example note that the usual connection on
$\mathbb R^n\setminus\{0\}$ neither satisfies condition
(\ref{NablaCondEq1}) nor condition (\ref{NablaCondEq2}) for the
above choice of $\gamma$. Moreover, the canonical volume density on
$\mathbb R^n\setminus\{0\}$ from which $\starWe$ has been obtained,
is not invariant with respect to the considered action. Therefore,
it is not to be expected that $\starWe^{J_0}$ coincides with an
intrinsically constructed star product $\cc{\star}_{1/2}$ on
$T^*S^{n-1}$ using the induced connection $\cc{\nabla}$ and an
appropriate volume density on $S^{n-1}$.
\subsubsection*{The Counterexample of the BRST-Method}
In this paragraph we do not want to make an attempt to apply the
BRST-method in deformation quantization as developped in
\cite{BorHerWal00} to the example $(T^*Q,\omega_{B_0})$ in full
generality, since this would be far beyond the scope of the present
paper but might be an interesting topic for a future project.
Instead, we merely consider two concrete intimately related
examples also considered in \cite{BorHerWal00}, where in the first
one the BRST-method can be applied without any problems but turns
out to fail in the second example. In the first case we will show
in particular that our reduction procedure yields the same result
like the BRST-method. Moreover, we want to point out that due to
our slightly more restrictive definition of a quantum momentum map
the peculiarity occuring in the counterexample of the BRST-method
is avoided in our framework.

Let us consider $T^*(S^1\times S^1)\cong T^*S^1\times T^*S^1$ with
the canonical symplectic form. As symmetry group we choose $S^1$ which
acts on the base of the second factor by group multiplication.
Thus the reduced phase space is $T^*S^1$. Using the canonical flat
covariant derivative on $S^1$ we can equip each factor $T^*S^1$
with the $\kappa$-ordered star product $\stark$ yielding a
well-defined star product $\star_{\kappa_1,\kappa_2}$ on
$T^*(S^1\times S^1)$, where we may even use different ordering
parameters $\kappa_1,\kappa_2$ in each factor.
Denoting by $\Pi_1: T^*S^1\times T^*S^1 \to T^*S^1$ the projection
onto the first factor, one clearly has
$(\Pi_1^*f)\star_{\kappa_1,\kappa_2} (\Pi_1^*f')
= \Pi_1^*(f \star_{\kappa_1} f')$ for all $f,f'\in
\Cinf{T^*S^1}[[\nu]]$ by definition of the star product
$\star_{\kappa_1,\kappa_2}$. Using local coordinates $(\exp(\im
x^1),\exp(\im x^2))$ for $S^1\times S^1$ and the induced
coordinates $(\exp(\im q^1),\exp(\im q^2),p_1,p_2)$ of
$T^*(S^1\times S^1)$, it is evident that $\xi_{\mbox{\tiny
$S^1\times S^1$}}= - \im \xi \partial_{x^2}$ for $\xi \in \mathfrak
g = \im \mathbb R$ and that $\gamma=\im \d x^2$ defines a
connection one-form for the trivial $S^1$-bundle $p:S^1\times
S^1\to S^1$. Moreover, the canonical classical momentum map is
given by $\JN{\xi} = -\im \xi p_2$, which is easily verified to be
a quantum momentum map for $\star_{\kappa_1,\kappa_2}$. Now one can
compute the reduced star product
${\star_{\kappa_1,\kappa_2}}^{J_0,\mu}$ using that $l(f) = \Pi_1^*
f$ for all $f\in \mathcal P(S^1)[[\nu]]$, which is rather obvious
from the explicit choice of $\gamma$. Completely analogously to the
argumentation in the preceding paragraph we then find
\begin{equation}
f {\star_{\kappa_1,\kappa_2}}^{J_0,\mu} f' = f \star_{\kappa_1} f',
\quad \text{ $f,f'\in \mathcal P(S^1)[[\nu]]$}.
\end{equation}
Thus, the reduced star product just gives back the star
product we started from in the first factor of $T^*S^1\times
T^*S^1$, which is in perfect agreement with the results of
\cite{BorHerWal00}.

Following \cite{BorHerWal00}, we now consider the $S^1$-invariant
equivalence transformation $\mathcal T := \exp (-\nu p_1
\partial_{p_2})$ and the $S^1$-invariant star product
$\star'_{\kappa_1,\kappa_2}:= \mathcal T
\star_{\kappa_1,\kappa_2}$. According to our results of Proposition
\ref{IsoAutoDerRedProp}, it is immediately clear that the
corresponding reduced star product
${\star'_{\kappa_1,\kappa_2}}^{J',\mu}$ is equivalent to
${\star_{\kappa_1,\kappa_2}}^{J_0,\mu}$, where $J'(\xi) =
\mathcal T \JN{\xi} = \JN{\xi} + \im\nu\xi p_1$ is used as
$S^1$-equivariant quantum momentum map. Actually,
${\star'_{\kappa_1,\kappa_2}}^{J',\mu}$ even coincides with
${\star_{\kappa_1,\kappa_2}}^{J_0,\mu}$, since $\mathcal T^{-1}
\circ \Pi_1^*= \mathcal T \circ
\Pi_1^* = \Pi_1^*$, hence $(\Pi_1^* f)
\star'_{\kappa_1,\kappa_2} (\Pi_1^* f') = \Pi_1^*
(f \star_{\kappa_1} f')$ which implies
${\star'_{\kappa_1,\kappa_2}}^{J',\mu}=\star_{\kappa_1}$.

Now, the crucial point, which lets the BRST-method fail in this
case, is that $J_0$ is an allowed `quantum momentum map' for
BRST-quantization, since it satisfies
\[
\frac{1}{\nu}\left(\JN{\xi}\star'_{\kappa_1,\kappa_2}
\JN{\eta} - \JN{\eta}\star'_{\kappa_1,\kappa_2} \JN{\xi}\right)=
0 = \JN{[\xi,\eta]}
\]
by $\dim{(S^1)} =1$. But in contrast to the properties of $J'(\xi)$
we have $-\frac{1}{\nu}\ad_{\star'_{\kappa_1,\kappa_2}}(\JN{\xi})
\neq \Lie_{\xi_{\mbox{\tiny $T^*(S^1\times S^1)$}}}$, hence
the `quantum momentum map' $\JN{\xi}$ does not generate the
classical symmetry via the quasi-inner derivation with respect to
$\star'_{\kappa_1,\kappa_2}$, whereas it does so with respect to
$\star_{\kappa_1,\kappa_2}$. Thus, our slightly more restrictive
definition of a quantum momentum map -- which of course imposes
additional conditions to be satisfied a priori (cf.\ Section
\ref{InvQMMSec}) -- completely avoids the peculiarity appearing in
the BRST-method when using a non-appropriate `quantum momentum
map'.
\begin{appendix}
\section{Equivariance Properties of Certain Differential Operators}
Throughout this appendix, $\phi$ will always denote a
diffeomorphism of $Q$ and $\Phi=T^*(\phi^{-1})$ the diffeomorphism
lifted to $T^*Q$. Moreover, we
continue to use the notation as introduced in Section
\ref{PrelimSec}.
\begin{LEMMA}\label{pullbackConnLem}
Let $\nabla$ be a torsion free connection on $Q$. Then there is a
uniquely defined tensor field $S_\phi \in \Ginf{\bigvee^2
T^*Q\otimes TQ}$ such that $(\phi^*\nabla)_X Y -\nabla_X Y =
S_\phi(X,Y)$ for all $X,Y\in \Ginf{TQ}$. Furthermore, for all
$\beta\in\Ginf{\bigvee T^*Q}$ one has
\begin{equation}\label{pullbacksymKovEq1}
\phi^*\sKov (\phi^{-1})^* \beta = (\phi^*\sKov)\beta,
\end{equation}
where $\phi^*\sKov$ denotes the operator of symmetric covariant
derivation corresponding to the connection $\phi^*\nabla$ which is
explicitly given by
\begin{equation}\label{pullbacksymKovEq2}
\phi^*\sKov = \sKov - \d x^i\vee \d x^j \vee i_s(
S_\phi(\partial_{x^i},\partial_{x^j})).
\end{equation}
\end{LEMMA}
\begin{PROOF}The claim about the existence and uniqueness of
the tensor field $S_\phi$ and the fact that it is symmetric follows
immediately from the observation that $\phi^*\nabla$ is a torsion free
connection. In order to prove Eq.\ (\ref{pullbacksymKovEq1}), it is
enough to prove it for $\beta \in \Cinf{Q}$ and
$\beta\in \Ginf{T^*Q}$, since $\phi^*\sKov (\phi^{-1})^*$ and
$\phi^*\sKov$ are derivations with respect to $\vee$. It is straightforward
to verify the formula for these cases. Analogously, it suffices
to show that (\ref{pullbacksymKovEq2}) is satisfied on $\Cinf{Q}$ and
$\Ginf{T^*Q}$. An easy computation shows this as well.
\end{PROOF}

Now we briefly recall some well-known basic definitions
concerning horizontal and vertical lifts of vector fields and
one-forms on $Q$ to vector fields on $T^*Q$. Moreover, we study
their behaviour with respect to pull-back by $\Phi$.
\begin{DEFINITION}\label{HorVerDef}
{\bf (cf.\ \cite[Def.\ 2]{BorNeuWal98})} Let $\nabla$ be a
connection on $Q$. Consider the connection mapping $K :
T(T^*Q)\to T^*Q$ defined by
\begin{equation}\label{KonnektorDefEq}
K\left(\left.\frac{d}{dt}\right|_{t=0}\zeta(t)\right):=
\left.\nabla^c_{\partial_t}\zeta\right|_{t=0}
\end{equation}
for a curve $\zeta$ in $T^*Q$, where $\nabla^c$ denotes the connection
pulled-back along the footpoint-curve $c = \pi \circ
\zeta$. Then $(T\pi\times K): T(T^* Q) \to TQ \oplus
T^*Q$ is a fibrewise isomorphism. The horizontal
and vertical lifts with respect to $\nabla$ then are well-defined and unique
by the following. The section
$\Hor{\nabla}(X) \in \Ginf{T(T^*Q)}$ is called the horizontal lift
of $X\in \Ginf{TQ}$, if and only if
\begin{equation}\label{HorLifDefEq}
T\pi \,\Hor{\nabla}(X) = X \circ \pi \quad\textrm{ and } \quad
K(\Hor{\nabla}(X))=0,
\end{equation}
and $\Ver_\nabla(\beta) \in \Ginf{T(T^*Q)}$ is the
vertical lift of $\beta \in \Ginf{T^*Q}$, if and only if
\begin{equation}\label{VerLifDefEq}
K(\Ver_\nabla(\beta)) = \beta \circ \pi \quad\textrm{ and } \quad
T\pi \,\Ver_\nabla(\beta)=0.
\end{equation}
\end{DEFINITION}

Working in a local bundle chart one finds that
\begin{equation}\label{HorVerLocEq}
\Hor{\nabla}(X) = (\pi^*X^i) \partial_{q^i} + p_j \pi^*(X^k
\Gamma^j_{ki})\partial_{p_i} \quad\textrm{ and } \quad
\Ver_\nabla(\beta) = (\pi^*\beta_i) \partial_{p_i},
\end{equation}
where the $\Gamma^j_{ki}$ denote the Christoffel symbols of $\nabla$
and $X^i$ resp.\ $\beta_k$ the components of $X$ resp.\
$\beta$ in some chart of $Q$. In particular, it turns
out that the vertical lift does not depend on the connection,
henceforth we will simply denote it by $\Ver$.
\begin{LEMMA}\label{HorVerFopEquiLem}
Let $X\in \Ginf{TQ}$, $\beta\in \Ginf{T^*Q}$ and $\beta' \in
\Ginf{\bigvee T^*Q}$. Then one has the following equivariance
properties of $\Hor{\nabla}$, $\Ver$ and $\mathsf F$ with respect
to pull-back by $\Phi$:
\begin{equation}\label{HorVerEquiEq}
\Phi^* \Hor{\nabla}(X) = \Hor{\phi^*\nabla}(\phi^*X)\quad
\textrm{ and }\quad \Phi^* \Ver(\beta) = \Ver(\phi^*\beta).
\end{equation}
Moreover, one has for all $f \in \Cinf{T^*Q}$
\begin{equation}\label{FopEquiEq}
\Phi^* \Fop{\beta'}f = \Fop{\phi^*\beta'}\Phi^*f.
\end{equation}
\end{LEMMA}
\begin{PROOF}
The proof of Eqs.\ (\ref{HorVerEquiEq}) consists of a straightforward
computation using the definitions of horizontal and vertical
lifts. Observe that the second of these identities implies
$\Phi^*\Fop{\beta}f =
\Fop{\phi^*\beta}\Phi^*f$, since $\Fop{\beta}= \Lie_{\Ver(\beta)}$.
But since $\mathsf F$ is compatible with the
$\vee$-product and since $\Phi^*
\Fop{\chi}f= \Phi^*((\pi^*\chi) f) = (\pi^*\phi^*\chi) \Phi^*f =
\Fop{\phi^*\chi} \Phi^* f$ for all $\chi\in \Cinf{Q}$, this also
implies (\ref{FopEquiEq}).
\end{PROOF}

\end{appendix}
\begin{small}

\end{small}
\end{document}